\documentclass[12pt,a4paper,leqno,titlepage]{amsart}

\usepackage{amsmath,amscd,amsfonts,amsthm,amssymb,dsfont,bbm,manfnt,amsaddr}
\usepackage{hyperref,url,color,appendix,eucal}
\usepackage[usenames,dvipsnames]{xcolor}
\usepackage{multicol}
\usepackage[OT2,T1]{fontenc}
\usepackage[utf8]{inputenc}
\usepackage[all]{xy}
\usepackage{graphicx}

\def\got{\mathfrak}

\newenvironment{pf}
{\medskip\noindent {\it Proof.  }}
{\hfill\nobreak $\Box$ \par\bigbreak}

\newcommand{\isomo}{\overset{\sim}{\rightarrow}}

\newcommand{\GL}{\mathrm{GL}}
\newcommand{\SL}{\mathrm{SL}}
\newcommand{\ps}{\par \smallskip}

\newcommand{\Z}{\mathbb{Z}}

\newcommand{\Q}{\mathbb{Q}}

\newcommand{\R}{\mathbb{R}}    
\newcommand{\C}{\mathbb{C}}

\newcommand{\AAA}{\mathbb{A}}

\newtheorem{thm}[subsection]{Theorem}
\newtheorem*{thm*}{Theorem}
\newtheorem{lemme}[subsection]{Lemma}

\newtheorem*{lemme*}{Lemma}
\newtheorem{remark}[subsection]{Remark}

\newtheorem{cor}[subsection]{Corollary}

\newtheorem{prop}[subsection]{Proposition}

\newtheorem{example}[subsection]{Example}

\newtheorem{definition}[subsection]{Definition}

\newtheorem{scholie}[subsection]{Scholium}

\newtheorem{thmintro}{Theorem}

\makeatletter

\@addtoreset{equation}{section}
\makeatother

\hypersetup{colorlinks=true,linkcolor=black,citecolor=black,urlcolor=Periwinkle,pdfborderstyle={/S/U/W 1}}

\title[]{{\sc Subgroups of ${\rm Spin}(7)$ or ${\rm SO}(7)$ with each element conjugate to some element of ${\rm G}_2$ \\ and applications to automorphic forms}}

\author{Ga\"etan Chenevier}

\date\today

\thanks{\noindent Pendant la r\'edaction de ce travail, l'auteur a \'et\'e financ\'e par le C.N.R.S. et a re\c{c}u le soutien du projet ANR-14-CE25 (PerCoLaTor). The author thanks Wee Teck Gan, Philippe Gille and Gordan Savin for their remarks.}

\address{{\small Laboratoire de Math\'ematiques d'Orsay, Universit\'e Paris-Sud, CNRS, \\Universit\'e Paris-Saclay, 91405 Orsay, France,\\ 
{\rm \href{mailto:gaetan.chenevier@math.cnrs.fr}{gaetan.chenevier@math.cnrs.fr}}}}

\begin{document}

\begin{abstract} As is well-known, the compact groups ${\rm Spin}(7)$ and ${\rm SO}(7)$ both have a single conjugacy class of compact subgroups of exceptional type ${\bf G}_2$. We first show that if $\Gamma$ is a subgroup of ${\rm Spin}(7)$, and if each element of $\Gamma$ is conjugate to some element of ${\rm G}_2$, then $\Gamma$ itself is conjugate to a subgroup of ${\rm G}_2$. The analogous statement for ${\rm SO}(7)$ turns out be false, and our main result is a classification of all the exceptions. They are the following groups, embedded in each case in ${\rm SO}(7)$ in a very specific way: ${\rm GL}_2(\Z/3\Z)$, ${\rm SL}_ 2(\Z/3\Z)$, $\Z/4\Z \times \Z/2\Z$, as well 
as the nonabelian subgroups of ${\rm GO}_2(\C)$ with compact closure, similitude factors group $\{\pm1\}$, and which are not isomorphic to the dihedral group of order $8$. More generally, we consider the analogous problems in which the Euclidean space is replaced by a quadratic space of dimension $7$ over an arbitrary field.\ps
This type of questions naturally arises in some formulation of a converse statement of Langlands' global functoriality conjecture, to which the results above have thus some applications. Moreover, we give necessary and sufficient local conditions on a cuspidal algebraic regular automorphic representation of $\GL_7$ over a  totally real number field so that its associated $\ell$-adic Galois representations can be conjugate into ${\rm G}_2(\overline{\Q_\ell})$.

\end{abstract}

\maketitle

\section*{Introduction}

For $n\geq 1$, let ${\rm Spin}(n)$ denote the Spin group of the standard Euclidean space $\R^n$ and ${\rm SO}(n)$ its special orthogonal group. As is well-known there is a unique conjugacy class of compact subgroups of ${\rm Spin}(7)$ (resp. ${\rm SO}(7)$) which are connected, semisimple, and whose root system is of type ${\bf G}_2$. We shall call such a subgroup a {\it ${\rm G}_2$-subgroup}.  \ps\ps

\begin{thmintro} \label{thmintroa} Let $\Gamma \subset {\rm Spin}(7)$ be a subgroup such that each element of $\Gamma$ is contained in a ${\rm G}_2$-subgroup.  Then $\Gamma$ is contained in a ${\rm G}_2$-subgroup. 
\end{thmintro}

This perhaps surprising result has a very simple proof. Indeed, let $W$ be a spin representation of ${\rm Spin}(7)$, an $8$-dimensional real vector space, and let $E$ be its standard $7$-dimensional real representation, which factors through ${\rm SO}(7)$. If $H$ is a ${\rm G}_2$-subgroup of ${\rm Spin}(7)$, then we have an isomorphism of $\R[H]$-modules $$W \simeq 1 \oplus E.$$ As a consequence, if $\Gamma$ is as in the statement of Theorem \ref{thmintroa}, we must have the equality $\det (t - \gamma_{|W}) = (t-1) \det(t- \gamma_{|E})$ for all $\gamma \in \Gamma$. But this implies that the (necessarily semisimple) $\R[\Gamma]$-modules $W$ and $1 \oplus E$ are isomorphic. In particular, $\Gamma$ fixes a nonzero vector in $W$. We conclude as the ${\rm G}_2$-subgroups of ${\rm Spin}(7)$ are well-known to be exactly the stabilizers of the nonzero elements of $W$. \ps\ps

Theorem \ref{thmintroa} admits a generalization over an arbitrary field, that we prove in \S \,\ref{pfspin7}. Let $k$ be a field and $E$ a $7$-dimensional nondegenerate (or even ``regular'', see \S \ref{regqsp}) quadratic space over $k$. If $C$ is an octonion $k$-algebra whose quadratic subspace of pure octonions is isometric to $E$, then the automorphism group of $C$ naturally embeds both in ${\rm Spin}(E)$ and in ${\rm SO}(E)$ (\S \,\ref{g2emb}). We call such a subgroup of ${\rm Spin}(E)$ or ${\rm SO}(E)$ a {\it ${\rm G}_2$-subgroup}. These subgroups may not exist, but in any case they are all conjugate (by the action of ${\rm SO}(E)$, see Proposition\,\ref{carG2sbgp}). When $k=\R$ and $E$ is Euclidean they are exactly the subgroups introduced above. When $k$ is algebraically closed, they also do exist and coincide with the closed connected algebraic subgroups which are simple of type ${\bf G}_2$ (see Proposition\, \ref{algcarG2}). 

\begin{thmintro} \label{thmintrob} Let $k$ be a field, $E$ the quadratic space of pure octonions of some octonion $k$-algebra, $W$ a spinor module for $E$, and $\Gamma$ a subgroup of ${\rm Spin}(E)$. Assume that the $k[\Gamma]$-module $W$ is semisimple. Then the following assertions are equivalent: \ps 
\begin{itemize}
\item[(a)] each element of $\Gamma$ is contained in a ${\rm G}_2$-subgroup of ${\rm Spin}(E)$, \ps\ps
\item[(b)] for each $\gamma \in \Gamma$ we have $\det(1-\gamma_{|W})=0$, \ps \ps
\item[(c)] $\Gamma$ is contained in a ${\rm G}_2$-subgroup of ${\rm Spin}(E)$. \ps\ps
\end{itemize}
\end{thmintro}
\ps\ps

This theorem is Corollary \ref{corthmbsecondev} of Theorem \ref{thmbsecondev}, of which Theorem \ref{thmintroa} is the special case $k=\R$ (see Remark \ref{examplethmAR}). Note that Theorem \ref{thmintroa} also follows from the case $k=\C$ of Theorem \ref{thmintrob} by standard arguments from the theory of complexifications of compact Lie groups (see \S \ref{pfthmc}). \ps\ps
The naive analogue of Theorem \ref{thmintroa} to ${\rm SO}(7)$ turns out that to be false. Let us introduce three embeddings that will play some role in the correct statement. We fix a complex nondegenerate quadratic plane $P$, denote by ${\rm GO}(2)$ the unique maximal compact subgroup of the Lie group of orthogonal similitudes of $P$, and by ${\rm O}(2)^{\pm} \subset {\rm GO}(2)$ the subgroup\footnote{The group ${\rm O}(2)^{\pm}$ is isomorphic to the quotient of $\mu_4 \times {\rm O}(2)$ by the diagonal order $2$ subgroup, the similitude factor being the square of the factor $\mu_4$.}  of elements with similitude factor $\pm 1$. Recall that $E$ is the standard $7$-dimensional real representation of ${\rm SO}(7)$.  Then:\ps\ps

---  there is a group homomorphism $\alpha : \Z/4\Z \times \Z/2\Z \rightarrow {\rm SO}(7)$, unique up to conjugacy, such that the associated representation of $\Z/4\Z \times \Z/2\Z$ on $E \otimes \C$ is the direct sum of its $7$ nontrivial characters;\ps\ps

---  there is a morphism $\beta : {\rm GL}_2(\Z/3\Z) \rightarrow {\rm SO}(7)$, unique up to conjugacy, such that the representation of ${\rm GL}_2(\Z/3\Z)$ on $E \otimes \C$ is the direct sum of its $4$ nontrivial irreducible representations of dimension $\leq 2$;\ps\ps

---  there is a group morphism $\gamma : {\rm O}(2)^{\pm} \rightarrow {\rm SO}(7)$, unique up to conjugacy, such that the representation of ${\rm O}(2)^{\pm}$ on $E \otimes \C$ is isomorphic to the direct sum of $P$, $P^\ast$ and of the three order $2$ characters of ${\rm O}(2)^{\pm}$. \ps\ps

In each of these cases, if $\Gamma \subset  {\rm SO}(7)$ denotes the image of the given morphism, we prove in \S \ref{parexamples} that each element of $\Gamma$ belongs to a ${\rm G}_2$-subgroup of ${\rm SO}(7)$. Let ${\rm D}_8$ denote the dihedral group of order $8$. \ps\ps

\begin{thmintro} \label{thmintroc} Let $\Gamma$ be a subgroup of ${\rm SO}(7)$ such that each element of $\Gamma$ is contained in a ${\rm G}_2$-subgroup. Then exactly one of the following assertions holds: \ps 
\begin{itemize}
\item[(i)] $\Gamma$ is contained in a ${\rm G}_2$-subgroup of ${\rm SO}(7)$,\ps\ps
\item[(ii)] $\Gamma$ is conjugate in ${\rm SO}(7)$ to one of the groups 
$$\alpha(\Z/4\Z \times \Z/2\Z), \, \, \,\,\beta({\rm GL}_2(\Z/3\Z)),\, \, \,\,\beta({\rm SL}_2(\Z/3\Z))\, \, \, \, {\rm or}\,\,\, \, \gamma(S),$$
where $S \subset {\rm O}(2)^{\pm}$ is a nonabelian subgroup, nonisomorphic to ${\rm D}_8$, and whose similitude factors group is $\{\pm 1\}$. 
\end{itemize}
\end{thmintro}
\ps\ps

Let us discuss the main steps of the proof of this theorem, which turned out to be much harder than the one of Theorem \ref{thmintroa}. Let $\Gamma$ be a subgroup of ${\rm SO}(7)$. We view $E$ as a (semisimple) representation of $\Gamma$. We first show in \S \ref{basicid} that each element of $\Gamma$ is contained in a ${\rm G}_2$-subgroup of ${\rm SO}(7)$ if, and only if, we have a $\R[\Gamma]$-linear isomorphism\ps\ps
{\flushleft{${\rm (\ast)}$ \hspace{3 cm}\,\,\,\,  $\Lambda^3 \, E\, \simeq \,E \,\oplus\, {\rm Sym}^2\, E.$}}\ps\ps
It is an exercise to check that the images of $\alpha, \beta$ and $\gamma$ do satisfy ${\rm (\ast)}$. Note also that identity ${\rm (\ast)}$ implies that $\Gamma$ fixes a nonzero alternating trilinear form $f$ on $E$. In the special case $E \otimes \C$ is irreducible, the classification of trilinear forms given in \cite{cohen} shows that the stabilizer of $f$ is a ${\rm G}_2$-subgroup, and we are done. For reducible $E \otimes \C$, we proceed quite differently (and do not rely on any classification result).  \ps\ps

Our starting point is the following necessary and sufficient condition for $\Gamma$ to belong to a ${\rm G}_2$-subgroup of ${\rm SO}(7)$. Let $\widetilde{\Gamma} \subset {\rm Spin}(7)$ denote the inverse image of $\Gamma$ under the natural surjective morphism ${\rm Spin}(7) \rightarrow {\rm SO}(7)$. We show that $\Gamma$ is contained in a ${\rm G}_2$-subgroup of ${\rm SO}(7)$ if, and only if, the restriction to $\widetilde{\Gamma}$ of the spin representation $W$ of ${\rm Spin}(7)$ contains an order $2$ character (Corollary \ref{carbetaspin}). Our strategy to prove Theorem \ref{thmintroc} is then to analyze the simultaneous structure of the $\C[\widetilde{\Gamma}]$-modules $E \otimes \C$ and $W \otimes \C$, when the identity ${\rm (\ast)}$ holds. \ps\ps

We argue by descending induction on the maximal dimension of a $\Gamma$-stable isotropic subspace of the complex quadratic space $E$, an invariant that we call the {\it Witt index} of $\Gamma$ (see \S \ref{defwindex}). This invariant can only increase when we replace $\Gamma$ by a subgroup. Moreover, for a given Witt index we also argue according to the possible dimensions (and selfduality) of the irreducible summands of $E \otimes \C$. This forces us to study in details quite a number of special cases, which is somewhat unpleasant, but leads naturally to the discovery the counterexamples $\alpha,\beta$ and $\gamma$ above. Several exceptional properties of low dimensional spinor modules play a role along the proof.   \ps\ps

Our arguments apply more generally to the special orthogonal group of a $7$-dimensional regular quadratic space over an arbitrary algebraically closed field $k$, that we assume for this introduction to have characteristic $\neq 2$. We define in \S \ref{parexamples} a natural analogue ${\rm O}_2^{\pm}(k)$ of the group ${\rm O}(2)^{\pm}$, as well as analogues of the morphisms $\alpha, \beta$ and $\gamma$ above (the morphism $\beta$ is actually only defined if the characteristic of $k$ is $\neq 3$). Our main theorem is then the following (see Theorem \ref{mainthmb}).\ps\ps

\begin{thmintro} \label{thmintrod} Let $k$ be an algebraically closed field, $E$ a $7$-dimensional nondegenerate quadratic space over $k$ and $\Gamma \subset {\rm SO}(E)$ a subgroup. Assume that each element of $\Gamma$ has the same characteristic polynomial as some element of some ${\rm G}_2$-subgroup of ${\rm SO}(E)$, and that $\Gamma$ satisfies the ``semisimplicity'' assumption denoted by ${\rm (S)}$ in \S \ref{parass}. \ps
\indent Then the same conclusion as the one of Theorem \ref{thmintroc} holds, with ${\rm SO}(7)$ replaced by ${\rm SO}(E)$, and ${\rm O}(2)^{\pm }$ replaced by ${\rm O}_2^{\pm}(k)$.
\end{thmintro}

We refer to {\it loc. cit.} for a discussion of assumption ${\rm (S)}$. Let us simply say here that ${\rm (S)}$ holds if the $k[\Gamma]$-module $E$ is semisimple and if the characteristic of $k$ is either $0$ or $>13$. The special case $k=\C$ of Theorem \ref{thmintrod} is actually equivalent to Theorem \ref{thmintroc} (see \S \ref{pfthmc}). \ps\ps

Let us end this group theoretic discussion by raising some natural questions. Let $G$ be a group and $H$ a subgroup of $G$. Denote by\label{defcalp} $\mathcal{P}(G,H)$ the following property of $(G,H)$: for all subgroup $\Gamma$ of $G$, if each element of $\Gamma$ is conjugate to an element of $H$, then $\Gamma$ is conjugate to a subgroup of $H$. We have explained that $\mathcal{P}({\rm Spin}(7),{\rm G}_2)$ holds and that $\mathcal{P}({\rm SO}(7),{\rm G}_2)$ does not.\ps\ps

{\bf Questions:} {\it Can we classify the couples $(G,H)$, with $G$ a compact connected Lie group and $H$ a closed connected subgroup, such that $\mathcal{P}(G,H)$ holds ? Are there ``remarkable'' couples of finite groups $(G,H)$ such that $\mathcal{P}(G,H)$ holds ?}\ps\ps

For instance, one can show\footnote{Hints: consider subgroups $\Gamma$ of the form ${\rm SO}(a+b-c) \times {\rm SO}(c)$ with $c =0,1,3$ to show that $\mathcal{P}({\rm SO}(a+b),{\rm SO}(a) \times {\rm SO}(b))$ implies $b=1$ and $a \in \{1,3\}$. To check ${\mathcal P}({\rm SO}(4),{\rm SO}(3))$, observe that if $V$ is the restriction to ${\rm SO}(3)$ of the tautological $4$-dimensional real representation of ${\rm SO}(4)$, then we have an isomorphism $V \oplus V^\ast \simeq 1 \oplus 1 \oplus \Lambda^2\, V$ (a similar observation proves ${\mathcal P}({\rm SU}(4),{\rm SU}(3))$ as well).} that for integers $a \geq b \geq 1$, the property  $\mathcal{P}({\rm SO}(a+b),{\rm SO}(a) \times {\rm SO}(b))$ holds if, and only if, we have $b=1$ and $a \in \{1,3\}$. As another example, if $\got{S}_n$ denotes the symmetric group of $\{1,\dots,n\}$, then\footnote{ Hint: for $n=4,5$ consider for $\Gamma$ the subgroup of the alternating subgroup of $\got{S}_{n+1}$ preserving $\{n,n+1\}$ (we have $\Gamma \simeq \got{S}_{n-1}$). } $\mathcal{P}(\got{S}_{n+1},\got{S}_n)$ holds if, and only if, we have $n\leq 3$.\ps\ps
\medskip
\ps \ps
{\sc Applications to automorphic and Galois representations} 
\ps\ps

Our original motivation for studying property $\mathcal{P}$ above is that it naturally arises in some formulation of a converse statement of Langlands' functoriality conjecture. Before focusing on the specific case of our study, we first briefly give the general context, assuming some familiarity of the reader with Langlands philosophy \cite{langlandspb,euler} and automorphic forms \cite{corvallis}. \ps\ps

Let $F$ be a number field, $G$ (resp. $H$) a connected semisimple linear algebraic group defined and split over $F$, $\widehat{G}$ (resp. $\widehat{H}$) a complex semisimple algebraic group dual to $G$ (resp. $H$) in the sense of Langlands, $\mathcal{G}$ (resp. $\mathcal{H}$) a maximal compact subgroup of $\widehat{G}$ (resp. $\widehat{H}$), $\rho : \mathcal{H} \rightarrow \mathcal{G}$ a continuous homomorphism, and $\pi$ a cuspidal tempered automorphic representation of $G(\AAA_F)$. Assume that for all but finitely many places $v$ of $F$, the Satake parameter ${\rm c}(\pi_v)$ of $\pi_v$, viewed as a well-defined conjugacy class in $\mathcal{G}$, is the conjugacy class some element in $\rho(\mathcal{H})$. If $\mathcal{P}(\mathcal{G},\rho(\mathcal{H}))$ holds, then we may expect the existence of a cuspidal tempered automorphic representation $\pi'$ of $H(\AAA_F)$ such that for almost all finite places $v$ of $F$, the $\mathcal{G}$-conjugacy class of $\rho({\rm c}(\pi'_v))$ coincides with ${\rm c}(\pi_v)$. \ps\ps

Indeed, here is a heuristical argument using the hypothetical Langlands group \cite{Lgl,kott,arthurconjlan}. \label{lanintro}Fix a finite set $S$ of places of $F$ such that for $v \notin S$ then $v$ is finite and $\pi_v$ is unramified, and let $\mathcal{L}_{F,S}$ be the Langlands group of $F$ unramified outside $S$. Following Kottwitz, this is a compact topological group, equipped for each $v \notin S$ with a distinguished conjugacy class ${\rm Frob}_v$. Langlands associates to $\pi$ some continuous morphism $\phi : \mathcal{L}_{F,S} \rightarrow \mathcal{G}$, whose image has a finite centralizer in $\mathcal{G}$, such that\footnote{As pointed out by Langlands long ago, note that these properties of $\phi$ may however not determine it uniquely up to $\mathcal{G}$-conjugacy in general (see {\it e.g.} \cite{larsen}).} the $\mathcal{G}$-conjugacy class of $\phi({\rm Frob}_v)$ is  ${\rm c}(\pi_v)$ for each $v \notin S$. The expected density (even equidistribution!) of $\cup_{v \notin S} {\rm Frob}_v$ in $\mathcal{L}_{F,S}$ implies that any element in $\phi(\mathcal{L}_{F,S})$ is conjugate to some element in $\rho(\mathcal{H})$, so that we have $\phi = \rho \circ \phi'$ for some continuous morphism $\phi' : \mathcal{L}_{F,S} \rightarrow \mathcal{H}$, by property $\mathcal{P}(\mathcal{G},\rho(\mathcal{H}))$. In turn, $\phi'$ is associated to some cuspidal tempered automorphic representation $\pi'$ of $H(\AAA_F)$. \ps\ps

Of course, the properties $\mathcal{P}({\rm Spin}(7),{\rm G}_2)$ and $\mathcal{P}({\rm SO}(7),{\rm G}_2)$ studied before correspond to the special cases where $H$ is of type ${\bf G}_2$ and $G$ is either ${\rm PGSp}_6$ or ${\rm Sp}_6$. As a first application, Theorems \ref{thmintroa} \& \ref{thmintroc} thus provide interesting conjectures in this context. We postpone the discussion of the case $G={\rm PGSp}_6$, about which more can be said,  to the end of this introduction (see Theorem \ref{tig2pgsp6}). Note that although $\mathcal{P}({\rm SO}(7),{\rm G}_2)$ does not hold, the argument above still shows that the existence of $\pi'$ is expected to hold in the case $G={\rm Sp}_6$, because none of the exceptions in the statement of Theorem \ref{thmintroc} has a finite centralizer in $\mathcal{G}\simeq {\rm SO}(7)$ (their Witt index is nonzero).  \ps\ps

When one knows how to associate to $\pi$ a compatible system of $\ell$-adic Galois representations, which requires at least some assumptions on the Archimedean components of $\pi$, the absolute Galois group of $F$ may be used as a substitute of the hypothetical Langlands group. This allows us to prove the following theorem (Corollary \ref{corcplx}). We denote by ${\rm W}_{M}$ the Weil group of the local field $M$, ${\rm G}_2$ a fixed split semisimple group over $\Q$ of type ${\bf G}_2$, and if $k$ is an algebraically closed field of characteristic $0$ we fix an irreducible polynomial representation $\rho : {\rm G}_2(k) \rightarrow {\rm GL}_7(k)$. \ps

\begin{thmintro} \label{thmintrof} Let $F$ be a totally real number field and $\pi$ a cuspidal automorphic representation of ${\rm GL}_7(\AAA_F)$ such that $\pi_v$ is algebraic regular for each real place $v$ of $F$. The following properties are equivalent: \begin{itemize}\ps
\item[(i)] for all but finitely many places $v$ of $F$, the Satake parameter of $\pi_v$ is the conjugacy class of an element in $\rho({\rm G}_2(\C))$,\ps\ps
\item[(ii)] for any finite place $v$ of $F$, there exists a continuous morphism $\phi_v : {\rm W}_{F_v} \times {\rm SU}(2) \rightarrow {\rm G}_2(\C)$, unique up to ${\rm G}_2(\C)$-conjugacy, such that $\rho \circ \phi_v$ is isomorphic to the Weil-Deligne representation attached to  $\pi_v$ by the local Langlands correspondence \cite{harristaylor}.\ps
\end{itemize}
\end{thmintro}
\ps\ps
We refer to \S \ref{autgalrep} for the unexplained terms of this statement. A key ingredient in the proof of this theorem is the following result, which is perhaps the main application of this paper (Corollary \ref{cor1thmgalois}).
\ps\ps

 \begin{thmintro} \label{thmintroe} Let $F$ be a totally real number field and $\pi$ a cuspidal automorphic representation of $\GL_7(\AAA_F)$ satisfying assumption (i) of Theorem \ref{thmintrof}, and such that $\pi_v$ is algebraic regular for each real place $v$ of $F$.  Let $E$ be a coefficient number field for $\pi$, $\ell$ a prime and $\lambda$ a place of $E$ above $\ell$. Then there exists a continuous semisimple morphism $$\widetilde{r}_{\pi,\lambda}: {\rm Gal}(\overline{F}/F) \longrightarrow {\rm G}_2(\overline{E_\lambda}),$$ unique up to ${\rm G}_2(\overline{E_\lambda})$-conjugacy, satisfying the following property: for each finite place $v$ of $F$ which is prime to $\ell$, and such that $\pi_v$ is unramified, the morphism $\widetilde{r}_{\pi,\lambda}$ is unramified at $v$ and we have the relation  $\det (t - \rho(\widetilde{r}_{\pi,\lambda}({\rm Frob}_v))) = \det (t - {\rm c}(\pi_v))$, where ${\rm c}(\pi)$ is the Satake parameter of $\pi_v$ viewed as a semisimple conjugacy class in ${\rm GL}_7(E)$.
\end{thmintro}

\ps\ps The proof of this theorem uses the existence and properties of the compatible system of $7$-dimensional $\ell$-adic Galois representations associated to $\pi$ \cite{harristaylor,shin,chharris,bc} and Theorem \ref{thmintrod}. We show that we are not in the exceptional cases of Theorem \ref{thmintrod} using the knowledge of the Hodge-Tate numbers of those representations. The uniqueness assertion is a consequence of a result of Griess \cite{griess}. \ps\ps

As promised earlier, we now go back to property ${\mathcal P}({\rm G}_2,{\rm Spin}(7))$. After the first version of this paper appeared on the arXiv, some discussions with Gan and Savin led to the conclusion that enough is known to prove unconditionally the conjecture mentioned above concerning the Langlands functorial lifting between ${\rm G}_2$ and ${\rm PGSp}_6$, using works of Arthur \cite{arthurbook}, Ginzburg-Jiang \cite{gjiang} and Xu \cite{Xu}. We are grateful to them to let us include this result in this paper: see Theorem \ref{thmCGS}. Recall that we may take $\widehat{{\rm G_2}}={\rm G}_2(\C)$ and $\widehat{{\rm PGSp}}_6={\rm Spin}_7(\C)$.\ps\ps

\begin{thmintro} \label{tig2pgsp6} Let $F$ be a number field and $\pi$ a cuspidal automorphic representation of ${\rm PGSp}_6(\AAA_F)$. Assume that $\pi$ is tempered, or more generally, that $\pi$ is {\rm nearly generic} (see \S \ref{localcarg2sp6}). Then the following properties are equivalent: \begin{itemize}\ps\ps
\item[(a)]  for all but finitely many places $v$ of $F$, the Satake parameter ${\rm c}(\pi_v)$ is the conjugacy class of an element of a ${\rm G}_2$-subgroup of ${\rm Spin}_7(\C)$,\ps\ps
\item[(b)] there exists a cuspidal automorphic representation $\pi'$ of ${\rm G}_2(\AAA_F)$ such that for all but finitely many finite places $v$ of $F$, the image in ${\rm Spin}_7(\C)$ of the Satake parameter ${\rm c}(\pi'_v)$ is conjugate to ${\rm c}(\pi_v)$.\ps\ps
\end{itemize}
\end{thmintro}
\newpage 
\tableofcontents

\section{Preliminaries on quadratic spaces and spinors}

\subsection{Regular quadratic spaces} \label{regqsp}Let $k$ be a field. A quadratic space over $k$ is a finite dimensional $k$-vector space $V$ equipped with a quadratic map ${\rm q} : V \rightarrow k$ \cite[\S 3 {\rm no}.4]{bou}. By definition, the map $\beta_V(x,y):= {\rm q}(x+y)-{\rm q}(x)-{\rm q}(y)$, $V \times V \rightarrow k$, is a symmetric $k$-bilinear form, and we have ${\rm q}(\lambda v)=\lambda^2 {\rm q}(v)$ for all $v \in V$ and $\lambda \in k$. We refer to Bourbaki \cite[\S 3, \S 4]{bou} and Knus \cite[Ch. 1]{knuscampinas} for the basic notions concerning quadratic spaces (isometries, orthogonal sums, base change, etc...). \ps\ps

A quadratic space $V$ over $k$ will be called {\it nondegenerate} if the bilinear form $\beta_V$ is a perfect pairing. If the characteristic of $k$ is $2$, then $\beta_V$ is alternate, thus $\dim V$ has to be even if $V$ is nondegenerate. It will be convenient to say that a quadratic space $V$ over $k$ is {\it regular} if either $V$ is nondegenerate, or if we are in the following situation: $\dim V$ is odd, the characteristic of $k$ is $2$, the kernel of $\beta_V$ is one dimensional and ${\rm q}$ is not identically zero on it.\footnote{We warn the reader that there does not seem to be any standard terminology for these classical notions. Several authors, such as Knus, use the term nonsingular for nondegenerate. Moreover, $V$ is $\frac{1}{2}-$regular is the sense of \cite{knuscampinas, knusquadherm} if and only if $\dim V$ is odd and $V$ is regular in our sense.} \ps\ps

If $I$ is a finite dimensional $k$-vector space, we denote by $I^\ast$ its dual vector space and by ${\rm H}(I)$ the quadratic space $I \oplus I^\ast$ with ${\rm q}(x+\varphi)=\varphi(x)$ (``hyperbolic quadratic space over $I$'').  This is a nondegenerate quadratic space. \ps\ps

The isometry group of a quadratic space $V$, denoted ${\rm O}(V)$, is the subgroup of elements $g \in {\rm GL}(V)$ such that ${\rm q} \circ g = g$. When $V$ is regular, there is a well-defined subgroup ${\rm SO}(V) \subset {\rm O}(V)$ of proper isometries: when $2 \in k^\times$ or $\dim V$ is odd, we set ${\rm SO}(V) = {\rm O}(V) \cap {\rm SL}(V)$, and in the remaining case ($\dim V$ even and $2 \not\in k^\times$) ${\rm SO}(V)$ is an index $2$ subgroup of ${\rm O}(V)$ defined as the kernel of the Dickson invariant of $V$ \cite[Ch. 6 p.59]{knuscampinas}.  We also denote by ${\rm GO}(V) \subset {\rm GL}(V)$ the subgroup of orthogonal similitudes of the quadratic space $V$, and when $V$ is regular, by ${\rm GSO}(V)$ its subgroup of proper similitudes (see e.g. the end of \S II.1 in \cite{cl}).\ps\ps
\subsection{Clifford algebras, Spin groups, and their relatives}\label{parclag}

Let $V$ be a quadratic space over $k$. Recall that the Clifford algebra ${\rm C}(V)$ of $V$ is a $k$-superalgebra \cite[\S 9]{bou} \cite[Ch. 4]{knuscampinas}  \cite{deligne}. It is equipped with a canonical injective morphism $k \oplus V \rightarrow {\rm C}(V)$, that we shall always view as an inclusion. The {\it Clifford group} of $V$, denoted $\Gamma(V)$, is the supergroup\footnote{A supergroup is a group $G$ equipped with a group homomorphism ${\rm p} : G \rightarrow \Z/2\Z$.} of homogeneous elements in ${\rm C}(V)^\times$ normalizing $V$. The normal subgroup of even elements in $\Gamma(V)$ is the {\it general spin group} of $V$ and denoted ${\rm GSpin}(V)$. The superaction of $\Gamma(V)$ on $V$, given by $(\gamma,v) \mapsto (-1)^{{\rm p}(\gamma)} \gamma v \gamma^{-1}$, defines a group homomorphism $$\pi_V : \Gamma(V) \longrightarrow {\rm O}(V),$$ \ps\ps
whose kernel is by definition the group of invertible homogeneous elements in the supercenter ${\rm Z}(V)$ of ${\rm C}(V)$. Following Chevalley, recall that ${\rm C}(V)$ is equipped with a unique anti-involution $x \mapsto x^{\rm t}$ fixing $V$ pointwise (hence preserving parity). 
It thus preserves $\Gamma(V)$ as well. If $\gamma \in \Gamma(V)$, the element $\nu(\gamma):=\gamma \gamma^{\rm t}$ belongs to ${\rm Z}(V)_0^\times$, which gives rise to a group homomorphism $$\nu_V : \Gamma(V) \rightarrow  {\rm Z}(V)_0^\times.$$
The kernel of $\nu_V$ is the so-called {\it pin} supergroup of $V$ and denoted ${\rm Pin}(V)$. The subgroup of even elements in ${\rm Pin}(V)$ is the {\it spin group} of $V$, denoted ${\rm Spin}(V)$.  We have $\nu_V (\lambda) = \lambda^2$ for all $\lambda \in k^\times$.\ps \ps

Assume from now on that $V$ is regular. If $\dim V$ is even (resp. odd) then ${\rm C}(V)$ (resp. ${\rm C}(V)_0$) is a central simple $k$-algebra \cite[\S 4]{bou} \cite[Ch. 4 \, Thm. 8]{knuscampinas}; in particular we have ${\rm Z}(V)_0=k$ in both cases. Applying the Skolem-Noether theorem, we can show that the homomorphism $\pi_V$ is surjective \cite[Prop. 6 Ch. 6]{knuscampinas} and induces 
an exact sequence
\begin{equation}\label{sexgspin} 1 \longrightarrow k^\times \longrightarrow {\rm GSpin}(V) \overset{\pi_V}{\longrightarrow} {\rm SO}(V) \longrightarrow 1.\end{equation}
Moreover, $\nu_V : {\rm GSpin}(V) \longrightarrow k^\times$ induces a homomorphism $\overline{\nu_V} : {\rm SO}(V) \rightarrow k^\times/{k^\times}^2$ called the {\it spinor norm}. We have thus an exact sequence 
\begin{equation}\label{sexspin} 1 \longrightarrow \mu_2(k) \longrightarrow {\rm Spin}(V) \overset{\pi_V}{\longrightarrow} {\rm SO}(V) \overset{\overline{\nu_V}}{\longrightarrow} k^\times/{k^\times}^2, \end{equation}
where $\mu_2(k) \subset k^\times$ is the subgroup of square roots of $1$. \ps\ps

Let $U,V$ be quadratic spaces over $k$ with $U$ nondegenerate and $V$ regular. Then $W=U \bot V$ is regular and the natural $k$-superalgebra isomorphism ${\rm C}(U) \otimes^{\rm gr} {\rm C}(V) \isomo {\rm C}(U \bot V)$ induces a group homomorphism 
\begin{equation}\label{rhoUV} \rho_{U; V} : {\rm GSpin}(U) \times {\rm GSpin}(V) \rightarrow {\rm GSpin}(U \bot V) \end{equation}
which satisfies $\nu_W \circ \rho_{U; V}(g,g') = \nu_{U}(g)\nu_{V}(g')$, $\pi_W \circ \rho_{U; V} = \pi_U \oplus \pi_{V}$ \cite[Ch. IV \S (6.5)]{knusquadherm}. By \eqref{sexgspin}, the image of $\pi_W \circ \rho_{U; V}$ is thus ${\rm SO}(U) \times {\rm SO}(V)$, and the kernel of $\rho_{U;V}$ is the subgroup of $k^\times \times k^\times$ whose elements $(\lambda,\lambda')$ satisfy $\lambda\lambda'=1$. \ps\ps


\begin{remark}\label{gpscheme} {\rm Most of the constructions of \S \ref{regqsp} and \S \ref{parclag} can be adapted in order to make sense over an arbitrary commutative ring $k$: see e.g. \cite{knusquadherm}, \cite[Ch. I \S 1]{cl} and \cite{cf}. That being done, each group of the form ${\rm G}(V)$ associated above to a regular quadratic space $V$ over $k$ appears as the group of $k$-points of a natural corresponding affine group scheme of finite type over $k$. We shall not need this point of view in the sequel, but in a few places, and when the field $k$ is algebraically closed, it will be convenient to see ${\rm G}(V)$ as a (reduced) linear algebraic group over $k$  in the classical sense \cite{humphreys1}; it is indeed immediate from the definitions above that each ${\rm G}(V)$ does have a linear algebraic group structure. As is well-known, the algebraic groups ${\rm SO}(V)$, ${\rm Spin}(V)$ and ${\rm GSpin}(V)$ are connected and reductive if $\dim V > 1$ (see e.g. \cite{cf}). 
}\end{remark}

\ps\ps
{\sc Notations:} Let $n\geq 0$ be an integer. We respectively denote by ${\rm O}_n(k)$, ${\rm SO}_n(k)$, ${\rm GO}_n(k)$, ${\rm GSO}_n(k)$, $\Gamma_n(k)$, ${\rm Spin}_n(k)$ and ${\rm GSpin}_n(k)$ the groups ${\rm O}(V)$, ${\rm SO}(V)$, ${\rm GO}(V)$, ${\rm GSO}(V)$, $\Gamma(V)$, ${\rm Spin}(V)$ and ${\rm GSpin}(V)$ with 
$$ V = \left\{  \begin{array}{ll} {\rm H}(k^{r}) & \,\,{\rm \,if}\,\,  n=2r, \\  {\rm H}(k^r) \bot k & \,\,{\rm \, if}\,\, n=2r+1.\end{array} \right. $$
 In this latter case, $k$ is viewed as a quadratic space with ${\rm q}(x)=x^2$. In both cases, $V$ is regular. When $k$ is algebraically closed, it is the unique regular quadratic space of dimension $n$ over $k$ up to isometry. Moreover, when $k=\R$, we shall also denote respectively by ${\rm O}(n)$, ${\rm SO}(n)$ and ${\rm Spin}(n)$ the corresponding groups with $V$ the standard Euclidean space $\R^n$ of dimension $n$. They are compact Lie groups in a natural way.

\subsection{Spinor modules}\label{spinormodules}

Let $V$ be a regular quadratic space over $k$. We say that the even Clifford algebra of $V$ is {\it trivial} if either $\dim V = 2r+1$ is odd and we have a $k$-algebra isomorphism ${\rm C}(V)_0 \simeq {\rm M}_{2^r}(k)$, or $\dim V = 2r$ is even and we have a $k$-algebra isomorphism ${\rm C}(V)_0 \simeq {\rm M}_{2^{r-1}}(k) \times {\rm M}_{2^{r-1}}(k)$.
Assuming that $V$ has a trivial even Clifford algebra, we can then define a spinor module for $V$ as follows. Again, there are two cases:\ps\ps
 
Case (a). If $\dim V=2r$ is even, then we have a graded algebra isomorphism ${\rm C}(V) \simeq {\rm M}_{2^{r-1}|2^{r-1}}(k)$. A spinor module for $V$ is a supermodule $W$ for the $k$-superalgebra ${\rm C}(V)$ such that the associated graded morphism ${\rm C}(V) \rightarrow {\rm End}(W)$ is an isomorphism. Such a module $W$ is simple of dimension $2^{r-1}|2^{r-1}$, and its restriction to the even algebra ${\rm C}(V)_0 \simeq {\rm End}(W_0) \times  {\rm End}(W_1)$ is the direct-sum of two nonisomorphic simple modules $W=W_0 \oplus W_1$ called half-spinor modules for $V$. The $k$-linear representations of the groups ${\rm GSpin}(V)$ and ${\rm Spin}(V)$ obtained by restriction of those modules $W$ and $W_i$ are respectively called spin and half-spin representations. \ps \ps
Case (b).  If $\dim V=2r+1$ is odd, a spinor module for $V$ is a simple module for the (matrix) $k$-algebra ${\rm C}(V)_0$. The restriction to ${\rm GSpin}(V)$ or ${\rm Spin}(V)$ of a spinor module for $V$ is by definition a spin representation. \ps\ps

Any hyperbolic quadratic space, or more generally any regular quadratic space containing a subspace of codimension $\leq 1$ which is hyperbolic, has a trivial even Clifford algebra. However, the converse does not hold:  if $k=\R$ and $V$ is positive definite, it is well-known that $V$ has a trivial even Clifford algebra if, and only if, we have $\dim V \equiv -1,0,1 \bmod 8$ \cite[table p. 103]{deligne}. \ps\ps

\begin{lemme}\label{lemmeregularhyp} Let $V$ be a nondegenerate quadratic space of even dimension, $e \in V$ an element such that ${\rm q}(e) \neq 0$, ${\rm O}(V)_e$ the stabilizer of $e$ in ${\rm O}(V)$, ${\rm SO}(V)_e={\rm SO}(V)\cap {\rm O}(V)_e$, and ${\rm GSpin}(V)_e$ the inverse image of ${\rm SO}(V)_e$ in ${\rm GSpin}(V)$ under the map $\pi_V$. \begin{itemize}
\item[(i)] The orthogonal $E$ of $e$ in $V$ is a regular quadratic subspace, and the natural map ${\rm O}(V)_e \rightarrow {\rm O}(E)$ induces a bijection $a: {\rm SO}(V)_e \isomo {\rm SO}(E)$.\ps \ps
\item[(ii)] The natural map $b : {\rm C}(E) \rightarrow {\rm C}(V)$ is injective and induces a group isomorphism ${\rm GSpin}(E) \rightarrow {\rm GSpin}(V)_e$. Furthermore, we have the equality $\pi_V \circ b = a^{-1} \circ \pi_E$ of maps ${\rm GSpin}(E) \rightarrow {\rm SO}(V)$. \ps \ps
\item[(iii)] If $V$ has a trivial even Clifford algebra, then so does $E$, and the restriction to $b: {\rm C}(E)_0 \rightarrow {\rm C}(V)_0$ of a half-spinor module for $V$ is a spinor module for $E$. 
\end{itemize}
\end{lemme}
\begin{pf} When $2$ is invertible in $k$, then $V = ke \bot E$ and the lemma is fairly standard. By lack of reference, we provide a proof which works in characteristic $2$ as well. We may and do choose a nondegenerate hyperplan $H \subset E$. Denote by $H^\bot$ the orthogonal of $H$ inside $V$. We have $\dim H^\bot = 2$, $V = H \oplus H^\bot$,  $e \in H^\bot$, and $E = H \oplus ke$. \ps\ps
The first statement of assertion (i) is immediate from the definitions. The second one is clear when $2 \in k^\times$. If $k$ has characteristic $2$, the surjectivity of the natural map ${\rm O}(V)_e \rightarrow {\rm O}(E)={\rm SO}(E)$ is a consequence of the Witt theorem (which is due to Arf in characteristic $2$ \cite[\S 4\, {\rm no}. 3 Thm. 1]{bou}). The kernel of this map is $(1 \times {\rm O}(H^\bot)) \cap {\rm O}(V)_e$, which is easily checked to be the group of order $2$ generated by the improper isometry $x \mapsto x - \frac{\beta(x,e)}{{\rm q}(e)} e$, which proves (i). \ps\ps
The map $b$ of the statement is the one induced by applying the ``Clifford functor'' to the inclusion $E \rightarrow V$, it is injective by \cite[\S 9 {\rm no}. 3 Cor. 3]{bou}. Let $\gamma \in {\rm GSpin}(E)$. We have $b(\gamma) \in {\rm C}(V)_0^\times$ and we claim that $b(\gamma) \in {\rm GSpin}(V)$. By surjectivity of $\pi_V$ in the sequence \eqref{sexgspin}, we may choose $h \in {\rm GSpin}(V)_e$ such that $\pi_V(h) = a^{-1}\circ \pi_E(\gamma)$. If we can show that the even element $\lambda := b(\gamma)^{-1} h \in {\rm C}(V)_0^\times$ is a scalar, then (ii) follows at once. The relation $\pi_V(h) = a^{-1}\circ \pi_E(\gamma)$ asserts that $\lambda$ commutes with $E=H \oplus k e$. The natural isomorphism of graded algebras ${\rm C}(H) \otimes^{\rm gr} {\rm C}(H^\bot) \isomo {\rm C}(V)$ identifies the supercentralizer of $H$ in ${\rm C}(V)$ to ${\rm C}(H^\bot)$. It follows that $\lambda$ is identified with an element of the graded quaternion algebra ${\rm C}(H^\bot)$ which is even and commutes with some nonzero odd element $e \in H^\bot = {\rm C}(H^\bot)_1$ (we have  $e^2 ={\rm q}(e) \neq 0$). But the commutator of such an element is $k[e]$, and $k[e] \cap {\rm C}(H^\bot)_0=k$, which proves $\lambda \in k^\ast$ and assertion (ii).\ps\ps
The map $b$ induces a $k$-algebra morphism ${\rm C}(E)_0 \rightarrow {\rm C}(V)_0$. If $V$ has trivial even Clifford algebra, then ${\rm C}(V)_0$ is isomorphic to a direct product $A_1 \times A_2$ of two matrix $k$-algebras $A_i$, each of them being of rank $2^{\dim V-2}=2^{\dim E-1}$. On the other hand, as $E$ is regular the $k$-algebra ${\rm C}(E)_0$ is central and simple of rank $2^{\dim E-1}$ as well. It follows that each projection ${\rm C}(E)_0 \rightarrow A_i$ is an isomorphism, which proves (iii). \end{pf}\ps\ps

For later use, we end this paragraph by stating two classical results.
 \ps
\begin{prop}\label{cliffordhyp} Let $I$ be a finite dimensional $k$-vector space, $V={\rm H(I)}$ the hyperbolic quadratic space over $I$ and $\rho_I : {\rm GL}(I) \rightarrow {\rm SO}({\rm H}(I))$ the natural homomorphism defined by $({\rho_I(g)}_{|I},{\rho_I(g)}_{|I^\ast})= (g,{}^{\rm t}\!g^{-1})$. Then the exterior superalgebra $\Lambda\, I$ is in a natural way a spinor module for $V$. Moreover, there is a natural homomorphism $$\widetilde{\rho}_I :  {\rm GL}(I) \rightarrow {\rm GSpin}({\rm H}(I))$$ such that $\pi_V \circ \widetilde{\rho}_I = \rho_I$, $\nu_V \circ \widetilde{\rho}_I = \det$, and such that the restriction of the spin representation $\Lambda \,I$ to $\widetilde{\rho}_I$ is the canonical representation of ${\rm GL}(I)$ on $\Lambda \,I$. \end{prop}
\ps\ps
\begin{pf}The assertion on the exterior algebra is well-known, and proved {\it e.g.} in \cite[\S 9 {\rm no} 4]{bou} or \cite[Ch. IV Prop. 2.1.1]{knusquadherm}. The other assertions are proved in \cite[Ch. IV \S 6.6]{knusquadherm}.\end{pf}

\ps\ps
\begin{prop}\label{decospin} Let $U,V$ be regular quadratic spaces over $k$ with $\dim U$ even and $\dim V$ odd. If  $U$ and $V$ have a trivial even Clifford algebra, then so does the regular quadratic space $U \bot W$. Moreover, the restriction to the canonical $k$-algebra morphism $${\rm C}(U)_0 \otimes {\rm C}(V)_0 \longrightarrow {\rm C}(U \bot V)_0$$ of a spinor module for $U \bot V$ is isomorphic to the (external) tensor product of a spinor module for $U$ and of a spinor module for $V$.
\end{prop}\ps
\begin{pf} Set $R={\rm C}(U)_0$, $S={\rm C}(V)_0$ and $T={\rm C}(U \bot V)_0$. We have a natural injective $k$-algebra homomorphism $R \otimes S \rightarrow T$.
Write $\dim U = 2a$ and $\dim V = 2b+1$. The we have $R \simeq {\rm M}_{2^{a-1}}(k) \times {\rm M}_{2^{a-1}}(k)$, $S \simeq  {\rm M}_{2^b}(k)$, and $T$ is central simple of dimension $2^{2a+2b}$.  As $T$ contains $R \otimes S$, which contains a $k$-algebra isomorphic to $k^{2^{a+b}}$, we have $T \simeq {\rm M}_{2^{a+b}}(k)$. We conclude using the following elementary fact: if $p$ and $q$ are integers, there is a unique ${\rm GL}_{p+q}(k)$-conjugacy class (resp.  ${\rm GL}_{pq}(k)$-conjugacy class) of $k$-algebra morphisms ${\rm M}_p(k) \times {\rm M}_q(k) \rightarrow {\rm M}_{p+q}(k)$ (resp. ${\rm M}_p(k) \otimes {\rm M}_q(k) \rightarrow {\rm M}_{pq}(k)$). 
\end{pf}

\ps\ps
\subsection{The Witt index of a representation in ${\rm O}(V)$}\label{defwindex}
\ps\ps

In this paragraph, we assume that $V$ is a nondegenerate quadratic space over the field $k$, that $\Gamma$ is a group and that
$\rho : \Gamma \rightarrow {\rm O}(V)$ is a group homomorphism. \ps\ps

\begin{definition}\label{defindex} The Witt index of $\rho$ is the maximal dimension of a subspace $I \subset V$ which is both totally isotropic {\rm (}{\rm i.e.} ${\rm q}(I)=0${\rm )} and stable by $\rho(\Gamma)$. When $\rho$ is the inclusion of a subgroup $\Gamma$ of ${\rm O}(V)$ we define the Witt index of $\Gamma$ as the one of $\rho$. 
\end{definition}

Denote by ${\rm Irr}(X)$ the set of isomorphism classes of the irreducible constituents of the $k[\Gamma]$-module $X$. As $V$ is isomorphic to $V^\ast$, the duality functor induces an action of $\Z/2\Z$ on ${\rm Irr}(V)$ (the ``duality'' action).

 \begin{prop} \label{propwindex} Assume $k$ is algebraically closed of characteristic $\neq 2$ and that $V$ is semisimple as a $k[\Gamma]$-module. Fix $\Phi \subset {\rm Irr}(V)$ a set of representatives for the duality action.  There exists a decomposition
\begin{equation} \label{candeco} V \,=\, ( I \oplus J )\, \bot \, W \end{equation}
where $I,J$ and $W$ are $\Gamma$-stable subspaces of $V$ such that: \begin{itemize}\ps\ps
\item[(i)]  $W$ is the orthogonal direct sum of nonisomorphic irreducible selfdual $k[\Gamma]$-modules, \ps\ps
\item[(ii)] $I$ and $J$ are totally isotropic, and we have  ${\rm Irr}(I) \subset \Phi$.\ps\ps
\end{itemize}
For any such decomposition, $\beta_V$ induces a $\Gamma$-equivariant isomorphism $I \isomo J^\ast$ and the dimension of $I$ is the Witt index of $\rho$. Morever, the isomorphism class of the $k[\Gamma]$-module $I$ does not depend on the choice of the decomposition \eqref{candeco}, and the isomorphism class of the $k[\Gamma]$-module $W$ depends neither on the choice of  \eqref{candeco} nor on the choice of $\Phi$. \end{prop}

{\scriptsize \begin{pf} (This is presumably well-known so we use a small font.) Assume first that there is a $\Gamma$-stable totally isotropic subspace $I \subset V$ such that $I^\bot = I$ (so that $\dim I = \frac{\dim V}{2}$ is the Witt index of $\rho$). We want to show the existence of a decomposition \eqref{candeco} satisfying (i) and (ii). As $V$ is semisimple we may find a $\Gamma$-stable subspace $F \subset V$ such that $V = I \oplus F$. The perfect pairing $\beta_V$ induces a $k[\Gamma]$-linear isomorphism $I \isomo F^\ast$, $x \mapsto (y \mapsto \beta_V(x,y))$. This shows that for $x \in F$, there exists a unique element $u(x) \in I$ such that $\beta_V(x,y) = \beta_V(u(x),y)$ for all $y \in F$. The map $u : F \rightarrow I$ is $k[\Gamma]$-linear. Observe that for all $x \in F$ we have 
$${\rm q}(x - \frac{u(x)}{2}) = {\rm q}(x) - \frac{1}{2}\,\beta_V(u(x),x) = \frac{1}{2} (\beta_V(x,x) - \beta_V(u(x),x)) = 0.$$
As a consequence, the subspace $J = \{x - \frac{u(x)}{2}, x \in F\}$ is totally isotropic, $\Gamma$-stable, and satisfies $I \oplus J = V$. This shows the existence of a decomposition \eqref{candeco} satisfying (i) and (ii) except perhaps ${\rm Irr}(I) \subset \Phi$. As $V$ is semisimple, we may find a $k[\Gamma]$-module decomposition $I=I_0 \oplus I_1$ (resp. $J = J_0 \oplus J_1$) where $I_0$ (resp. $J_0)$ is the sum of the irreducible nonselfdual subrepresentations of $I$ (resp. $J$) whose isomorphism class do not (resp. do) belong to $\Phi$. By construction, any $k[\Gamma]$-linear morphism $J_0 \rightarrow I_1^\ast$ is zero, so we have $\beta_V(J_0,I_1)=0$, and $\beta_V$ induces an isomorphism $I_0^\ast \rightarrow J_0$. As a consequence, $I'=J_0 \oplus I_1$ is a $\Gamma$-stable totally isotropic subspace of $V$ of dimension $\frac{\dim V}{2}$, hence such that $I'^\bot = I'$, which satisfies furthermore ${\rm Irr}(I') \subset \Phi$. We conclude the proof by applying to $I'$ instead of $I$ the previous argument. \ps\ps

We now claim the existence of a decomposition \eqref{candeco} satisfying (i), (ii) and such that $\dim I$ is the Witt index of $\rho$. Let $I$ be a $\Gamma$-stable totally isotropic subspace of $V$. By the semisimplicity assumption, we may choose a $\Gamma$-stable subspace $V' \subset I^\bot$ such that we have $I^\bot = I \oplus V'$. This forces $V'$ to be a nondegenerate subspace of $V$. By construction, we have $I^\bot \cap {V'}^\bot = I$. The previous paragraphs apply thus to the representation of $\Gamma$ on ${V'}^\bot$ and show the existence of a $\Gamma$-stable subspace $J$ such that ${V'}^\bot = I \oplus J$ satisfying (ii). Replacing $V$ by $V'$, we may thus assume that $V$ has no nonzero $\Gamma$-stable totally isotropic subspace. We have to show that $V$ is the direct sum of nonisomorphic irreducible selfdual $k[\Gamma]$-modules. \ps\ps

Let $U \subset V$ be any $\Gamma$-stable subspace. The isomorphism $\beta_V : V \isomo V^\ast$ induces a $\Gamma$-equivariant morphism $U \rightarrow U^\ast$ whose kernel is $U\cap U^\bot$. This later space is totally isotropic and $\Gamma$-stable, hence zero by assumption. It follows that each $\Gamma$-stable subspace of $V$ is nondegenerate, and in particular, selfdual as a $k[\Gamma]$-module. Assume that there are $\Gamma$-stable, orthogonal, subspaces $U,U' \subset V$ such that $U$ and $U'$ are isomorphic and irreducible as $k[\Gamma]$-modules. Choose a $k[\Gamma]$-linear isomorphism $f: U \rightarrow U'$. Both $\beta_V$ and $\beta_V \circ f$ define nondegenerate and $\Gamma$-invariant symmetric bilinear forms on $U$. As $U$ is irreducible and ${\rm char}\, k \neq 2$, there is a unique such form up to a scalar, so there is $\lambda \in k^\times$ such that we have ${\rm q}(u) = \lambda \,{\rm q}(f(u))$ for all $u \in U$. As $k$ is algebraically closed, there is $\mu \in k^\times$ such that $\lambda = -\mu^2$. It follows that $\{u + \mu f(u), u \in U\} \subset U \oplus U' \subset V$ is a $\Gamma$-stable totally isotropic subspace, a contradiction. This proves the claim. \ps\ps

We now prove the last assertions regarding an arbitrary decomposition \eqref{candeco} satisfying (i) and (ii). Observe that the $k[\Gamma]$-module $W$ is isomorphic to the direct sum of the irreducible selfdual representations of $\Gamma$ which occur in $V$ with an odd multiplicity. In particular, its isomorphism class does not depend on the choice of \eqref{candeco} or $\Phi$. As a consequence, the isomorphism class of $I \oplus I^\ast$ has the same property, and the one of $I$ is uniquely determined by $\Phi$. In particular, $\dim I$ does not depend either on the choice of \eqref{candeco}. This concludes the proof, as we have shown the existence of a decomposition \eqref{candeco} satisfying (i), (ii) and such that $\dim I$ is the Witt index of $\rho$. \end{pf}}

\begin{cor} \label{corwindex} Under the assumptions of Proposition \ref{propwindex}, the ${\rm O}(V)$-conjugacy class of $\rho$ only depends on the isomorphism class of the $k[\Gamma]$-module $V$. If furthermore $\dim V$ is odd, the same assertion also holds with ${\rm O}(V)$ {\rm (}$ = \{\pm 1\} \times  {\rm SO}(V)${\rm )} replaced by ${\rm SO}(V)$. \ps\ps
\end{cor}

{\scriptsize \begin{pf} For $i=1,2$ let $V_i$ be a nondegenerate quadratic space over $k$ (algebraically closed of characteristic $\neq 2$) and $\rho_i : \Gamma \rightarrow {\rm O}(V_i)$ a group homomorphism. We assume that the $k[\Gamma]$-modules $V_i$ are semisimple and isomorphic. We have to show that  there is an isometry $u :  V_1 \rightarrow V_2$ such that $\rho_2 (\gamma) \circ u  = u \circ \rho_1(\gamma)$ for all $\gamma \in \Gamma$. By Proposition \ref{propwindex}, we may assume that we have $V_i = {\rm H}(I_i) \bot W_i$, where 
 $I_i,I_i^\ast$ and $W_i$ are preserved by $\rho_i(\Gamma)$ and the $W_i$ have Witt index $0$. By that proposition, we have $W_1 \simeq W_2$ and $I_1 \simeq I_2$ as $k[\Gamma]$-modules. Any $k[\Gamma]$-equivariant isomorphism $I_1 \isomo I_2$ induces a $\Gamma$-equivariant isometry ${\rm H}(I_1) \isomo {\rm H}(I_2)$. We may thus assume $V_i=W_i$ for each $i$. In this case, the isotypical components of the $V_i$ are nondegenerate, so we may even assume that each $V_i$ is irreducible. This implies that the $\Gamma$-invariants of $({\rm Sym}^2 V_i)^\ast$ have dimension $1$. It follows that if $u$ is any $k[\Gamma]$-linear isomorphism $V_1 \isomo V_2$, there is $\lambda \in k^\times$ such that $\lambda u$ is an isometry $V_1 \isomo V_2$. \end{pf}}

\ps

\begin{remark}\label{windexgo}
{\rm Assume more generally that $\rho$ is a morphism $\Gamma \rightarrow {\rm GO}(V)$, with similitude factor $\mu$. Then the considerations of this paragraph extend verbatim to this setting if one replaces the duality functor $X \mapsto X^\ast$ on $k[\Gamma]$-modules by the functor $X \mapsto X^\ast \otimes \mu$. 
}
\end{remark}
\ps\ps
\section{Octonion algebras and ${\rm G}_2$-subgroups}
\ps\ps
\subsection{The Clifford algebra of an octonion algebra} \label{octoalg}Let $C$ be an octonion $k$-algebra. Recall that this is a nondegenerate $8$-dimensional quadratic space over $k$ equipped with an element $e \in C$ and a $k$-bilinear map $C \times C \rightarrow C$, $(x,y) \mapsto xy$, such that for all $x,y \in C$ we have ${\rm q}(xy)={\rm q}(x){\rm q}(y)$ and $xe=ex=x$ \cite{blijspringer1}\,\cite[Ch 1]{springer}. \ps\ps

In order to discuss some features of the spin representations of the group ${\rm Spin}(C)$, let us recall the description of the Clifford algebra of $C$ given in \cite[\S 3.6 p. 61]{springer} and \cite[Prop. 4.4 ]{kps}. Consider the $k$-supervector space $\mathcal{C}$ with $\mathcal{C}_0=\mathcal{C}_1=C$ (so that $\dim \mathcal{C} = 8|8$). 
If $x \in C$, denote by $\ell(x)$ the odd endomorphism of $\mathcal{C}$ sending $(a,b)$ to $(xb,\overline{x}a)$, where $x \mapsto \overline{x}:= x - \beta_C(x,e)e$ is the canonical involutive automorphism of $C$ \cite{blijspringer1}. Then we have $\ell(x) \circ \ell(x) = {\rm q}(x)\,  {\rm Id}$, so the $k$-linear map $\ell : C \rightarrow {\rm End}(\mathcal{C})$ 
extends to a $k$-superalgebra morphism 
\begin{equation}\label{isolc}
\ell_C : {\rm C}(C) \rightarrow {\rm End}(\mathcal{C}).
\end{equation}
This is an isomorphism as the left-hand side is central simple and both $k$-algebras have the same dimension.  \ps\ps

\begin{scholie} \label{scholtrivclif}If $C$ is an octonion $k$-algebra then $\ell_C$ is an isomorphism of $\Z/2\Z$-graded $k$-algebras. In particular, $C$ has a trivial even Clifford algebra, and $\mathcal{C}$ (resp. each copy of $C$ in it) is in a natural way a spinor module (resp. half-spinor modules) for $C$.
\end{scholie}\ps\ps

View $\mathcal{C}$ as a nondegenerate quadratic space for the quadratic form ${\rm q} \bot {\rm q}$ (orthogonal sum). Let $v,w \in C$ be two nonisotropic elements and consider their product $\gamma_{v,w}:=v \otimes w \in {\rm GSpin}(C)$.  We have for all $a,b \in C$ the formula
\begin{equation}\label{formoctell} \ell_C(\gamma_{v,w})(a,b) = ( v(\overline{w}a), \overline{v}(wb)).\end{equation}
In particular, $\ell_C(\gamma_{v,w})$ is a similitude of $\mathcal{C}$ of factor ${\rm q}(v){\rm q}(w)=\nu_C(\gamma_{v,w})$.  Recall that the Cartan-Dieudonn\'e theorem \cite[Prop. 8 \& Prop. 14]{dieudonne} asserts that the nonisotropic elements of $C$ generate $\Gamma(C)$, which implies that the elements of the form $\gamma_{v,w}$ generate ${\rm GSpin}(C)$. \ps\ps

\begin{scholie} For $i=0,1$, the space $\mathcal{C}_i=C$ viewed as a half-spinor module for $C$ defines a morphism $\rho_{C ; i} : {\rm GSpin}(C) \rightarrow {\rm GSO}(C)$ with similitude factor $\nu_C$.
\end{scholie}

Define a group homomorphism $$\tau_C : {\rm GSpin}(C) \rightarrow {\rm SO}(C) \times {\rm GSO}(C) \times {\rm GSO}(C)$$ by the formula $\tau_C(\gamma)=(\pi_C(\gamma),\rho_{C; 0}(\gamma),\nu_C(\gamma)^{-1}\rho_{C;1}(\gamma))$ (note the scalar factor  $\nu_C(\gamma)^{-1}$ in the last component). Following \cite{springer}, we say that a triple $(t_1,t_2,t_3)  \in {\rm SO}(C) \times {\rm GSO}(C) \times {\rm GSO}(C)$ is {\it related} if we have $$t_1(xy)=t_2(x)t_3(y)\, \,\, \, \, \, \, \forall x,y \in C.$$As is shown in \cite[Thm. 3.2.1 (i) \& (iii)]{springer}, the subset ${\rm GRT}(C)$ of related triples is a subgroup of ${\rm SO}(C) \times {\rm GSO}(C) \times {\rm GSO}(C)$. The following proposition is a modest variant of \cite[Prop. 3.6.3]{springer}.\ps\ps

\begin{prop} \label{groupeGRT} Let $C$ be an octonion $k$-algebra. The morphism $\tau_C$ induces an isomorphism ${\rm GSpin}(C) \isomo {\rm GRT}(C)$. It restricts to an isomorphism between ${\rm Spin}(C)$ and the subgroup of all related triples in ${\rm SO}(C)^3$. \end{prop}
\ps\ps

\begin{pf} The kernel of $\tau_C$ is trivial as $\ker \ell_C = 0$. By \cite[Thm. 3.2.1]{springer} assertion (iii) and formula \eqref{formoctell}, we have $\tau_C(\gamma_{v,w}) \in {\rm GRT}(C)$, hence ${\rm Im} \tau_C \subset {\rm GRT}(C)$. The same theorem shows that the first projection ${\rm pr}_1 : {\rm GRT}(C) \rightarrow {\rm SO}(C)$, $(t_1,t_2,t_3) \mapsto t_1$, is surjective, and that $\ker \, {\rm pr}_1$ is the subgroup of elements of the form $(1,m_{\lambda},m_{\lambda^{-1}})$ where $\lambda \in k^\times$ and $m_{\lambda} : x \mapsto \lambda x$. But this kernel is exactly $\tau_C(k^\times)$ and the map ${\rm pr}_1 \circ \tau_C = \pi_C$ is surjective. This proves ${\rm Im}\, \tau_C = {\rm GRT}(C)$ and the last assertion follows. \end{pf}

\ps\ps
\subsection{The automorphism group of an octonion algebra} \label{g2emb} Let $C$ be an octonion $k$-algebra. We denote by ${\rm G}_{\tiny 2}(C)$ the automorphism group of the octonion $k$-algebra $C$. By \cite[Cor. 2.2.5]{springer},  ${\rm G}_{\tiny 2}(C) \subset {\rm O}(C)$ is actually a subgroup of ${\rm SO}(C)$. We shall denote by $${\rm i}_C : {\rm G}_{\tiny 2}(C) \longrightarrow {\rm SO}(C)$$
this inclusion.  Consider the ``diagonal'' group homomorphism $\iota_C : {\rm G}_{\tiny 2}(C) \rightarrow {\rm SO}(C) \times {\rm SO}( C)$ sending an automorphism $\phi$ of $C$ to $(a,b) \mapsto (\phi(a),\phi(b))$ (an even $k$-linear automorphism of $\mathcal{C}$). The isomorphism \eqref{isolc} allows to define a map $\widetilde{{\rm i}}_C:=\ell_C^{-1} \circ \iota_C : {\rm G}_2(C) \rightarrow {\rm C}(C)_0^\times.$

 \ps 
\begin{prop} \label{iclift} The map $\widetilde{{\rm i}}_C$ is a group homomorphism ${\rm G}_2(C) \rightarrow {\rm Spin}(C)$ such that $\pi_C \circ \widetilde{{\rm i}}_C = {\rm i}_C$.  
\end{prop}

\begin{pf} Observe that for $x \in C$ and $\phi \in {\rm G}_{\tiny 2}(C)$, we have the relation $\iota_C(\phi) \ell_C(x) \iota_C(\phi)^{-1} = \ell_C(\phi(x))$. It follows that $\widetilde{{\rm i}}_C(\phi):=\ell_C^{-1} \circ \iota_C (\phi)$ is an element of ${\rm GSpin}(C)$ which satisfies $\pi_C(\widetilde{{\rm i}}_C(\phi))={\rm i}_C(\phi)$. We have $\nu_C(\widetilde{{\rm i}}_C(\phi))=1$ as the Chevalley anti-involution of ${\rm C}(C)$ corresponds via $\ell_C$ to the adjonction on ${\rm End}(C \oplus C)$ with respect to $\beta_{C \oplus C}$ by \cite[Prop. 4.4]{kps}. \end{pf}  \ps\ps

The subspace $P \subset C$ of {\it pure octonions} is defined as the orthogonal of the unit element $e$ for the bilinear form $\beta_C$. As ${\rm q}(e)=1 \neq 0$, we are in the context of Lemma~\ref{lemmeregularhyp}. In particular, $P$ is a regular quadratic subspace of $C$ of dimension $7$ and we have natural isomorphisms $a : {\rm SO}(C)_e \isomo {\rm SO}(P)$ and $b : {\rm GSpin}(P) \isomo {\rm GSpin}(C)_e$ satisfying $\pi_C \circ b = a^{-1} \circ \pi_P$. As ${\rm G}_2(C)$ preserves $e$ by definition, we get natural injective morphisms $${\rm j}_C:=a \circ {\rm i}_C : {\rm G}_2(C) \longrightarrow {\rm SO}(P) \, \, \, {\rm and}\, \, \, \widetilde{{\rm j}}_C= b^{-1} \circ \widetilde{{\rm i}}_C : {\rm G}_2(C) \longrightarrow {\rm Spin}(P)$$
satisfying $\pi_P \circ \widetilde{{\rm j}}_C = {\rm j}_C$. By Lemma~\ref{lemmeregularhyp} (iii), a half-spinor module for $C$ restricts via $b$ to a spinor module for $P$, and we shall denote by $$\rho_C := \rho_{C ; 0} \circ b : {\rm GSpin}(P) \longrightarrow {\rm GSO}(C)$$ the (canonical) associated spin representation of ${\rm GSpin}(P)$. \ps\ps

\begin{prop}\label{propgspin7} Let $C$ be an octonion algebra over $k$, $P \subset C$ the $7$-dimensional quadratic space of pure octonions, $\widetilde{{\rm j}}_C : {\rm G}_2(C) \longrightarrow {\rm Spin}(P)$ be the canonical morphism, and $\rho_C : {\rm GSpin}(P) \longrightarrow {\rm GSO}(C)$ the canonical spin representation of ${\rm GSpin}(P)$. Then: \begin{itemize}
\item[(a)] $\rho_C \circ  \widetilde{{\rm j}}_C$ coincides with the canonical morphism ${\rm i}_C : {\rm G}_2(C) \rightarrow {\rm SO}(C)$; \ps \ps
\item[(b)] the action of ${\rm GSpin}(P)$ on the set ${\rm Q}=\{v \in C, {\rm q}(v) \neq 0\}$, defined by $(\gamma,v) \mapsto \rho_C(\gamma)v$, is transitive;\ps\ps
\item[(c)] the stabilizer of the unit element $e \in {\rm Q}$ in ${\rm GSpin}(P)$ is ${\rm i}_C({\rm G}_2(C))$, and the stabilizer of $ke \subset C$ in ${\rm GSpin}(P)$ is $k^\times \times  {\rm i}_C({\rm G}_2(C))$. \ps\ps
\end{itemize}
\end{prop}
\ps\ps
\begin{pf} Assertion (a) is obvious from the definitions. To show (b), observe first that the action of the statement is well-defined as we have $\rho_C({\rm GSpin}(P)) \subset {\rm GSO}(C)$. In \cite[Lemma 3.4.2]{springer}, the authors show that for any $v \in C$ such that ${\rm q}(v) \neq 0$, there is a related triple $(t_1,t_2,t_3) \in {\rm GRT}(C)$ such that $t_1 \in {\rm SO}(C)_e$ and $t_2(e)=v$. This is exactly assertion (b), by Proposition \ref{groupeGRT}. Moreover, they also show in Proposition 3.4.1 {\it loc. cit.} that for an element $(t_1,t_2,t_3)$ in ${\rm GRT}(C)$ such that $t_1 \in {\rm SO}(C)_e$, we have $t_2(e) \subset k e$ if, and only if, we have $t_1 \in {\rm G}_2(C)$. It follows that the stabilizer of $e \in {\rm Q}$ in ${\rm GSpin}(P)$ is included in $k^\times \cdot \widetilde{{\rm j}}_C({\rm G}_2(C))$. But this stabilizer contains 
$\widetilde{{\rm j}}_C({\rm G}_2(C))$ by (a), and intersects $k^\times$ trivially by construction, which proves (c). 
\end{pf}
\ps\ps
\begin{example} \label{examplecayley}{\rm  Assume $k=\R$ and $C$ is the Euclidean (Cayley-Graves') octonion $\R$-algebra. Then ${\rm Spin}(P)$ is isomorphic to the compact Lie group ${\rm Spin}(7)$, the map ${\rm GSpin}(P) \rightarrow {\rm Spin}(P)$ is surjective, and ${\rm SO}(C)$ is isomorphic to ${\rm SO}(8)$. The spin representation ${\rm Spin}(P) \rightarrow {\rm SO}(C)$ defines a group homomorphism ${\rm Spin}(7)\rightarrow {\rm SO}(8)$. The proposition asserts that ${\rm Spin}(P)$ acts transitively on the Euclidean unit sphere $S$ ($ \simeq \mathbb{S}^7$) of $C$, and the stabilizers of this action form a conjugacy class of compact connected subgroups of ${\rm Spin}(P)$ of type ${\bf G}_2$. See \cite[Thm. 5.5]{adams} for another proof of this classical result. Let us add that any compact connected subgroup $H \subset  {\rm Spin}(P)$ which is of type ${\rm G}_2$ fixes a point in $S$ (hence belongs to the aforementioned conjugacy class). Indeed, there is up to isomorphism a unique nontrivial irreducible complex representations of $H$ of dimension $\leq 8$, which is the $7$-dimensional representation defined by $\pi_P$ on $P \otimes \C$ (and thus defined over $\R$), so the natural representation of $H$ on $C$ contains the trivial representation. In the next paragraph, we prove more a general result along those lines, which applies to an arbitrary field $k$. 
} 
\end{example}
\ps\ps

\subsection{${\rm G}_2$-subgroups of ${\rm SO}(P)$ and ${\rm Spin}(P)$ with $\dim P=7$} \label{defg2sb}

\begin{definition} Let $P$ be a regular quadratic space of dimension $7$ over a field $k$ and let $\Gamma$ be a subgroup of ${\rm SO}(P)$ {\rm (}resp. ${\rm Spin}(P)${\rm )}. We say that $\Gamma$ is a ${\rm G_2}$-subgroup if there exists an octonion $k$-algebra $C$ and an isometric embedding of $P \hookrightarrow C$ onto the space of pure octonions of $C$, such that we have $\Gamma = {\rm j}_C({\rm G}_2(C))$ {\rm (}resp. $\Gamma= \widetilde{{\rm j}}_C({\rm G_2}(C))${\rm )}.
\end{definition}\ps\ps

Note that there may be no such subgroup in general, {\it i.e.} no octonion $k$-algebra with pure subspace isometric to $P$. They do exists, for instance, if $P$ is the orthogonal of an element of norm $1$ in ${\rm H}(k^4)$ (hence if $k$ is algebraically closed), as ${\rm H}(k^4)$ may be endowed with a structure of split octonion $k$-algebra. In other words, ${\rm SO}_7(k)$ and ${\rm Spin}_7(k)$ have ${\rm G}_2$-subgroups. 
\ps\ps
\begin{prop}\label{carG2sbgp}  Let $P$ be a regular quadratic space of dimension $7$ over $k$, the set of ${\rm G}_2$-subgroups of ${\rm SO}(P)$ {\rm (}resp. ${\rm Spin}(P)${\rm )} form a single ${\rm SO}(P)$-conjugacy class {\rm (}resp. ${\rm GSpin}(P)$-conjugacy class{\rm )}. 
\end{prop}
\ps\ps
\begin{pf} For $i=1,2$, consider an octonion $k$-algebra $C_i$, with neutral element $e_i$, and an isometry $\mu_i : P \isomo e_i^\bot$.  We show first that $C_i$ are isomorphic as octonion $k$-algebras. By \cite[Assertion (2.3)]{blijspringer1}, it is enough to show that they are isometric as quadratic spaces. This is obviously true if ${\rm char}\, k \,\neq 2$ as we have $C_i \simeq P \bot \langle 1\rangle $. If ${\rm char}\, k \, = 2$, and if $H$ is a nondegenerate hyperplane of $P$, then $C_i$ is isometric to the orthogonal sum of $H$ and of a $2$-dimensional nondegenerate quadratic space over $k$ which represents $1$. We conclude from the triviality of the Arf invariant of $C_i$, which in turn follows from the isomorphism \eqref{isolc}. \ps\ps
Unravelling the definitions, the ${\rm SO}(P)$ case of the statement is now obvious. In the ${\rm Spin}(P)$ case, it remains to show the following fact: let $\alpha \in {\rm SO}(P)$ and ${\rm C}(\alpha)$ be the associated even automorphism of ${\rm C}(P)$ defined by functoriality of the Clifford algebra construction, then ${\rm C}(\alpha)$ is an inner automorphism defined by some element in ${\rm GSpin}(P)$. But this follows from the Skolem-Noether theorem as: ${\rm C}(P)_0$ is a matrix $k$-algebra (Scholium \ref{scholtrivclif}), ${\rm C}(\alpha)$ acts trivially on the center $Z$ of ${\rm C}(P)$ (as $\alpha$ belongs to ${\rm SO}(P)$), and the $k$-algebra ${\rm C}(P)$ is generated by ${\rm C}(P)_0$ and $Z$.\end{pf}
\ps\ps

\begin{prop}\label{algcarG2} Let $P$ be a regular quadratic space of rank $7$ over an algebraically closed field $k$, and let $\Gamma \subset {\rm SO}(P)$ {\rm (}resp. $\Gamma \subset {\rm Spin}(P)${\rm )} be a subgroup. The following properties are equivalent: 
\begin{itemize}\ps\ps
\item[(i)] $\Gamma$ is a ${\rm G}_2$-subgroup,\ps\ps
\item[(ii)] $\Gamma$ is a closed connected subgroup, which is a simple linear algebraic group of type ${\bf G}_2$.\ps\ps
\end{itemize}
\end{prop}
\ps\ps
{\scriptsize
\begin{pf} If $C$ is an octonion $k$-algebra, its automorphism group ${\rm G}_2(C)$ has a unique structure of linear algebraic group such that the morphism ${\rm i}_C : {\rm G}_2(C) \rightarrow {\rm SO}(C)$ is a closed immersion.  By \cite[Thm. 2.3.5]{springer}, it is connected, simple, of type ${\rm G}_2$ and the morphism ${\rm j}_C : {\rm G}_2(C) \rightarrow {\rm SO}(P)$ is a closed immersion as well. The morphism $\widetilde{{\rm j}}_C : {\rm G}_2(C)  \rightarrow {\rm Spin}(P)$ defined above, which is a morphism of algebraic groups, satisfies $\pi_P \circ \widetilde{{\rm j}}_C = {\rm j}_C$, hence is a closed immersion as well. This proves (i) $\Rightarrow$ (ii). \ps\ps

Before proving (ii) $\Rightarrow$ (i) let us recall some properties of the (finite dimensional, algebraic) linear representations of reductive groups over a field $k$ of arbitrary characteristic, following \cite[Ch. XI]{humphreys1} and \cite[Chap. 2 \& 3]{humphreys2} (see also \cite{rags}). Let $G$ be a semisimple algebraic group over $k$. If $\lambda$ is a dominant weight of $G$, the category of linear representations of $G$ equipped with an extremal vector of weight $\lambda$ has an initial object denoted $V(\lambda)$ and called the Weyl module of $\lambda$. The Weyl-module $V(\lambda)$ has a unique irreducible quotient $L(\lambda)$, and any irreducible representation of $G$ is isomorphic to a $V(\lambda)$ for a unique dominant weight $\lambda$. If $k$ has characteristic $0$, we have $V(\lambda)=L(\lambda)$, but not in general. Nevertheless, the dimension of $V(\lambda)$ is independent of the characteristic of $k$, and given by the classical ``Weyl dimension formula''. \ps\ps

The following lemma is well-known in characteristic $0$, by lack of a reference we provide a proof in arbitrary characteristic. \ps\ps
\ps
\begin{lemme} \label{lemmerepg2} Let $k$ be an algebraically closed field and $G$ a connected linear algebraic group over $k$ of which is simple of type ${\bf G}_2$. Let $\omega_1$ {\rm (}resp. $\omega_2${\rm )} be the fundamental weight of $G$ which is a short (resp. long) root. Let $V$ be an irreducible, nontrivial, linear representation of $G$ of dimension $\leq 11$. Then we have: \ps \ps \begin{itemize}
\item[(i)] {\rm (}${\rm char}\, k \neq 2,3${\rm )}\, $\dim V=7$ and $V \simeq V(\omega_1)  = L(\omega_1)$, \ps\ps
\item[(ii)] {\rm (}${\rm char}\, k = 3${\rm )} \,$\dim V=7$ and we have either $V \simeq V(\omega_1)=L(\omega_1)$ or $V \simeq  L(\omega_2)$, and in the latter case $G\, \rightarrow {\rm GL}(V)$ is not a closed immersion. \ps\ps
\item[(iii)] {\rm (}${\rm char}\, k = 2${\rm )} \,$\dim V=6$ and $V \simeq L(\omega_1)$. \ps\ps
\end{itemize}
Moreover  there is up to a scalar a unique quadratic form on the $7$-dimensional vector space $V(\omega_1)$ which is preserved by $G$ and such that the resulting quadratic space over $k$ is regular.
\end{lemme}

\ps\ps

\begin{pf} 
Denote by $\alpha_1$ (resp. $\alpha_2$) the simple short (resp. long) root of $G$ and by $W$ its Weyl group. We have $\langle \omega_i, \alpha_j^\vee \rangle = \delta_{i,j}$, $\omega_1 = 2\alpha_1 +  \alpha_2$ and $\omega_2 = 3 \alpha_1 +2 \alpha_2$. Let $X$ be the set of nonzero dominant weights $\mu$ of $G$ such that $|W \cdot \mu|<|W|=12$. Observe that we have $X= \coprod_{i=1,2} \Z_{>0} \omega_i$ and $|W \cdot \mu| =6$ for each $\mu \in X$. If $\mu$ is a weight of $V$, then so is $w \cdot \mu$ for all $w \in W$; as a consequence, if $\lambda$ denotes the dominant weight of $G$ such that $V \simeq V(\lambda)$, then the dominant weights of $V$ are in $\{\lambda,0\}$. If $k$ has characteristic zero, then as a general fact, a dominant weight $\mu$ occurs in $V(\lambda)$ if, and only if, $\mu \prec \lambda$, i.e. $\lambda-\mu$ is a sum of positive roots. As $0 \prec \omega_1 \prec \omega_2$, this implies $\lambda = \omega_1$. We have $\dim V(\omega_1)=7$. \ps \ps
	Assume from now on that $k$ has prime characteristic $p$. The analysis above and Steinberg tensor product theorem show that $\lambda$ is $p$-restricted (as $\dim V(\lambda) < 6^2 =36$): we have $\langle \lambda, \alpha^\vee \rangle <p$ for each simple root $\alpha$. But then \cite[Chap. 2, Prop. 2.11]{rags} shows that $\lambda - n \alpha$ is a weight of $V(\lambda)$ for each simple root $\alpha$ and each integer $n$ such that $0 \leq n \leq \langle \lambda, \alpha^\vee \rangle$. The relations $2 \omega_1 = \omega_2 + \alpha_1$ and $2 \omega_2 = 3 \omega_1 + \alpha_2$ implies thus $\lambda = \omega_i$ with $i=1,2$. By \cite[Cor. 4.3 \& Table 1]{springerweyl}, we have $V(\omega_i)=L(\omega_i)$ if $p>3$. As $\dim V(\omega_2)=12$ (in fact, $V(\omega_2)$ is the adjoint representation of $G$) this excludes the case $i=2$ if ${\rm char} \,p\, >3$, hence prove assertion (i). \ps\ps
	If $p=3$, the same reference shows $L(\omega_1)=V(\omega_1)$, $\dim L(\omega_2)=7$, that $V(\omega_2)$ is an extension of $L(\omega_2)$ by $L(\omega_1)$, and that $L(\omega_2)$ is isomorphic to the restriction of $L(\omega_1)$ via the inseparable isogeny $G \rightarrow G$. If $p=2$, then $V(\omega_2)$ (which is the adjoint representation of $G$) is actually irreducible, hence we must have $V \simeq L(\omega_1)$. This shows (ii) and (iii). \ps\ps
	We now prove the last assertion. View $G$ as a ${\rm G}_2$-subgroup of ${\rm SO}(P)$, with $P$ regular of dimension $7$ over $k$. We have $P \simeq V(\omega_1)$ by \cite[\S 2.3]{springer}; these authors prove in particular that $P$ is irreducible, except in characteristic $2$ in which case the only invariant subspace is the kernel of $\beta_P$: see Thm. 2.3.3 {\it loc. cit.} The last statement follows in a standard way from this irreducibility property if ${\rm char} \,k\, \neq 2$. \ps\ps
	
	Assume ${\rm char} \,k\, =\, 2$, set $P=V(\lambda_1)$, and let $q_1,q_2 : P \rightarrow k$ be two $G$-invariant quadratic forms on $P$ which both gives the vector space $P$ a structure of regular quadratic space, say $P_1$ and $P_2$. Then the unique $G$-invariant line $L$ in $P$ must be ${\rm ker} \,\beta_{P_1} \,=\, {\rm ker} \beta_{P_2}$. Replacing $q_2$ by some scalar multiple if necessary we may assume that $q_1$ and $q_2$ coincide on this line. Then $q_1-q_2$ factors through a $G$-invariant quadratic form on $P/L$. As $P/L$ is an irreducible $G$-module, the kernel of $\beta_{P_1}-\beta_{P_2}$ is either $L$ or $P$. Assume it is $L$. Then $P/L$ is a nondegenerate quadratic space for $q_1-q_2$. As $G$ is simply connected, we thus have an injective algebraic group morphism $G \rightarrow {\rm Spin}(P/L) \simeq {\rm Spin}_6(k)$. But the halfspinor modules of ${\rm H}(k^3)$ have dimension $4 <6$, hence are trivial as $G$-modules by (iii), a contradiction. It follows that $\beta_{P_1}=\beta_{P_2}$, i.e. $q_1-q_2$ is a Frobenius semilinear form. As $P$ has no $G$-invariant subspace of dimension $6$, and as $k$ is perfect, this must be the $0$ form, and we are done. \end{pf}

We now prove $(ii) \Rightarrow (i)$ of Proposition \ref{algcarG2}. Let $G$ be a simple algebraic group of type ${\bf G}_2$.  Let $\rho : G \rightarrow {\rm SO}(P)$ be algebraic group homomorphism that is furthermore a closed immersion. If we can prove that the representation $P$ of $G$ is isomorphic to $V(\omega_1)$, then the last assertion of Lemma \ref{lemmerepg2} shows that $\rho$ is unique up to conjugation, hence that $\rho(G)$ is a ${\rm G}_2$-subgroup by Proposition \ref{carG2sbgp}. By Lemma \ref{lemmerepg2} (i) and (ii), we are done if ${\rm char} \, k \neq 2$. If ${\rm char}\, k \, =2$, then Lemma \ref{lemmerepg2} (iii) shows that $\omega_1$ is a highest weight of $P$, hence we have a nonzero morphism $f: V(\omega_1) \rightarrow P$. If $f$ is not an isomorphism, then $P$ has a subrepresentation $P'$ isomorphic to $L(\omega_1)$, and $P = P' \bot \ker \beta_P$. The space $P'$ is thus nondegenerate, but this is absurd as by last paragraph of the proof of Lemma \ref{lemmerepg2}, there is no nontrivial algebraic group homomorphism $G \rightarrow {\rm SO}_6(k)$. \ps \ps

Assume now that $G$ is a closed subgroup of ${\rm Spin}(P)$. Consider the central isogeny $\pi_P : {\rm Spin}(P) \rightarrow {\rm SO}(P)$. Then $G'=\pi_P(G)$ is a ${\rm G}_2$-subgroup by the previous paragraph.  We have seen that the morphism $\pi_P : {\rm Spin}(P) \rightarrow {\rm SO}(P)$ admits a section $j$ over such a subgroup $G'$, so that we have a direct product $\pi_P^{-1}(G') = \mu_2(k) \times j(G')$. But we have an inclusion $G \subset \pi_P^{-1}(G')$, and the projection of $G$ on the finite $\mu_2(k)$ factor is necessarily trivial as $G$ is simple, so we have $G \subset j(G')$ and then $G=j(G')$.
\end{pf}
}
\ps\ps
\subsection{Polynomials of type ${\bf G}_2$}
\ps\ps
\begin{definition}\label{defpg2} Let $k$ be a field and $P(t) \in k[t]$. We say that the polynomial $P$ is {\it of type $G_2$} if there are elements $x,y$ in an algebraic closure of $k$ such that 
$$P(t) = (t-1)(t-x)(t-y)(t-xy)(t-x^{-1})(t-y^{-1})(t-x^{-1}y^{-1}).$$
In particular, $P$ is a monic polynomial and $t^7 P(1/t) = -P(t)$.
\end{definition}
\ps\ps
\begin{prop}\label{typeg2ss} Let $P$ be a regular quadratic space of dimension $7$ over $k$ and $\gamma \in {\rm SO}(P)$. If $\gamma$ belongs to a ${\rm G}_2$-subgroup of ${\rm SO}(P)$, then the characteristic polynomial of $\gamma$ is of type ${\rm G}_2$.  The converse holds if $k$ is algebraically closed and $\gamma$ is semisimple.
\end{prop}
\ps\ps
\begin{pf} We may assume that $k$ is algebraically closed and, by Jordan decomposition, that $\gamma$ is semisimple. In this case, we have recalled in the proof of Proposition \ref{algcarG2} that the linear representation $P$ of a ${\rm G}_2$-subgroup $G$ of ${\rm SO}(P)$ is isomorphic to the Weyl module $V(\omega_1)$. In particular, its weights are $0$ and the six short roots of $G$, which are $\pm \alpha_1, \pm (\alpha_2+\alpha_1), \pm (\alpha_2+2\alpha_1)$ in the notations of the proof {\it loc. cit}, which shows the first assertion. The last assertion is then an easy consequence of the fact that the conjugacy class of $\gamma$ in ${\rm SO}(P)$ is uniquely determined by $\det(t-\gamma)$ (a simple special case of Corollary \ref{corwindex}).
\end{pf}
\ps\ps
\begin{remark}\label{rempotg2}{\rm  Let $P \in k[t]$ be a monic polynomial such that $t^7P(1/t)=-P(t)$. There is a unique monic polynomial $Q \in k[t]$ of degree $3$ such that $t^{-3}P(t) = (t-1)Q(t+t^{-1})$. Write $Q = t^3 - at^2 + bt - c$ with $a,b,c \in k$ (some simple universal integral polynomials in the coefficients of $P$).  It is not difficult to check, {\it e.g.} using the equivalence between (i) and (iii) in Lemma \ref{calcpcar} below, that $P$ is of type ${\bf G}_2$ if, and only if, we have the relation $a^2 = 2b+c+4$. }
\end{remark}

\section{Proof of Theorems \ref{thmintroa} \& \ref{thmintrob} }\label{pfspin7}

The following lemma is well-known (at least in characteristic $0$).

\begin{lemme}\label{pcarspin} Let $k$ be an algebraically closed field, $V$ a regular quadratic space of odd dimension $2r+1$ over $k$, $(W,\rho)$ a spin representation of ${\rm GSpin}(V)$ and $\pi_V : {\rm GSpin}(V) \rightarrow {\rm SO}(V)$ the canonical map. For any element $\gamma \in {\rm GSpin}(V)$, there are elements $x_0, x_1,\dots,x_r \in k^\times$ such that: \begin{itemize} \ps \ps
\item[(i)] $\det(t - \rho(\gamma)) = \prod_{(\epsilon_i) \in \{-1, 1\}^r} (t -  x_0 \prod_{i=1}^r x_i^{\epsilon_i})$,\ps\ps
\item[(ii)] $\det(t - \pi_V(\gamma)) = (t-1) \prod_{i=1}^r (t-x_i^2)(t-x_i^{-2})$, \ps\ps
\item[(iii)] $\nu_V(\gamma) = x_0^2$.
\end{itemize}
\end{lemme}

\begin{pf} Both ${\rm GSpin}(V)$ and ${\rm SO}(V)$ are connected linear algebraic $k$-groups in a natural way \cite[p. 40]{springer}, and both $\pi_V$ and $\rho$ are $k$-group homomorphisms. By considering a Jordan decomposition of $\gamma$, we may assume that $\gamma$ is semisimple, in which case so is $\pi_V(\gamma)$. 
We may thus assume that $V = {\rm H}(I) \bot k e$ and that $\pi_V(\gamma)$ preserves $I$, $I^\ast$ and fixes $e$. Choose $x_1,\dots,x_r \in k^\times$ whose squares are the eigenvalues of $\pi_V(\gamma)$ on $I$, then $\det(t - \pi_V(\gamma))$ is as in assertion (ii). The last exact sequence of \S \ref{parclag} shows that we may write $\gamma= \rho_{{\rm H}(I);ke}(\gamma' \times 1)$ with $\gamma' \in {\rm GSpin}({\rm H}(I))$, and Proposition \ref{decospin} implies that $\det(t - \rho(\gamma))$ is the characteristic polynomial of $\gamma'$ in a spinor module for ${\rm H}(I)$. Proposition \ref{cliffordhyp} shows that there is a unique element $\lambda \in k^\times$ such that we have $\gamma' = \lambda \cdot \widetilde{\rho}_I(\pi_{{\rm H}(I)}(\gamma'))$ and $\det(t-\rho(\gamma))=\det(t-\lambda \,\pi_V(\gamma) \, | \, \Lambda\, I)$. This concludes the proof with $x_0 = \lambda x_1\cdots x_r$. \end{pf}
\ps\ps

Polynomials {\it of type ${\bf G}_2$} have been introduced in Definition \ref{defpg2}. \ps\ps

\begin{lemme}\label{calcpcar} Let $k$ be a field, $V$ a regular quadratic space of dimension $7$ with trivial even Clifford algebra, $(W,\rho)$ a spin representation of ${\rm GSpin}(V)$, $\pi_V : {\rm GSpin}(V) \rightarrow {\rm SO}(V)$ the natural map, and $\gamma \in {\rm GSpin}(V)$.  The following properties are equivalent: \ps\begin{itemize}\ps\ps
\item[(i)] the characteristic polynomial of $\pi_V(\gamma)$ is of type ${\bf G}_2$, \ps\ps
\item[(ii)] $\rho(\gamma)$ has an eigenvalue whose square is $\nu_V(\gamma)$, \ps\ps
\item[(iii)] $\nu_V(\gamma)$ is an eigenvalue of $\rho(\gamma^2)$, \ps\ps
\item[(iv)] $\det(t - \nu_V(\gamma)^{-1} \rho(\gamma^2)) = (t-1) \det (t- \pi_V(\gamma^2))$. \ps \ps
\end{itemize}
Moreover, if $\rho(\gamma) \in \GL(W)$ has an eigenvalue $\lambda$ such that $\lambda^2 = \nu_V(\gamma)$, then we have  
$\det(t - \lambda^{-1} \rho(\gamma)) = (t-1) \det (t- \pi_V(\gamma))$.
\end{lemme}
\ps\ps
\begin{pf} Note that the equivalence (ii) $\Leftrightarrow$ (iii) and the implication (iv) $\Rightarrow$ (iii) are obvious. Let $K$ be an algebraic closure of $k$. By Lemma \ref{pcarspin} applied to $\gamma$ viewed as an element of ${\rm GSpin}(V \otimes_k K)$, we may fix $x_0,x_1,x_2,x_3$ in $K^\times$ as in that lemma. \ps\ps

Assume assertion (i) of Lemma~\ref{calcpcar} holds.  Up to replacing some of the $x_i$, with $i\geq 1$, by $x_i^{-1}$ if necessary,  we may assume that we have $x_3^2 = x_2^2 x_1^2$. In particular, if we set $\epsilon = x_1x_2x_3^{-1}$, then $\epsilon=\pm 1$ and $\epsilon x_0$ is an eigenvalue of $\rho(\gamma)$, hence assertion (ii) above holds.\ps\ps
 Conversely, assume that $\epsilon \in \{\pm 1\}$ and that $\epsilon x_0$ is an eigenvalue of $\rho(\gamma)$. Again, up to replacing some of the $x_i$, with $i\geq 1$, by $x_i^{-1}$ if necessary, we may assume that we have $x_1x_2x_3^{-1}=\epsilon$. In particular, we have $x_3^2 = x_2^2x_1^2$ and assertion (i) holds. This shows (ii) $\Rightarrow$ (i). Moreover, the relation $\epsilon x_1x_2x_3^{-1} = 1$ implies $\epsilon x_1x_2x_3 = x_1^2x_2^2$, $\epsilon x_1x_2^{-1}x_3 = x_1^2$ and $\epsilon x_1x_2^{-1}x_3^{-1} =  x_2^2$, hence the equality 
$$\det(t - x_0^{-1} \epsilon \rho(\gamma)) = (t-1) \det (t- \pi_V(\gamma)),$$
holds, which proves the last statement. This statement applied to $\gamma^2$ shows in turn (iii) $\Rightarrow$ (iv), which concludes the proof.\end{pf}
\ps\ps
\begin{thm} Let $V$ be a regular quadratic space of dimension $7$ over $k$ with trivial even Clifford algebra and $(W,\rho)$ be a spin representation of ${\rm Spin}(V)$. Let $\Gamma \subset {\rm Spin}(V)$ be a subgroup. Assume that: \ps \ps
\begin{itemize}
\item[(i)] $W$ is semisimple as a $k[\Gamma]$-module, \ps\ps
\item[(ii)] for each $\gamma \in \Gamma$, the element $\rho(\gamma)$ has the eigenvalue $1$, \ps\ps
\end{itemize}
then there is an element $w \in W-\{0\}$ such that $\rho(\gamma) w = w$ for all $\gamma \in W$. 
\end{thm}
\ps\ps
\begin{pf} Consider $E_1=W$ and $E_2=V \oplus k$ as $k$-linear representations of $\Gamma$. Assumption (ii) and the last assertion of Lemma~\ref{pcarspin} show that for each $\gamma \in \Gamma$, we have $\det ( t - \gamma\, |\,  E_1) = \det (t - \gamma \,|\, E_2)$. The Brauer-Nesbitt theorem\footnote{The Brauer-Nesbitt theorem is the following statement \cite[\S 20, {\rm no} 6, Thm. 2 \& Cor. 1]{bou8}. Let $k$ be a field, $\Gamma$ a group, and $E_i$ a finite dimensional, semisimple, $k$-linear representation of $\Gamma$ for $i=1,2$. Assume we have $\det ( t - \gamma\, |\,  E_1) = \det (t - \gamma \,|\, E_2)$ for each $\gamma \in \Gamma$. Then $E_1$ is isomorphic to $E_2$. } implies thus that $E_1$ and $E_2$ have isomorphic semisimplification. The theorem follows as $E_1$ is semisimple by assumption, and $E_2$ contains the trivial representation. 
\end{pf}
\ps\ps
\begin{thm}\label{thmbsecondev} Let $C$ be an octonion $k$-algebra, $P \subset C$ the quadratic space of pure octonions of $C$, $Q \subset C$ the subset of element $v \in C$ such that ${\rm q}(v) \neq 0$, and $\widetilde{{\rm j}}_C: {\rm G}_2(C) \rightarrow {\rm Spin}(P)$ the canonical morphism. Recall that the quadratic space $C$ itself may be endowed with a canonical structure of a spinor module for $P$ and denote by $\rho_C : {\rm GSpin}(P) \rightarrow {\rm GSO}(C)$ the associated spin representation. \ps \ps

Let $\Gamma \subset {\rm Spin}(P)$ be a subgroup such that $\rho_C(\Gamma)$ acts in a semisimple way on the $k$-vector space $C$. The following assertions are equivalent: \ps\ps
\begin{itemize}
\item[(a)] for each $\gamma \in \Gamma$, there is $g \in {\rm GSpin}(P)$ such that $g \gamma g^{-1} \in \widetilde{{\rm j}}_C( {\rm G}_2(C) )$, \ps\ps
\item[(b)] for each $\gamma \in \Gamma$, there is $s \in Q$ such that $\rho_C(\gamma)s = s$, \ps\ps
\item[(c)] for each $\gamma \in \Gamma$, the element $\rho_C(\gamma) \in \GL(C)$ has the eigenvalue $1$, \ps\ps
\item[(d)] there is an element $g \in {\rm GSpin}(P)$ such that $g\, \Gamma\, g^{-1} \subset  \widetilde{{\rm j}}_C( {\rm G}_2(C) )$, \ps\ps
\item[(e)] there is an element $s \in Q$ such that $\rho_C(\gamma)s=s$ for all $\gamma \in \Gamma$.
\end{itemize}
\end{thm}

\begin{pf} The equivalences $(a) \Leftrightarrow (b)$ and $(d) \Leftrightarrow (e)$ follow from Proposition \ref{propgspin7}, and the implications $(e) \Rightarrow (b) \Rightarrow (c)$  are obvious. Assume that assertion (c) holds. The previous theorem ensures that the trivial representation of $\Gamma$ occurs in the semisimple spinor module $C$.  To conclude the proof it suffices to show that for any group $G$ and any nondegenerate quadratic space $U$, if we have a group homomorphism $G \rightarrow {\rm O}(U)$ whose underlying $k$-linear representation $U$ is semisimple and contains the trivial representation, then there exists a $G$-invariant vector in $U$ which is nonisotropic. But as $U$ is semisimple, we have $U \simeq U^G \oplus W$, where $-^G$ denotes $G$-invariants and $W$ is some $G$-stable subspace of $V$. As the trivial representation is selfdual and $W$ is semisimple, we have $(W^\ast)^G = W^G = 0$. It follows that $U^G$ is a nondegenerate subspace of $U$. We conclude as any nonzero nondegenerate quadratic space obviously contains some nonisotropic vector. 
\end{pf}

Given the definition of ${\rm G}_2$-subgroups in \S \ref{defg2sb}, we get the following:

\begin{cor}\label{corthmbsecondev} Let $P$ be a regular quadratic space of dimension $7$ over a field $k$ such that ${\rm Spin}(P)$ possesses ${\rm G}_2$-subgroups. Let $\Gamma \subset {\rm Spin}(P)$ be a subgroup acting in a semisimple way on a spinor module $W$ for $P$. The following assertions are equivalent: \ps\ps
\begin{itemize}
\item[(i)] each element of $\Gamma$ belongs to some ${\rm G}_2$-subgroup of ${\rm Spin}(P)$,\ps\ps
\item[(ii)] each element of $\Gamma$ has the eigenvalue $1$ in $W$, \ps\ps
\item[(iii)] there is a ${\rm G}_2$-subgroup of ${\rm Spin}(P)$ containing $\Gamma$. \ps\ps
\end{itemize}
\end{cor}

\begin{cor}\label{corcar2} Let $k$ be a perfect field of characteristic $2$ and $P$ a regular quadratic space of dimension $7$ over $k$ such that ${\rm SO}(P)$ possesses ${\rm G}_2$-subgroups. Let $\Gamma \subset {\rm SO}(P)$ whose inverse image in ${\rm Spin}(P)$ acts in a semisimple way on a spinor module for $P$. The following assertions are equivalent: \ps\ps
\begin{itemize}
\item[(i)] each element of $\Gamma$ belongs to some ${\rm G}_2$-subgroup of ${\rm SO}(P)$,\ps\ps
\item[(ii)] for each $\gamma \in \Gamma$, then $\det(t-\pi_P(\gamma))$ is of type ${\bf G}_2$, \ps\ps
\item[(iii)] there is a ${\rm G}_2$-subgroup of ${\rm SO}(P)$ containing $\Gamma$. \ps\ps
\end{itemize}
\end{cor}

\begin{pf} The map $\pi_P : {\rm Spin}(P) \rightarrow {\rm SO}(P)$ is surjective by the exact sequence \eqref{sexspin}. Moreover, an element $\gamma$ of ${\rm Spin}(P)$ possesses the eigenvalue $1=\pm 1$ in a spinor module for $P$ if, and only if,  $\det(t-\pi_P(\gamma))$ is of type ${\bf G}_2$, by Lemma \ref{calcpcar}. The corollary follows thus from Corollary~\ref{corthmbsecondev} applied to $\pi_P^{-1}(\Gamma)$.
\end{pf}

\begin{cor} \label{corconnected} Let $k$ be an algebraically closed field,  $P$ be a regular quadratic space over $k$, and $\Gamma \subset {\rm SO}(P)$ a subgroup whose inverse image in ${\rm Spin}(P)$ acts in a semisimple way on a spinor module for $P$. Assume that: \begin{itemize}\ps
\item[(i)] the characteristic polynomial of each element of $\Gamma$ is of type ${\bf G}_2$, \ps\ps
\item[(ii)] $\Gamma$ is a Zariski closed and connected subgroup of ${\rm SO}(P)$. \ps\ps
\end{itemize}
Then there is a ${\rm G}_2$-subgroup of ${\rm SO}(P)$ containing $\Gamma$. 
\end{cor}

\begin{pf} Let $G$ denote the connected component of the identity of the inverse image of $\Gamma$ in ${\rm Spin}(P)$. This is a connected linear algebraic subgroup of ${\rm Spin}(P)$ such that $\pi_P(G)=\Gamma$, as $\Gamma$ is irreducible as a variety. It is enough to show that $G$ is contained in a ${\rm G}_2$-subgroup of ${\rm Spin}(P)$. By Lemma \ref{calcpcar}, for each element $g \in G$, the element $g^2$ has the eigenvalue $1$, acting on a spinor module $W$ for $P$. But as $G$ is connected, any semisimple element of $G$ is a square, as tori are divisible as groups. As $G$ has a Zariski-dense subset consisting of semisimple elements, it follows that $\det(1-g | W)=0$ for all $g \in G$. We conclude by Corollary \ref{corthmbsecondev}. 
\end{pf}

In the following corollary, we denote by $-1 \in {\rm Spin}(E)$ the generator of the natural central subgroup $\mu_2(k)$. \ps

\begin{cor} \label{carbetaspin}  Let $k$ be an algebraically closed field,  $P$ be a regular quadratic space over $k$, and $\Gamma \subset {\rm SO}(P)$ a subgroup whose inverse image $\widetilde{\Gamma}$ in ${\rm Spin}(P)$ acts in a semisimple way on a spinor module $W$ for $P$. Then $\Gamma$ is contained in a ${\rm G}_2$-subgroup of ${\rm SO}(P)$ if, and only if, there is a group homomorphism $\beta  : \widetilde{\Gamma} \rightarrow k^\times$ which occurs in the restriction of $W$ to $\widetilde{\Gamma}$ and satisfies $\beta^2=1$. 
\end{cor}

\begin{pf} Under the assumption on $k$, there exists ${\rm G}_2$-subgroups and we have $\pi_P({\rm Spin}(P))={\rm SO}(P)$. We already explained in the last paragraph of the proof of Proposition \ref{algcarG2} that the inverse image in ${\rm Spin}(P)$ of a ${\rm G}_2$-subgroup of ${\rm SO}(P)$ has the form $\mu_2(k) \times H$ where $H$ is a ${\rm G}_2$-subgroup of ${\rm Spin}(P)$. Note that the generator $-1$ of $\mu_2(k)$ acts by $-{\rm id}$ on $W$.
If $\Gamma$ is contained in a ${\rm G}_2$-subgroup of ${\rm SO}(P)$, the existence of the character $\beta$ as in the statement follows {\it e.g.} from Proposition \ref{propgspin7}. Conversely, if there is a character $\beta$ as in the statement, and if $\Gamma^1 \subset \widetilde{\Gamma}$ denotes the kernel of $\beta$, then we have $\widetilde{\Gamma} = \mu_2(k) \times \Gamma^1$ and Corollary \ref{corthmbsecondev} applied to $\Gamma^1$ concludes the proof. 
\end{pf}

\begin{remark} \label{examplethmAR} {\rm Let $k=\R$ as in Example \ref{examplecayley}. Then any subgroup $\Gamma \subset {\rm Spin}(7)$ acts in a semisimple way on the Euclidean space $C$ (see e.g. the proof of Proposition \ref{critasss} (i) below). Corollary~\ref{corthmbsecondev} thus asserts that if an arbitrary subgroup $\Gamma \subset {\rm Spin}(7)$ has the property that each element $\gamma \in \Gamma$ fixes a point of the unit sphere ${\mathbb S}^7$ (or equivalently, belongs to some ${\rm G}_2$-subgroup of ${\rm Spin}(7)$), then $\Gamma$ itself has a fixed point in this sphere (or equivalently, is included in some ${\rm G}_2$-subgroup of ${\rm Spin}(7)$). Moreover, if $\Gamma' \subset {\rm SO}(7)$ is a compact connected subgroup, and if the characteristic polynomials of the elements of $\Gamma'$ are all of type ${\bf G}_2$, then Corollary \ref{corconnected} shows that $\Gamma'$ is contained in a ${\rm G}_2$-subgroup of ${\rm SO}(7)$.}
\end{remark}

\section{Proof of Theorem \ref{thmintrod} }\label{proofthmb}

\subsection{Basic identities and the irreducible case}\label{basicid} Let $k$ be a field and $\Gamma$ a group. We shall denote by ${\rm R}_k(\Gamma)$ the Grothen\-dieck ring of finite dimensional $k$-linear representations of $\Gamma$. If $X$ is a finite dimensional $k$-linear representation of $\Gamma$, we denote by $[X]$ its class in ${\rm R}_k(\Gamma)$. The Adams endomorphism $\psi^2$ of the additive group of ${\rm R}_k(\Gamma)$ is defined by $\psi^2 [X] = [X \otimes X] - 2 [\Lambda^2 X] = [{\rm Sym}^2 X] - [\Lambda^2 X]$ for all $X$ as above. We denote by $1$ the trivial representation, as well as its class in ${\rm R}_k(\Gamma)$. \ps\ps

\begin{lemme}\label{lemmesymlambda2} Let $E$ be a $7$-dimensional regular quadratic space over a field $k$, with trivial even Clifford algebra, and let $W$ a spin representation of ${\rm GSpin}(E)$. The representation $W\otimes W \otimes \nu_{E}^{-1}$ of ${\rm GSpin}(E)$ factors via $\pi_E$ through a representation of ${\rm SO}(E)$ and we have 
$$[({\rm Sym}^2 \,W) \otimes \nu_E^{-1}] = [\Lambda^3 \,E] + 1 \hspace{.5 cm} {\rm and} \hspace{.5 cm} [(\Lambda^2\, W) \otimes \nu_E^{-1}] = [\Lambda^2 (1 \oplus E)]$$ in ${\rm R}_k({\rm SO}(E))$, as well as the following relation in ${\rm R}_k({\rm GSpin}(E))$: 
$$\psi^2 [W] \otimes \nu_E^{-1}- \psi^2 (1+[E]) = [\Lambda^3\, E] - [E] - [{\rm Sym}^2 \,E].$$ 
If furthermore ${\rm char}\, k \neq 2$ then we have $k[{\rm GSpin}(E)]$-module isomorphisms $({\rm Sym}^2 \,W) \otimes \nu_E^{-1} \simeq 1 \oplus \Lambda^3\, E$ and $(\Lambda^2 \,W) \otimes \nu_E^{-1} \simeq  \Lambda^2 (1 \oplus E)$. 
\end{lemme}

\begin{pf} The second identity is a consequence of the first two ones and of the natural isomorphisms $\Lambda^2 (E \oplus k )\simeq \Lambda^2 E \oplus E$ and ${\rm Sym}^2 (E \oplus k) \simeq {\rm Sym}^2 E \oplus E$. To check the first two displayed equalities, it is enough to check that the characteristic polynomial of any element of ${\rm GSpin}(E)$ is the same on both representations in the two cases, by the Brauer-Nesbitt theorem. We omit the details, but this follows at once from the formulas given in Lemma~\ref{pcarspin} (see Remark \ref{decendc} for another argument). \ps\ps
To prove the last assertion, we may and do assume that $k$ is algebraically closed. Assume ${\rm char}\, k \neq 2$. It is well-known that the $k[{\rm SO}(E)]$-modules $\Lambda^i\, E$ are then irreducible for $0 \leq i \leq 3$ (see e.g. the tables in \cite{lubeck}). Recall from \S \ref{octoalg} that $W$ may be given a structure of a nondegenerate quadratic space over $k$ such that ${\rm GSpin}(E)$ acts by orthogonal similitudes of factor $\nu_{E}$. In particular, we have an isomorphism $W \simeq W^\ast\otimes \nu_{E}$, and both representations $({\rm Sym}^2 \,W)\otimes \nu_{E}^{-1}$ and $(\Lambda^2\, W)\otimes \nu_{E}^{-1}$ of are selfdual (for the first one it also uses $2 \in k^\ast$). As $\Lambda^i \, E$ is irreducible and selfdual for $0 \leq i \leq 3$, the (proven) first two displayed equalities of the statement imply
$({\rm Sym}^2 \,W)\otimes \nu_{E}^{-1} \simeq 1 \oplus \Lambda^3 \, E$ and $(\Lambda^2 \,W)\otimes \nu_{E}^{-1} \simeq E \oplus \Lambda^2\,E$. 
 \end{pf}

\begin{remark}\label{decendc}{\rm  Let $E$ be a regular quadratic space of odd dimension with trivial even Clifford algebra, and $W$ a spinor module for $E$. Then ${\rm GSpin}(E)$ acts by conjugation on ${\rm C}(E)_0 \simeq {\rm End}\, W$. This action factors through its quotient ${\rm SO}(E)$ and is induced by the obvious action of the latter on ${\rm C}(E)$ (defined by functoriality of the Clifford construction). This action preserves the canonical filtration of ${\rm C}(E)$ whose associated graded $k$-algebra is the exterior algebra $\Lambda \, E$. In particular, we obtain a more conceptual explanation of the identity $[{\rm End}\, W]= \sum_{0 \leq 2i < \dim E} [\Lambda^{2i}\, E]$ in ${\rm R}_k({\rm SO}(E))$. }\end{remark}
\ps\ps
\begin{prop}\label{propfondso7g2} Let $E$ be a $7$-dimensional regular quadratic space over a field $k$, $K$ an algebraic closure of $k$, and $(W,\rho)$ a spin representation of ${\rm GSpin}(E \otimes_k K)$. Let $\Gamma \subset {\rm SO}(E)$ be a subgroup and set $\widetilde{\Gamma}:=\pi_E^{-1}(\Gamma) \subset {\rm GSpin}(E)$. The following properties are equivalent:\ps \begin{itemize}
\item[(i)] for all $\gamma \in \Gamma$, the polynomial $\det (t - \gamma \,|\, E)$ is of type ${\bf G}_2$, \ps \ps
\item[(ii)] for all $\gamma \in \widetilde{\Gamma}$, $\rho(\gamma)$ has an eigenvalue whose square is $\nu(\gamma)$, \ps\ps
\item[(iii)] the equality $\psi^2 [W] \otimes \nu_E^{-1} = 1 + \psi^2 [E \otimes_k K]$ holds in ${\rm R}_K(\widetilde{\Gamma})$,\ps\ps
\item[(iv)] the equality $[\Lambda^3 E ] = [E] + [{\rm Sym}^2 E]$ holds in ${\rm R}_k(\Gamma)$. \ps\ps
\end{itemize}
In particular, if $\Lambda^3 E$ is semisimple as a representation of $\Gamma$, and if these properties hold, then $E$ admits a nonzero, $\Gamma$-invariant, alternating trilinear form.
\end{prop}

\begin{pf} The equivalences $(i) \Leftrightarrow (ii) \Leftrightarrow (iii)$ are exactly the equivalences $(i) \Leftrightarrow (ii) \Leftrightarrow (iv)$ of Lemma~\ref{calcpcar} (plus the Brauer-Nesbitt theorem). The equivalence between (iii) and (iv) follows from Lemma \ref{lemmesymlambda2}. The last assertion comes from the fact that $({\rm Sym}^2 E)^\ast$ contains the nonzero $\Gamma$-invariant $\beta_E$ and from the relations  in ${\rm R}_k(\Gamma)$ : $$[({\rm Sym}^2 E)^\ast]=[{\rm Sym}^2 E^\ast]=[{\rm Sym}^2\, E].$$\end{pf}

 \begin{thm} \label{thmbirrcase} Let $E$ be a $7$-dimensional regular quadratic space over an algebraically closed field $k$ and $\Gamma \subset {\rm SO}(E)$ a subgroup. Assume that the representation $E$ of $\Gamma$ is irreducible, with $\Lambda^3\, E$ and ${\rm Sym}^2\, E$ semisimple. If the characteristic polynomial of each element of $\Gamma$ is of type ${\bf G}_2$ then there is a unique ${\rm G}_2$-subgroup of ${\rm SO}(E)$ containing $\Gamma$.
\end{thm}

\begin{pf} Assume first that $E$ is an arbitrary $7$-dimensional vector space over $k$. The group $\GL(E)$ naturally acts on $\Lambda^3 E$; if $f \in  \Lambda^3 E$, we set $H_f = \{g \in {\rm SL}(E), g \cdot f = f\}$, which is a closed subgroup of ${\rm SL}(E)$. By \cite[Cor 2.3]{cohen}, if $H_f$ acts irreducibly on $E$ then, as an algebraic group, $H_f$ is simple of type ${\bf G}_2$. It follows from Lemma \ref{lemmerepg2} that in this case, $k$ has characteristic $\neq 2$, $E$ is isomorphic to the representation $V(\omega_1)$ of $H_f$, and there is a structure of regular quadratic space on $E$ preserved by $H_f$. \ps 

Let now $E$ and $\Gamma$ be as in the statement. By Proposition \ref{propfondso7g2}, there is a nonzero $\Gamma$-invariant element $f \in \Lambda^3 E$. In particular, $H_f \supset \Gamma$ acts irreducibly on $E$, hence the discussion above applies. We have thus ${\rm char} \,k\, \neq 2$ and there is a nondegenerate quadratic form on $E$ preserved by $H_f$. Such a form is preserved by $\Gamma$, which acts irreducibly on $E$, hence this form must be proportional to the quadratic form of $E$. In other words, we have $H_f \subset {\rm SO}(E)$ and $H_f$ is a ${\rm G}_2$-subgroup of ${\rm SO}(E)$ by Proposition \ref{algcarG2}. \ps

If $H \subset {\rm SO}(E)$ is a ${\rm G}_2$-subgroup containing $\Gamma$, then we have $H = H_{f'}$ for some nonzero $f' \in \Lambda^3 E$ since $H$ is conjugate to $H_f$. It is thus enough to show that the form $f$ chosen above is unique up to a scalar. Under the assumptions we have $\Lambda^3 \,E\, \simeq \,E \,\oplus\, {\rm Sym}^2\, E$ as $k[\Gamma]$-modules by Proposition \ref{propfondso7g2}. But as $E$ is irreducible and ${\rm char} \,k\,  \neq 2$, the subspace of $\Gamma$-invariants is trivial in $E$ and one dimensional in ${\rm Sym}^2 \,E$. 
\end{pf}
\ps\ps

\subsection{Some (counter-)examples: ${\rm O}_2^{\pm}(k)$ and $\GL_2(\Z/3\Z)$} \label{parexamples}\ps\ps
 As we shall see, the irreducibility assumption in Theorem \ref{thmbirrcase} is actually necessary. We introduce now two important examples of subgroups $\Gamma$ of ${\rm SO}(E)$ (with $\dim E = 7$) such that each element of $\Gamma$ has a characteristic polynomial of type ${\bf G}_2$, but which will turn out not to be included in any ${\rm G}_2$-subgroup of ${\rm SO}(E)$. 

\ps\ps\ps
{\sc The group ${\rm O}_2^{\pm}(k)$} \ps\ps

Let $k$ be an algebraically closed field of characteristic $\neq 2$. Consider the similitude group ${\rm GO}_2(k)$ of the quadratic space $P:={\rm H}(k)$ which is hyperbolic of dimension $2$ over $k$. Denote by $\mu : {\rm GO}_2(k) \rightarrow k^\times$ the similitude factor (so ${\rm ker}\, \mu = {\rm O}_2(k)$) and set $\epsilon = \det :  {\rm GO}_2(k) \rightarrow k^\times$. The character $\epsilon \mu^{-1}$ has order $2$ and its kernel is the subgroup ${\rm GSO}_2(k) \subset  {\rm GO}_2(k)$ of proper similitudes, which is also the stabilizer of each of the two isotropic lines in $P$. Fix such a line and denote by $\chi : {\rm GSO}_2(k) \rightarrow k^\times$ the character defined by the action on this line; then ${\rm GSO}_2(k)$ acts by $\mu \chi^{-1}$ on the other line and the morphism $\chi \times \mu : {\rm GSO}_2(k) \rightarrow k^\times \times k^\times$ is an isomorphism. We set $${\rm O}_2^\pm (k) = \{ g \in {\rm GO}_2(k), \, \, \, \mu^2(g) = 1\}={\rm ker}\, \mu^2 = {\rm ker}\, \epsilon^2.$$ 
The orthogonal group ${\rm O}_2(k)$ is an index $2$ subgroup of ${\rm O}_2^{\pm}(k)$, and as $-1$ is a square in $k$ we have an exact sequence $$ 1 \longrightarrow \mu_2(k) \overset{\lambda \mapsto (\lambda,\lambda \,{\rm id}_P)}{\longrightarrow} \mu_4(k) \times {\rm O}_2(k) \longrightarrow {\rm O}_2^\pm(k) \longrightarrow 1$$
(with $\mu_n(k)=\{\lambda \in k^\times, \lambda^n=1\}$). The three characters $\mu, \epsilon, \mu\epsilon$ are distinct of order $2$ when restricted to ${\rm O}_2^\pm (k)$ and, as a $k[{\rm O}_2^\pm (k)]$-module, the space $P$ is irreducible and satisfies $P^\ast \simeq P \otimes \epsilon \simeq P \otimes \mu$. \ps\ps

\begin{prop} \label{o2pm} Let $k$ be an algebraically closed field of characteristic $\neq 2$ and $E$ a nondegenerate quadratic space of dimension $7$ over $k$. There is a unique conjugacy class of injective group homomorphisms $$\rho : {\rm O}_2^\pm(k) \rightarrow {\rm SO}(E)$$ such that the $k[{\rm O}_2^\pm(k)]$-module $E$ is isomorphic to $P \oplus P^\ast \oplus \epsilon \oplus \mu \oplus \epsilon\mu$, where $P$, $\mu$ and $\epsilon$ are defined as above. For such a $\rho$, and any $\gamma \in {\rm O}_2^\pm(k)$, then $\det(t - \rho(\gamma))$ is of type ${\bf G}_2$. 
\end{prop}
\ps\ps
\begin{pf} The existence of such a morphism $\rho : {\rm O}_2^\pm(k) \rightarrow {\rm SO}(F)$ with $F = {\rm H}(P) \bot k^3$ and $k^3$ endowed with the standard quadratic form $(x_i) \mapsto \sum_i x_i^2$, is immediate from the equalities $\epsilon^2=\mu^2=1$. The existence part of the first assertion follows as we have $F \simeq E$ since $k$ is algebraically closed, as well as the uniqueness part (up to conjugation) by Corollary \ref{corwindex}. It remains to show that for any $\gamma \in {\rm O}_2^\pm(k)$ then $\det(t - \rho(\gamma))$ is of type ${\bf G}_2$. Set $H = {\rm O}_2^\pm(k) \cap {\rm GSO}_2(k) = \ker \epsilon\mu$ and choose $\chi : H \rightarrow k^\ast$ as in the discussion above. We have an isomorphism $$E \simeq \chi \oplus \mu  \chi^{-1} \oplus \chi^{-1} \oplus \mu^{-1} \chi \oplus \mu \oplus \mu \oplus 1$$
of $k[H]$-modules, which shows that $\det(t-\rho(\gamma))$ is of type ${\bf G}_2$ for all $\gamma \in H$ (recall $\mu^2=1$). Let $\gamma \in {\rm GO}_2(k)-{\rm GSO}_2(k)$, i.e. $\mu(\gamma) \neq \epsilon(\gamma)$. A simple computation shows that we have $\gamma^2 = \mu(\gamma) {\rm id}_P$ and $\det(t-\gamma)=t^2-\mu(\gamma)$. In particular, if we have $\gamma \in {\rm O}_2^\pm(k)$, then the element $s:=\mu(\gamma)$ is $\pm 1$, $\epsilon(\gamma)=-s$, and we have $$\det (t - \rho(\gamma)) = (t^2-s)^2(t-s)(t+s)(t+1)=(t^2\pm 1)^2(t+1)^2(t-1)$$
which is easily seen to be of type ${\bf G}_2$. \end{pf}
\ps
A consequence of our main result (Theorem \ref{mainthmb}) will be that if $\rho$ is as in the statement of the proposition above, then $\rho({\rm O}_2^\pm(k))$ is not contained in any ${\rm G}_2$-subgroup of ${\rm SO}(E)$. We shall now give another argument for this fact. \ps\ps

\begin{lemme} \label{lemmacentor2} Let $\rho$ be as in Proposition \ref{o2pm}. There is a unique conjugacy class of subgroups $\Gamma \subset {\rm O}_2^\pm(k)$ such that $\Gamma \simeq \Z/4\Z \times \Z/2\Z$ and such that the $k[\Gamma]$-module $E$ is the direct sum of all nontrivial characters of $\Gamma$.  These subgroups are exactly the centralizers of the order $2$ elements $s \in {\rm O}_2^\pm(k)$ such that $\mu(s)=1$ and $\epsilon(s)=-1$. 
\end{lemme}
\ps\ps
Indeed, this follows from an inspection of the centralizers of the order $2$ elements of ${\rm O}_2^\pm(k)$ that we leave as an exercise to the reader. \ps\ps

\begin{prop} Let $\Gamma$ be a subgroup of ${\rm O}_2^\pm(k)$ isomorphic to $\Z/4\Z \times \Z/2\Z$ and as in Lemma \ref{lemmacentor2}. Then $\rho(\Gamma)$ {\rm (}hence {\it a fortiori} $\rho({\rm O}_2^{\pm}(k))${\rm )} is not contained in any ${\rm G}_2$-subgroup of ${\rm SO}(E)$. 
\end{prop}
\ps

\begin{pf} (see Proposition \ref{lemmaindex2} for a rather different proof) Let $H$ be a ${\rm G}_2$-subgroup of ${\rm SO}(E)$ and $s \in H$ of order $2$ (and semisimple as ${\rm char}\, k \neq 2$). Then the centralizer $C$ of $s$ in $H$ is isomorphic to ${\rm SO}({\rm H}(k^2))$ and we have a decomposition $E \,= \,E_4\, \bot\, E_3$ where $E_4 = {\rm ker}(s+{\rm id})$ gives rise to the natural representation representation of $C$ of dimension $4$ (we omit the details here). In particular, we have $\det E_4 = 1$ as $k[C]$-module. On the other hand, assume that we have $s=\rho(s')$ with $s' \in {\rm O}_2^\pm(k)$ such that $\mu(s')=1$ and $\epsilon(s')=-1$, and that $\Gamma$ is the centralizer of $s'$. As $s'$ is not a square in ${\rm O}_2^\pm(k)$, the $4$ characters $\chi : \Gamma \rightarrow k^\times$ such that $\chi(s') = -1$  are of the form $c^i c'$, $0 \leq i \leq 3$, where $c$ has order $4$, $c'$ has order $2$, and  $c(s') = 1 = - c'(s')$. In particular, their product is $c^2 \neq 1$, a contradiction.
\end{pf}
\ps\ps

The group ${\rm O}_2^\pm(k)$ will appear in practice using the following lemma. 
\ps\ps

\begin{lemme}\label{crito2pm} Let $\Gamma$ be a group, $k$ an algebraically closed field with ${\rm char}\, k\, \neq 2$, $Q$ an irreducible $k[\Gamma]$-module of dimension $2$, and $c : \Gamma \rightarrow k^\times$ a character with $c^2=1$ occurring in ${\rm Sym}^2 Q$. Then there is a unique $ {\rm O}_2(k)$-conjugacy class of morphisms $\Gamma \rightarrow {\rm O}_2^\pm(k)$ whose similitude factor is $c$ and such that the $k[\Gamma]$-module ${\rm H}(k)$ is isomorphic to $Q$.
\end{lemme}

\begin{pf} As ${\rm char}\, k \neq 2$, we have natural isomorphisms ${\rm Hom}_k(Q^\ast,Q) \simeq Q \otimes Q \simeq {\rm Sym}^2 Q \oplus \det Q$; moreover, these $k[\Gamma]$-modules are semisimple as so is $Q$, by \cite{serre2}. By assumption, there is thus a nonzero $k[\Gamma]$-linear map ${\rm Sym}^2\, Q \rightarrow c$. Such a symmetric bilinear form on $Q$ is necessarily nondegenerate, as $Q$ is simple. It gives $Q$ a structure of nondegenerate quadratic space, necessarily isomorphic to ${\rm H}(k)$ as $k$ is algebraically closed. This proves the existence part of the lemma. Assume there are two $k[\Gamma]$-linear isomorphisms $u_i : Q \isomo {\rm H}(k)$, $i=1,2$, such that $\Gamma$ acts as similitude of factor $c$ on the nondegenerate quadratic form ${\rm q} \circ u_i$. By the first line of this proof, the $\Gamma$-invariants in ${\rm Hom}_k({\rm Sym}^2\,Q,c)$ have dimension $1$, so those two forms are proportional, i.e. $u_1 \circ u_2^{-1} \in {\rm GO}_2(k)$. We conclude as ${\rm GO}_2(k) = k^\times \cdot {\rm O}_2(k)$.
\end{pf}
\ps\ps

\begin{example}\label{exempled8} {\rm Consider the dihedral group ${\rm D}_8$ of order $8$. Denote by $\mathcal{E}$ the set of characters ${\rm D}_8 \rightarrow k^\times$; we have $|\mathcal{E}|=4$. There is a unique $\eta_0 \in \mathcal{E}$ whose kernel is the unique subgroup of ${\rm D}_8$ of order $4$. We claim that there are exactly $3$ conjugacy classes of morphisms $r : {\rm D}_8 \rightarrow {\rm O}_2^\pm(k)$, distinguished by the character $\mu \circ r$, which can be any element of $\mathcal{E} -\{\eta_0\}$. Indeed, up to isomorphism, there is a unique faithful $2$-dimensional $k[{\rm D}_8]$-module $J$. It is irreducible and satisfies $\det J = \eta_0$. In particular, $J$ is selfdual and we have $J \simeq J \otimes \eta$ for all $\eta \in \mathcal{E}$, so we have ${\rm Sym}^2 J \simeq \bigoplus_{\eta \in \mathcal{E}-\{\eta_0\}} \eta$ and Lemma \ref{crito2pm} proves the claim. Note however that the morphisms $\rho \circ r: {\rm D}_8 \rightarrow {\rm SO}(E)$, when $r$ varies as above, are all conjugate under ${\rm SO}(E)$ by Corollary \ref{corwindex}.}
\end{example}

\ps\ps\ps

{\sc The group ${\rm GL}_2(\Z/3\Z)$}\ps\ps

 The natural action of ${\rm GL}_2(\Z/3\Z)$ on the set of lines in $\Z/3\Z \times \Z/3\Z$, numbered in an arbitrary way, defines a surjective morphism $${\rm GL}_2(\Z/3\Z) \rightarrow \got{S}_4$$ with kernel $(\Z/3\Z)^\times$, which realizes ${\rm GL}_2(\Z/3\Z)$ as a central extension of $\got{S}_4$ by $\Z/2\Z$. This allows to view any $\got{S}_4$-module as a ${\rm GL}_2(\Z/3\Z)$-module by inflation.\ps\ps

Let $k$ be an algebraically closed field of characteristic $\neq 2,3$, so that $|{\rm GL}_2(\Z/3\Z)| \in k^\times$. Up to isomorphism, the irreducible $k[\got{S}_4]$-modules are of the form $1, c, H, V$ and $V \otimes c$, where $c$ is the signature, $H$ is a $2$-dimensional representation inflated from a surjective homomorphism $\got{S}_4 \rightarrow \got{S}_3$ and $V$ is $3$-dimensional with determinant $1$. Moreover, still up to isomorphism, the irreducible $k[{\rm GL}_2(\Z/3\Z)]$-modules with nontrivial central character are of the form $P$, $P^\ast$ and $P \otimes H$, where $P$ has dimension $2$. We have $\det H = \det P=c$, $H^\ast \simeq H \simeq H \otimes c$, $P^\ast \simeq P \otimes c$, $V^\ast \simeq V$, as well as the following identities:
\begin{equation} \label{idgl2f3}{\rm Sym}^2 H  \simeq 1 \oplus H,\hspace{.4 cm}  H \otimes V \simeq V \oplus V\otimes c \hspace{.4 cm} {\rm and} \hspace{.4 cm} P \otimes P^\ast \simeq 1 \oplus V.\end{equation}
\ps

\begin{prop} \label{rhogl23}Let $k$ be an algebraically closed field of characteristic $\neq 2,3$ and $E$ a nondegenerate quadratic space of dimension $7$ over $k$. There is a unique conjugacy class of injective morphisms $\rho : \GL_2(\Z/3\Z) \rightarrow {\rm SO}(E)$ such that the $k[\GL_2(\Z/3\Z)]$-module $E$ is isomorphic to $P \oplus P^\ast \oplus c \oplus H$, where $P$, $c$ and $H$ are as above.  For such a $\rho$, and any $\gamma \in {\rm GL}_2(\Z/3\Z)$, then $\det(t - \rho(\gamma))$ is of type ${\bf G}_2$. 
\end{prop}

\begin{pf} Recall that $c$ has order $2$ and that we have $\det H = c$ and $H^\ast  \simeq H$, so there is a structure of nondegenerate quadratic space on $H$ such that $\got{S}_4$ acts as orthogonal isometries. The existence of $\rho$ (of Witt index $2$) follows, and its uniqueness up to conjugacy is a consequence of Corollary \ref{corwindex}. To check the last assertion we check that assertion (iv) of Proposition \ref{propfondso7g2} holds. Using the identities \eqref{idgl2f3}, it is tedious but straightforward to check that in ${\rm R}_k({\rm GL}_2(\Z/3\Z))$ both elements $[E] + [{\rm Sym}^2 E]$ and $[\Lambda^3 E]$ are equal to $3 + c +  3\,[H] + 2\,[P] +2\,[P^\ast]  + [V] + 2 \,[V \otimes c] + 2\, [H \otimes P]$.
\end{pf}
\ps\ps
A consequence of our main Theorem \ref{mainthmb} (or even Proposition \ref{lemmaindex3}) will be that there is no ${\rm G}_2$-subgroup of ${\rm SO}(E)$ containing $\rho(\SL_2(\Z/3\Z))$ (hence $\rho(\GL_2(\Z/3\Z))$ {\it a fortiori}).  \ps\ps
\begin{remark}\label{rhosl23}
{\rm The group ${\rm SL}_2(\Z/3\Z)$, which is the kernel of $c$, is the universal central extension of $\got{A}_4$ by $\Z/2\Z$.  Its abelianization is isomorphic to $\Z/3\Z$, and if $\eta : {\rm SL}_2(\Z/3\Z) \longrightarrow k^\times$ is an order $3$ character, then we have  $H_{|{\rm SL}_2(\Z/3\Z)} \simeq \eta \oplus \eta^{-1}$. Moreover, the restriction of $P$ to ${\rm SL}_2(\Z/3\Z)$ is, up to isomorphism, the unique two-dimensional irreducible $k[{\rm SL}_2(\Z/3\Z)]$-module of determinant $1$. 
}
\end{remark}
\ps\ps

\subsection{Assumption {\rm (S)} and statement of the main theorem}\label{parass}

\begin{definition}\label{defiasss} Let $E$ be a regular quadratic space over an algebraically closed field $k$, $\Gamma$ a subgroup of ${\rm SO}(E)$, $\widetilde{\Gamma}$ the inverse image of $\Gamma$ under $\pi_E : {\rm Spin}(E) \rightarrow {\rm SO}(E)$, and $W$ a spinor module for $E$. We shall say that {\it $\Gamma$ satisfies assumption {\rm (S)}} if  $E$, ${\rm Sym}^2\,E$, $\Lambda^3 \, E$ and $W$ are semisimple $k[\widetilde{\Gamma}]$-modules. 
\end{definition}

This technical assumption {\rm (S)} (where $S$ stands for {\rm ``semisimplicity''}) will be quite natural in our proof of Theorem \ref{mainthmb}. It is of course satisfied if $\Gamma$ is finite with $2\,|\Gamma| \in k^\ast$. The next proposition shows that assumption {\rm (S)} is not as strong as it may seem.\ps\ps

\begin{prop} \label{critasss}Let $E$ be a nondegenerate quadratic space of dimension $7$ over an algebraically closed field $k$, and $\Gamma \subset {\rm SO}(E)$ a subgroup. Assume that at least one of the following assumptions holds:\ps\ps
\begin{itemize} 
\item[(i)] $k=\C$ and the closure of $\Gamma$ in the Lie group ${\rm SO}(E)$ is compact, \ps\ps
\item[(ii)] ${\rm char}\, k = 0$ or ${\rm char}\, k >13$, and the $k[\Gamma]$-module $E$ is semisimple, \ps\ps
\item[(iii)] the $k[\Gamma]$-modules $E \otimes E$ and $\Lambda^3 E$ are semisimple.\ps\ps
\end{itemize}
Then $\Gamma$ satisfies assumption {\rm (S)}.
\end{prop}

\begin{pf} Let $G \subset \GL_n(\C)$ be a subgroup with compact closure $\overline{G}$. A standard argument using a Haar measure of $\overline{G}$ shows that $\C^n$ possesses a Hermitian inner product which is $\overline{G}$-invariant, hence $G$-invariant as well. In particular the $\C[G]$-module $\C^n$ is semisimple. Under assumption (i), $\widetilde{\Gamma}$ has also a compact closure in the Lie group ${\rm SO}(E)$, which proves that $\Gamma$ satisfies assumption {\rm (S)}. \ps\ps

Observe that the assumption on $E$ implies ${\rm char}\, k \neq 2$. By Lemma \ref{lemmesymlambda2} we have $k[\widetilde{\Gamma}]$-linear isomorphisms \begin{equation}\label{2isolemmsim} \Lambda^2 \, W \simeq E \oplus \Lambda^2 E \hspace{1 cm}{\rm and}\hspace{1 cm} {\rm Sym}^2 W \simeq 1 \oplus \Lambda^3 E.\end{equation}
We also have $W \otimes W \simeq {\rm Sym}^2\, W \oplus \Lambda^2\, W$ and $E \otimes E \simeq {\rm Sym}^2 E \oplus \Lambda^2 E$. By Serre \cite{serre2}, the semisimplicity of $W \otimes W$ (resp $E \otimes E$) imply that of $W$ (resp. $E$). As a consequence, assumption (iii) implies that $\Gamma$ satisfies assumption {\rm (S)}. \ps\ps

Assume that the $k[\Gamma]$-module $E$ is semisimple. If ${\rm char}\, k = 0$, a classical result of Chevalley ensures that $E^{\otimes n}$ is semisimple for each integer $n\geq 0$, so assumptions (iii) and {\rm (S)} are satisfied. Set $p = {\rm char}\, k$ and assume $p>13$. As $p>2 \dim E -2$, $\Lambda^2 E$ and ${\rm Sym}^2$ are semisimple by the main result of Serre \cite{serre}. In particular, $\Lambda^2 \, W$ is also semisimple by the first isomorphism of \eqref{2isolemmsim}. As $8 = \dim W \not \equiv 2,3 \bmod p$, another result of Serre \cite{serre2} ensures the semisimplicity of $W$. But by \cite{serre} again this implies that ${\rm Sym}^2\, W$ is semisimple as $p> 2 \dim W - 2$. We deduce that $\Lambda^3\, E$ is semisimple by the second isomorphism of \eqref{2isolemmsim}. We have proved that assumption (iii) implies assumption {\rm (S)}.\end{pf}
\ps\ps
\begin{remark} \label{hypSex}{\rm Let $\rho : {\rm O}_2^{\pm}(k) \longrightarrow {\rm SO}(E)$ be as in Proposition \ref{o2pm} (in particular, we have ${\rm char}\, k \neq 2$) and let $\Gamma \subset {\rm O}_2^{\pm}(k)$ be a subgroup. We claim that $\rho(\Gamma)$ satisfies assumption {\rm (S)}. Indeed, it is enough to show that $\rho(\Gamma')$ satisfies {\rm (S)}, where $\Gamma'$ is the kernel of $\mu\epsilon$, as $|\Gamma/\Gamma'| \in k^\times$. But this kernel is a subgroup of ${\rm GSO}_2(k)$, which acts as a direct sum of two characters on $P$. It follows that $E^{\otimes n}$ is a direct sum of characters of $\Gamma'$ for all $n\geq 0$, and we are done by Proposition \ref{critasss} (iii).}
\end{remark}
\ps\ps

We can now state our main theorem.  \ps\ps

\begin{thm} \label{mainthmb} Let $E$ be a $7$-dimensional nondegenerate quadratic space over the algebraically closed field $k$ and $\Gamma \subset {\rm SO}(E)$ a subgroup. Assume that: \ps\ps
\begin{itemize}
\item[(i)] $\Gamma$ satisfies assumption {\rm (S)}, \ps\ps
\item[(ii)] each element of $\Gamma$ has a characteristic polynomial of type ${\bf G}_2$,\ps\ps
\end{itemize}
Then exactly one of the following assertions holds: \begin{itemize} \ps \ps
\item[(a)] there is a ${\rm G}_2$-subgroup of ${\rm SO}(E)$ containing $\Gamma$,\ps\ps
\item[(b)] ${\rm char}\, k \neq 3$ and we have $\Gamma = \rho({\rm GL}_2(\Z/3\Z))$ or $\Gamma = \rho({\rm SL}_2(\Z/3\Z))$, for some $\rho : {\rm GL}_2(\Z/3\Z) \longrightarrow {\rm SO}(E)$ as in Proposition \ref{rhogl23},\ps\ps
\item[(c)] we have $\Gamma \simeq \Z/2\Z \times \Z/4\Z$ and the $k[\Gamma]$-module $E$ is the direct sum of the nontrivial characters of $\Gamma$,\ps\ps
\item[(d)] we have $\Gamma = \rho(\Gamma')$ where $\rho : {\rm O}_2^\pm(k) \longrightarrow {\rm SO}(E)$ is as in Proposition \ref{o2pm}, and $\Gamma' \subset {\rm O}_2^{\pm}(k)$ is a nonabelian subgroup, nonisomorphic to ${\rm D}_8$, and such that $\mu(\Gamma')=\{\pm1\}$.\ps\ps
\end{itemize}
In particular, we are in case {\rm (a)} if the Witt index of $\Gamma$ is $\leq 1$.
\end{thm}
\ps\ps

Note that the assumption on $E$ in the statement above implies ${\rm char}\, k \neq 2$ (but when ${\rm char}\, k\,=2$ we can apply Corollary \ref{corcar2} instead). Observe also that the ``exceptional'' cases (b), (c) and (d) do actually occur, as is shown by Propositions \ref{o2pm} \& \ref{rhogl23}, Lemma \ref{lemmacentor2} and Remarks \ref{rhosl23} \& \ref{hypSex}. \ps\ps

The end of this section is devoted to the proof of Theorem \ref{mainthmb}. {\it From now on and until the end of \S \ref{proofthmb}, the letter $k$ will always denote an algebraically closed field of characteristic $\neq 2$}. Note that under assumption {\rm (S)}, Proposition \ref{propfondso7g2} shows that assumption (ii) is equivalent to the existence of a $k[\Gamma]$-module isomorphism 
\begin{equation}\label{eqfundg2} \Lambda^3 E \simeq E \oplus {\rm Sym}^2 E.\end{equation}
We shall study this identity and argue according to the Witt index of $\Gamma$, which is an element in $\{0,1,2,3\}$, in the decreasing order. One reason for that is that if $\Gamma' \subset \Gamma$ is a subgroup satisfying {\rm (S)}, then $\Gamma'$ satisfies the assumptions of the theorem and its Witt index is at least the one of $\Gamma$. If $\Gamma$ satisfies (b) then its Witt index is $3$ if $\Gamma \simeq {\rm SL}_2(\Z/3\Z)$ and $2$ if $\Gamma \simeq {\rm GL}_2(\Z/3\Z)$. If $\Gamma$ satisfies (c) its Witt index is $2$. If $\Gamma$ satisfies (d) then its Witt index is $3$ if $\epsilon(\Gamma')=1$ and $2$ otherwise. This explains the last assertion. \ps

\subsection{The case of Witt index $3$}\label{parcasewindex3}
\ps
\begin{prop} \label{lemmaindex3} Let $E$ be a $7$-dimensional nondegenerate quadratic space over $k$ and $\Gamma \subset {\rm SO}(E)$ a subgroup of Witt index $3$ satisfying assumptions (i) and (ii) of the statement of Theorem \ref{mainthmb}. Then exactly one of the following assertions holds: \ps 
\begin{itemize}
\item[(i)] $\Gamma$ is contained in a ${\rm G}_2$-subgroup of ${\rm SO}(E)$, \ps\ps
\item[(ii)]  there is an order $2$ character $c$ and a $2$-dimensional irreducible representation $J$ of $\Gamma$ such that $[E]=1+2c+2[J]$ in ${\rm R}_k(\Gamma)$, $J \simeq J \otimes c$ and $\det J =1$,\ps\ps \ps
\item[(iii)] ${\rm char}\, k \neq 3$, $\Gamma \simeq {\rm SL}_2(\Z/3\Z)$ and $[E]=2[J]+1+c+c^{-1}$ where $J$ is ``the'' $2$-dimensional irreducible representation of $\Gamma$ with determinant $1$ and $c^{\pm 1}$ are the two order $3$ characters of $\Gamma$.
\end{itemize}
\end{prop}
\ps
Let us first introduce a notation. If $X$ is a $k[\Gamma]$-module of finite dimension over $k$, we shall denote by ${\rm h}(X)$ the $k[\Gamma]$-module $X \oplus X^\ast$. \ps\ps

\begin{pf} We may assume that $E = {\rm H}(I) \bot k $ with $\dim I = 3$ and that $\Gamma$ preserves the subspaces $I$ and $I^\ast$. In particular, we have a $k[\Gamma]$-linear isomorphism $$E \simeq {\rm h}(I) \oplus 1.$$ 
Let $\widetilde{\Gamma}$ denote the inverse image of $\Gamma$ under $\pi_E : {\rm Spin}(E) \rightarrow {\rm SO}(E)$.
By Propositions \ref{cliffordhyp} and \ref{decospin},  $\widetilde{\Gamma}$ is included in the natural subgroup $\rho_{{\rm H}(I);k}(\widetilde{\rho_I}({\rm GL}(I)) \times k^\times)$ of ${\rm GSpin}(E)$, and if $\alpha : \widetilde{\Gamma} \rightarrow k^\times$ denotes the projection on the $k^\times$-factor (which is also the center of ${\rm GSpin}(E)$), then we have $(\nu_E)_{|\widetilde{\Gamma}} = \alpha^2\, \det\,I=1$ and ${W}_{|\widetilde{\Gamma}} \simeq \,\alpha \,\otimes \,\Lambda\, I$. If we set $c = \Lambda^3 \, I = \det I$, we have thus $c \alpha^2 =1 $ and $\Lambda^2 I \simeq I^\ast \otimes c$, so we have a $k[\widetilde{\Gamma}]$-linear isomorphism 
\begin{equation}\label{decowind3}W \simeq  {\rm h}(\alpha  \oplus \alpha \otimes I).\end{equation} By assumption (S) and Corollary~\ref{carbetaspin},  $\Gamma$ is contained in a ${\rm G}_2$-subgroup of ${\rm SO}(E)$ if, and only if, there is a character $\beta : \widetilde{\Gamma} \rightarrow k^\times$ occurring in $W_{|\widetilde{\Gamma}}$ and such that $\beta^2=1$. Note that $W_{|\widetilde{\Gamma}}$ contains $\alpha$ as well as the relation $\alpha^2 = c^{-1}$. \ps\ps

We have $\Lambda^2 I \simeq I^\ast \otimes c$, $\Lambda^2 I \otimes I^{\ast} \simeq {\rm Sym}^2 I^\ast \otimes c \oplus I$ and $E \simeq {\rm h}(I) \oplus 1$. A straightforward expansion shows that we have $k[\Gamma]$-linear isomorphisms
$$\Lambda^3\, E \,\simeq \,{\rm h}(c \oplus ({\rm Sym}^2\, I^\ast) \otimes c \oplus I \oplus I^\ast \otimes c \, ) \oplus I \otimes I^\ast,$$
$$ {\rm Sym}^2\, E \,\simeq \,{\rm h}({\rm Sym}^2\, I \oplus I) \oplus I \otimes I^\ast \oplus 1.$$
In particular, $I$, $I \otimes I^\ast$ and ${\rm Sym}^2 I$ are semisimple by assumption ({\rm S}). Equation \eqref{eqfundg2} is thus equivalent to
$${\rm h}(x) \simeq {\rm h}(c^{-1} \otimes x) \, \, \, \,{\rm with}\, \,\, \, x = 1 \oplus I \oplus {\rm Sym}^2\, I.$$

  \ps \ps
In particular, the trivial representation $1$ of $\Gamma$ is a summand of $c^{-1}x$, i.e. the character $c$ is a summand of $x$. \ps \ps

{\it Case 1 : the $k[\Gamma]$-module $I$ is irreducible}. It follows that either $c=1$ or $c$ is a summand of ${\rm Sym}^2 I$. This second possibility implies that there is a nonzero $\Gamma$-equivariant map $I^\ast \rightarrow I \otimes c^{-1}$. As $I$ is irreducible, such a map is an isomorphism, and taking determinant gives the relation $c^{-1} = c c^{-3}=c^{-2}$. So we have $c=1$ in all cases, hence $\alpha^2=1$ and $\Gamma$ is contained in a ${\rm G}_2$-subgroup of ${\rm SO}(E)$. \ps\ps

{\it Case 2 : $I \simeq \oplus_i \chi_i$ is a sum of $3$ characters $\chi_i : \Gamma \rightarrow k^\times$, $i=1,2,3$.} We have $x \simeq 1 \oplus \bigoplus_i (\chi_i \oplus \chi_i^2 \oplus c \chi_i^{-1})$, hence also $c^{-1} x \simeq c^{-1} \oplus \bigoplus_i (c^{-1} \chi_i \oplus c^{-1}\chi_i^2 \oplus \chi_i^{-1})$. The relation ${\rm h}(x) \simeq {\rm h}(c^{-1}x)$ is thus equivalent to 
 $${\rm h}(1\oplus \bigoplus_i \chi_i^2) \simeq {\rm h}(c^{-1}\oplus \bigoplus_{i} c^{-1}\chi_i^2).$$ As a consequence, $1$ is a summand of $c^{-1}\oplus \bigoplus_{i} c^{-1}\chi_i^2$, and we have either $c=1$ or $c=\chi_i^2$ for some $i$. We conclude the proof by choosing accordingly $\beta=\alpha$ or $\beta = \alpha \chi_i$, which occurs in $W$ and satisfies $\beta^2=1$. \ps\ps

{\it Case 3 : $I \simeq \chi \oplus J$ with $J$ an irreducible representation of $\Gamma$ of dimension $2$.} We have $\det J = c \chi^{-1}$. The isomorphism \eqref{decowind3} shows that $\Gamma$ is contained in a ${\rm G}_2$-subgroup if, and only if, we have $c=1$ or $c= \chi^2$. Assume $c \neq 1$ and $c \neq \chi^2$. We have 
$$x \simeq 1 \oplus \chi \oplus \chi^2 \oplus J \oplus \chi \otimes J \oplus  {\rm Sym}^2 J,$$ so either $c=\chi$ or $c$ is a summand of ${\rm Sym}^2 J$. In this latter case, we have $J \simeq J^\ast \otimes c$ and by taking the determinant $c \chi^{-1} = c^2 c^{-1} \chi$, i.e. $\chi^2=1$. If follows that $1$ occurs twice in $x$, hence twice in $c^{-1}x$ as well.  As ${\rm Hom}_{k[\Gamma]}(J,J^\ast \otimes c)$ has dimension $1$ by irreducibility of $J$, it follows that $c=\chi$ in all cases. In particular we have $\det J=1$ and $J^\ast \simeq J$. The equation ${\rm h}(x) \simeq {\rm h}(c^{-1}x)$ simplifies as 
\begin{equation} \label{eqind3bis} {\rm h}(c^2) \oplus  ({\rm Sym}^2 J)^{\oplus 2} \simeq {\rm h}(c) \otimes (1 \oplus {\rm Sym}^2 J).\end{equation}\ps\ps

 {\it Subcase {\rm (a)}.} Assume first that either $c$ or $c^{-1}$ occurs in ${\rm Sym}^2 J$. As $J$ is selfdual this implies $J \simeq J \otimes c$, hence $c$ has order $2$, and we are in case (ii). \ps\ps
 {\it Subcase {\rm (b)}.}  Assume now that $c^{\pm 1}$ does not occur in ${\rm Sym}^2 J$. As $c \neq 1$, equation \eqref{eqind3bis} shows $c=c^{-2}$, {\it i.e.} $c$ has order $3$, as well as a $k[\Gamma]$-linear isomorphism
 \begin{equation} \label{eqind3bisbis} {\rm Sym}^2\, J \oplus  \, {\rm Sym}^2\, J\,\simeq \, (c\oplus c^{-1}) \otimes {\rm Sym}^2 \, J. \end{equation}
 We will eventually show that we are in case (iii). As $c$ has order $3$, observe that we have $3 \in k^\times$. \ps\ps
 
 Let us show that ${\rm Sym}^2 \, J$ is irreducible. Otherwise, some character $\mu$ of $\Gamma$ occurs in ${\rm Sym}^2 \, J$ (which is $3$-dimensional and semisimple). The relation \eqref{eqind3bisbis} shows that each of the three distinct characters $\mu\, c^i$, $i \in \{-1,0,1\}$, occur in ${\rm Sym}^2 \, J$, so that we have ${\rm Sym}^2\, J \simeq \mu \oplus \mu c \oplus \mu c^{-1}$. As $J$ is selfdual and irreducible of determinant $1$, those characters must all have order $2$. But ${\rm Sym}^2\, J$ has determinant $(\det J)^3 = 1$, so we also have $\mu^3=1$: a contradiction. So ${\rm Sym}^2 \, J$ is irreducible and equation \eqref{eqind3bisbis} implies $${\rm Sym}^2 \, J \simeq c \otimes {\rm Sym}^2 \,J.$$ 
 
 Consider the subgroup $\Gamma^0 := {\rm ker}\, c$, which is normal of index $3$ in $\Gamma$. Clifford theory\footnote{By this we mean the following classical facts. Let $k$ be an algebraically closed field (of arbitrary characteristic), $\Gamma$ a group, $c : \Gamma \rightarrow k^\times$ a character of finite order $n$, and $V$ a finite dimensional irreducible $k$-linear representation of $\Gamma$ such that $V \simeq V \otimes c$. Then $n$ divides $\dim V$ (apply $\det$), and if $\Gamma^0$ denotes the kernel of $c$ then $V_{|\Gamma^0}$ is a direct sum of $n$ irreducible representations of $\Gamma^0$ of dimension $\frac{\dim V}{n}$, whose isomorphism classes are distinct and permuted transitively by the outer action of $\Gamma/\Gamma^0 \simeq \Z/n\Z$ on $\Gamma^0$. Moreover, if $W$ denotes any such representation we have $V \simeq {\rm Ind}_{\Gamma^0}^\Gamma W$. }  shows that $({\rm Sym}^2\, J)_{|\Gamma^0}$ is a direct sum of $3$ distinct characters permuted transitively by the outer action of $\Gamma/\Gamma^0 \simeq \Z/3\Z$ on $\Gamma^0$. As $\Gamma^0$ is a normal subgroup of index $3$ in $\Gamma$, and as $\dim J = 2$, note that $J_{|\Gamma^0}$ is irreducible. It follows that each of the $3$ characters in $({\rm Sym}^2\, J)_{|\Gamma^0}$ has order $2$, the product of all of them being $1$. In particular, if $\Gamma^1$ denotes the kernel of the morphism $\Gamma \rightarrow {\rm GL}({\rm Sym}^2 J)$, then we have $\Gamma^1 \subset \Gamma^0$ and $\Gamma/\Gamma^1$ is an extension of $\Gamma/\Gamma^0 \simeq \Z/3\Z$ by $\Gamma^0/\Gamma^1 \simeq (\Z/2\Z)^2$ such that the associated outer action of $\Z/3\Z$ on $(\Z/2\Z)^2$ permutes transitively the three nonzero elements.  It follows that $\Gamma/\Gamma^1$ is isomorphic to $\got{A}_4$ (and that its irreducible $k$-linear representations are, up to isomorphism, ${\rm Sym}^2 J$, $1$, $c$ and $c^2$). \ps\ps
 
As $\Gamma^1$ acts trivially on ${\rm Sym}^2 \, J$, it acts by scalars of square $1$ in $J$, so $|\Gamma^1|\leq 2$.  As $\got{A}_4$ has no irreducible $k$-linear representation of dimension $2$, we have $\Gamma^1 \neq 1$, {\it i.e.} $\Gamma^1 \simeq \Z/2\Z$. As $E \simeq J \oplus J \oplus 1 \oplus c \oplus c^{-1}$, one sees that $\Gamma^1$ is central in $\Gamma$. We know from Schur that there is a unique nontrivial central extension of $\got{A}_4$ by $\Z/2\Z$, sometimes denoted $\widetilde{\got{A}_4}$, and which is isomorphic to ${\rm SL}_2(\Z/3\Z)$; moreover, such a group has exactly three $2$-dimensional irreducible $k$-linear representations, twists of each others, which are distinguished by their determinant (which can be any of the three characters of $\got{A}_4$). As a consequence, we have $\Gamma \simeq \SL_2(\Z/3\Z)$ and assertion (iii) holds.  \end{pf}

\begin{pf} Let us prove Theorem \ref{mainthmb} when $\Gamma$ has Witt index $3$. Given Proposition \ref{lemmaindex3} it is enough to show that: \ps\ps 
(1) we are in case (ii) of the proposition if, and only if, we are in case (d) of the statement of Theorem \ref{mainthmb} with $\epsilon(\Gamma')=1$. \ps\ps
(2) we are in case (iii) if, and only if, we are in case (b) of Theorem \ref{mainthmb} and $\Gamma=\rho({\rm SL}_2(\Z/3\Z))$.\ps\ps

Assertion (2) is an immediate consequence of Proposition \ref{rhogl23} and Remark \ref{rhosl23}. Let us check assertion (1). Observe first that a subgroup $\Gamma' \subset {\rm O}_2^{\pm}(k)$ acts irreducibly on $P={\rm H}(k)$ if, and only if, it is nonabelian, as $P$ is a $2$-dimensional semisimple $k[\Gamma']$-module by Remark \ref{hypSex}. This shows that if we are in case (d) of the statement of Theorem \ref{mainthmb} with $\epsilon(\Gamma')=1$, then we are in case (ii) of the statement of Proposition \ref{lemmaindex3} (with $c$ playing the role of the character $\mu$). Assume conversely that we are in case (ii).  As $c$ occurs in $J \otimes J^\ast$, the natural morphism $\Gamma \rightarrow {\rm GL}(J)$ is injective. Moreover, we have $c \neq 1$, $\det J = 1$ and $J \simeq J^\ast$, hence $c$ occurs in ${\rm Sym}^2 J$, and we conclude the proof by Lemma \ref{crito2pm} and Proposition \ref{o2pm}. Note that $\epsilon(\Gamma')=1$ implies $\Gamma' \not\simeq {\rm D}_8$ by the discussion in Example \ref{exempled8}. \end{pf} 
\ps\ps
\subsection{The case of Witt index $2$}
\ps\ps
We shall now consider the Witt index $2$ case. If $\Gamma$ is a group and $P$ a finite dimensional $k$-linear representation of $\Gamma$, we shall denote by ${\rm ad}\, P\, \subset {\rm Hom}(P,P) \simeq P^\ast \otimes P$ the subspace of trace $0$ endomorphisms. If $\dim P =2$, we natural isomorphisms $P^\ast \otimes P \simeq 1 \oplus {\rm ad} \,P$ and ${\rm Sym}^2\, P\, \simeq \,\det P \otimes {\rm ad} \,P$ (recall ${\rm char}\, k \neq 2$). 

\ps\ps
\begin{prop} \label{lemmaindex2} Let $E$ be a $7$-dimensional nondegenerate quadratic space over $k$ and $\Gamma \subset {\rm SO}(E)$ a subgroup of Witt index $2$ satisfying assumptions (i) and (ii) of the statement of Theorem \ref{mainthmb}. Then exactly one of the following assertions holds: \ps 
\begin{itemize}
\item[(i)] $\Gamma$ is contained in a ${\rm G}_2$-subgroup of ${\rm SO}(E)$, \ps\ps
\item[(ii)]  $\Gamma \simeq \Z/4\Z \times \Z/2\Z$, and as a $k[\Gamma]$-module $E$ is isomorphic to the direct sum of the seven nontrivial characters of  $\Gamma$,\ps \ps
\item[(iii)] there is an irreducible, nonselfdual, $2$-dimensional representation $P$ of $\Gamma$, as well as two distinct order $2$ characters $c,\epsilon$ of $\Gamma$, such that $\det P = c$, $P \simeq P \otimes \epsilon$ and $[E] = [P] + [P^\ast] + c + \epsilon + c\epsilon$ in ${\rm R}_k(\Gamma)$. \ps\ps
\item[(iv)] ${\rm char}\, k \neq 3$, $\Gamma \simeq {\rm GL}_2(\Z/3\Z)$, and as a $k[\Gamma]$-module $E$ is isomorphic to the direct sum of the three $2$-dimensional irreducible representations of $\Gamma$ and of its nontrivial (order $2$) character.  \ps\ps
\end{itemize}
In case (i), we have $[E] = [P]+[P^\ast] + [{\rm ad}\, P]$ in ${\rm R}_k(\Gamma)$, where $P \subset E$ denotes any $2$-dimensional totally isotropic subspace stable by $\Gamma$, and $P$ is irreducible.
\end{prop}
\ps
\begin{pf} Let $P \subset E$ be a $\Gamma$-stable, two dimensional totally isotropic subspace of $E$. We may assume that we have $E = {\rm H}(P) \bot V$ and that $\Gamma$ preserves $P,P^\ast$ and $V$, with $V$ of Witt index $0$ and dimension $3$. In particular, we have a $k[\Gamma]$-linear isomorphism $E \simeq {\rm h}(P)\oplus V$ and $\det V = 1$. Set $c = \det P$. Equation \eqref{eqfundg2} on $\Gamma$ is easily seen to be equivalent to
\begin{equation} \label{eqind2} V \otimes (c \oplus c^{-1} \oplus {\rm ad} \, P) \simeq (1\oplus c \oplus c^{-1}) \otimes {\rm ad}\, P\oplus {\rm Sym}^2\, V,
\end{equation}
both sides being semisimple $k[\Gamma]$-modules. Let $\widetilde{\Gamma}$ denote the inverse image of $\Gamma$ under $\pi_E : {\rm Spin}(E) \rightarrow {\rm SO}(E)$. The sequence \eqref{rhoUV} and Proposition \ref{cliffordhyp} show that the group $\widetilde{\Gamma}$ is contained in the subgroup $$\rho_{{\rm H}(P);V}(\widetilde{\rho_P}({\rm GL}(P)) \times {\rm GSpin}(V)) \simeq {\rm GL}(P) \times {\rm GSpin}(V)$$
of ${\rm GSpin}(E)$ such that $\nu_E = \det P \,\,\cdot \,\,\nu_V = 1$. Denote by $W_0$ a spinor module for $V$. It is well-known that the natural morphism ${\rm GSpin}(V) \rightarrow {\rm GL}(W_0)$ is an (exceptional) isomorphism, and that we have $\det W_0 = \nu_V$ and $V = \,{\rm ad} \,W_0$. If $W$ denotes a spinor module for $E$, it follows from Propositions \ref{cliffordhyp} and \ref{decospin} that we have a $k[\widetilde{\Gamma}]$-linear isomorphism 
\begin{equation} \label{eqind2spin} 
W \simeq \Lambda P \otimes  W_0 \simeq (1 \oplus c \oplus P) \otimes W_0. \end{equation}
 \ps\ps
Observe that $W_0$ is irreducible as a representation of $\widetilde{\Gamma}$. Indeed, it is semisimple, being a direct summand of $W$. If $W_0$ is a sum of two characters of $\widetilde{\Gamma}$, then $V = \,{\rm ad} \,W_0$ contains the trivial character. But then we have a $\Gamma$-stable decomposition $V = k e \bot  V'$ with $V'$ of dimension $2$ and determinant $1$. Such a $V'$ is reducible (as it is so as a representation of ${\rm SO}(V') \simeq k^\times$), of Witt index $0$, hence a sum of two order $2$ characters. But these characters are necessarily equal as $\det V'=1$, a contradiction. \ps\ps

We claim that the following properties are equivalent:\label{pp1} \begin{itemize} \ps\ps
\item[(a)] $\Gamma$ is contained in a ${\rm G}_2$-subgroup of ${\rm SO}(E)$. \ps\ps
\item[(b)] There exists a character $\beta : \widetilde{\Gamma} \rightarrow k^\times$ with $\beta^2=1$ such that $\beta$ is a summand of $W_{|\widetilde{\Gamma}}$, or which is the same, such that $W_0^\ast$ is isomorphic to $P \otimes \beta$,\ps\ps
\item[(c)] We have ${\rm ad}\, P \simeq V$ as $k[\Gamma]$-modules. \ps\ps
\end{itemize}
Indeed, given Corollary \ref{carbetaspin}, the irreducibility of $W_0$ and the isomorphism \ref{eqind2spin}, it only remains to explain (c) $\Rightarrow$ (b). Observe first that if $U$ is a two-dimensional $k$-vector space, viewed as a $k[{\rm GL}(U)]$-module, then the determinant on ${\rm End}\, U$ equip ${\rm ad} \,U$ with a structure of $3$-dimensional nondegenerate quadratic space over $k$, and the natural map ${\rm ad} : {\rm GL}(U)/k^\times \rightarrow {\rm SO}({\rm ad}\, U)$ is an isomorphism. Let $\rho_i : \Gamma \rightarrow \GL(U)$, $i=1,2$, be two representations such that ${\rm ad} \,\rho_i$ is semisimple and such that ${\rm ad} \rho_1 \simeq {\rm ad} \rho_2$. The last assertion of Proposition \ref{propwindex} shows that there is some $g \in {\rm GL}(U)$ such that ${\rm ad}\, g \,\,{\rm ad} \rho_1\,\, {\rm ad} g^{-1} \,=\, {\rm ad} \rho_2$, i.e. $g \rho_1(\gamma) g^{-1}  \rho_2(\gamma)^{-1} \in k^\ast$ for all $\gamma \in \Gamma$. There is thus a unique group homomorphism $\Gamma \rightarrow k^\ast$ such that $\rho_1 = \rho_2 \otimes \chi$. If $\det \rho_1 = \det \rho_2$ then we have $\chi^2 = 1$. These observations apply to the representations $\rho_i$ defined by $P$ and $W_0^\ast$, and prove the equivalence between (b) and (c). \ps\ps

{\it Case 1 :  $P$ is reducible}. In particular, $\Gamma$ is not included in any ${\rm G}_2$-subgroup of ${\rm SO}(E)$ by the criterion (b) above. By assumption, there exists a character $\chi$ of $\Gamma$ such that $P \simeq \chi \oplus c \chi^{-1}$, and we thus have ${\rm ad}\, P \simeq 1 \oplus \chi^2 c^{-1} \oplus \chi^{-2} c$.  Let ${\rm n}(V)$ be the number of $1$-dimensional summands of $V$. An inspection of  equation \eqref{eqind2} shows the inequality $5 {\rm n}(V) \geq 9$, which implies ${\rm n}(V) \geq 2$, hence ${\rm n}(V)=3$ as $\dim V = 3$. So $V$ is a sum of $3$ characters, necessarily distinct of order $2$. In particular we have ${\rm Sym}^2 V \simeq 1 \oplus 1 \oplus 1 \oplus V$ and equation \eqref{eqind2} takes then the form: 
$$V \otimes (c \oplus c^{-1} \oplus \chi^2 c^{-1} \oplus \chi^{-2} c) \simeq (1\oplus c \oplus c^{-1}) \otimes  (1 \oplus \chi^2 c^{-1} \oplus \chi^{-2} c) \oplus 1 \oplus 1 \oplus 1.$$
As $V$ has Witt index $0$, the character $c^{\pm 1}$ (resp. $\chi^2 c^{-1}$, $\chi^{-2}$) occurs at most once in $V$, and if it does we have $c^2= 1$ (resp. $\chi^2 c^{-1} = \chi^{-2}c$). 
As the trivial representation occurs at least $4$ times on the right-hand side of the equation above, it occurs exactly $4$ times on both sides. It follows that both $c$ and $\chi^2 c^{-1}$ occur in $V$, and that we have $\chi^2 \neq  1$ and $\chi^2 \neq c$.  In particular, $c$ and $\chi^2 c^{-1}$ are distinct of order $2$, $\chi$ has order $4$, we have $V \simeq c \oplus \chi^2 \oplus \chi^2 c$
and 
$E \simeq \chi \oplus \chi^{-1} \oplus \chi\, c \oplus \chi^{-1}c \oplus \chi^2 \oplus c \oplus c \chi^2$. 
We are thus in case (ii) of the proposition. \ps\ps

{\em Case 2 : $P$ is irreducible and ${\rm ad}\, P$ is not isomorphic to $V$.} In particular,  $\Gamma$ is not included in any ${\rm G}_2$-subgroup of ${\rm SO}(E)$. We will show that we are in case (iii) of the proposition if $V$ is a sum of $3$ characters, and in case (iv) otherwise. \ps\ps

As $P$ is irreducible, observe that for any character $\epsilon$ of $\Gamma$, the dimension of the space 
${\rm Hom}_{k[\Gamma]} (\epsilon, P \otimes P^\ast ) \simeq {\rm Hom}_{k[\Gamma]}(P \otimes \epsilon,P)$ is at most $1$, and if it is nonzero we have $P \simeq P \otimes \epsilon$ and $\epsilon^2=1$. In particular, the semisimple representation ${\rm ad}\, P$ does not contain $1$. Moreover, both representations $V$ and ${\rm ad}\, P$ are semisimple, $3$-dimensional, multiplicity free and with trivial determinant. As a consequence, $\dim {\rm Hom}_{k[\Gamma]}(V,{\rm ad} \,P)$ is the number of irreducible summands of $V$ which occur in ${\rm ad}\, P$,  and we have $\dim {\rm Hom}_{k[\Gamma]}(V,{\rm ad} \,P) \leq 1$ as $V$ and ${\rm ad}\, P$ are not isomorphic. \ps\ps

Let us show now that $c$ is a constituent of $V$. Observe that $1$ is a summand of ${\rm Sym}^2 V$, hence of the right-hand side of equation \eqref{eqind2} as well. Assume that $V \otimes c^{\pm 1}$ does not contain $1$. Then there is a nonzero $\Gamma$-equivariant morphism $V \rightarrow {\rm ad}\, P$. By assumption, this is not an isomorphism, so $V$ is reducible. It follows that ${\rm Sym}^2 V$ contains at least twice the trivial representation, and again by \eqref{eqind2}, that ${\rm Hom}_{k[\Gamma]}(V, {\rm ad}\, P)$ has dimension $\geq 2$, which is absurd by the previous paragraph. We conclude that $V$ does contain $c^{\pm 1}$, hence $c$ itself. We may thus write $$V \simeq c \oplus H,$$
and we have $c \neq 1$, $c^2 = 1$, $\det H = c$ and $H^\ast \simeq H \simeq H \otimes c$. Equation~\eqref{eqind2} reduces then to 
\begin{equation} \label{eqind2bis} 1 \oplus H \otimes ({\rm ad}\, P \oplus 1) \simeq {\rm ad}\, P \otimes (1 \oplus c) \oplus {\rm Sym}^2 H.\end{equation}

{\em Subcase {\rm (a)} :  $H$ is reducible.} Write $H = \chi_1 \oplus \chi_2$ for some distinct order $2$ characters $\chi_i$, with $\chi_1 \chi_2 =c$. Relation \eqref{eqind2bis} becomes $(\chi_1 \oplus \chi_2) \otimes ({\rm ad}\, P \oplus 1) \simeq ({\rm ad}\, P \oplus 1) \otimes (1 \oplus c)$. It follows that $1$ occurs in the left-hand side, so some $\chi_i$, call it $\epsilon$, occurs in ${\rm ad}\, P$. We have thus $V = c \oplus \epsilon \oplus \epsilon c$ and $P \simeq P \otimes \epsilon$. We have the following equivalences: ${\rm ad} \,P$ is not isomorphic to $V$, $c$ is not a constituent of ${\rm ad} \,P$,  $P$ is not isomorphic to $P \otimes c$,  $P$ is not selfdual (as we have $P^\ast \simeq P \otimes {\det P}^{-1} \simeq P \otimes c$). We are thus in case (iii) of the proposition. \ps\ps

{\em Subcase {\rm (b)} :  $H$ is irreducible.} Observe that if $H$ is a summand of ${\rm ad} \, P$, then ${\rm ad}\, P \simeq H \oplus \det H \simeq V$, which is not the case by assumption. The isomorphism $H \simeq H \otimes c$, and the relation $c^2=1$, show that $H$ is not a summand of ${\rm ad} \, P \otimes c$ either. Equation \eqref{eqind2bis} shows thus that $H$ is a summand of ${\rm Sym}^2 H$. This equation is then equivalent to the following system of isomorphisms: 
\begin{equation}\label{eqind2bisbis} \left\{ \begin{array}{c} {\rm Sym}^2 H \simeq H \oplus 1 \\ \\ H \otimes {\rm ad}\, P \simeq (1\oplus c) \otimes {\rm ad}\, P \end{array} \right. \end{equation}
We claim that ${\rm ad}\, P$ is irreducible. Indeed, if we can write ${\rm ad}\, P = \epsilon \oplus T$ for some character $\epsilon$, we have already seen that $\epsilon$ has order $2$, the relation $\det T =\epsilon$, and the isomorphisms $T \simeq T^\ast \simeq T \otimes \epsilon$. By irreducibility of $H$, the second equation of \eqref{eqind2bisbis} shows that $T$ is irreducible and that $H \otimes \epsilon$  is isomorphic to $T$ or $T \otimes c$. But the isomorphisms $T \simeq T \otimes \epsilon$ and $H \simeq H \otimes c$ show $H \simeq T$ in both cases, which implies $V \simeq {\rm ad}\, P$ again: a contradiction. \ps\ps

Set $\Gamma^0 = {\rm ker}\, c$, it is an index $2$ subgroup of $\Gamma$. 
The relation $H \simeq H \otimes c$ and Clifford theory show that there is some character $\lambda : \Gamma^0 \rightarrow k^\times$ such that $H \simeq {\rm Ind}_{\Gamma^0}^\Gamma \lambda$. As $\det H = c$, we have $H_{|\Gamma^0} \simeq \lambda \oplus \lambda^{-1}$ and $\lambda^2 \neq 1$. As $H \simeq H^\ast$, we also have ${\rm Sym}^2 H \simeq 1 \oplus {\rm Ind}_{\Gamma^0}^\Gamma  \lambda^2$. As a consequence, we have an isomorphism $H \simeq {\rm Ind}_{\Gamma^0}^\Gamma  \lambda^2$. This forces either $\lambda = \lambda^2$ or $\lambda^{-1} = \lambda^2$. The first case is absurd as $\lambda \neq 1$, so $\lambda$ is an order $3$ character of $\Gamma^0$ (and ${\rm char}\, k \neq 3$). The image of $\Gamma \rightarrow {\rm GL}(H)$ is isomorphic to $\got{S}_3$, and its kernel is the subgroup $\Gamma^1:={\rm ker}\, \lambda \subset \Gamma^0$.  \ps\ps

As $\Gamma^0$ has index $2$ in $\Gamma$, and ${\rm dim}\, {\rm ad}\, P = 3$, the irreducibility of ${\rm ad}\, P$ implies the irreducibility of $({\rm ad}\, P)_{|\Gamma^0}$, hence of $P_{|\Gamma^0}$ as well. As $P^\ast \simeq P \otimes c$, it follows that we have $[E]\,= \,2\, [P] + 1+ \lambda + \lambda^{-1}$ in ${\rm R}_k(\Gamma^0)$, with $\lambda$ of order $3$ and $P$ irreducible: we are thus in {\it subcase (b)} of {\it Case 3.} of the proof of Proposition~\ref{lemmaindex3} (for the triple $(\Gamma^0,P,\lambda)$ instead of $(\Gamma,J,c)$), and in particular in case (iii) of the statement of that proposition. It implies that $\Gamma^0$ is ``the'' nontrivial central extension of $\got{A}_4$ by $\Z/2\Z$, and that is center ${\rm Z}(\Gamma^0)$ is the subgroup of order $2$ in ${\rm SO}(E)$ acting by $\pm {\rm id}$ on $P \oplus P^\ast$ and by ${\rm id}$ on $V$. This description shows that ${\rm Z}(\Gamma^0)$ lies in the center ${\rm Z}(\Gamma)$ of $\Gamma$, hence ${\rm Z}(\Gamma^0)={\rm Z}(\Gamma)$ ($\simeq \Z/2\Z$). \ps\ps

Let $\Gamma^2 \subset \Gamma$ be the kernel of the morphism $\Gamma \rightarrow {\rm GL}({\rm ad} P)$, we have ${\rm Z}(\Gamma) \subset \Gamma^2$. The second equation of \eqref{eqind2bisbis} shows $H_{|\Gamma^2}  \simeq 1 \oplus c$, the subspace of $\Gamma^2$-invariants in $H$ is thus nonzero. But this subspace is stable by $\Gamma$ as $\Gamma^2$ is a normal subgroup, it follows that $\Gamma^2$ acts trivially in $H$, {\it i.e.} we have $\Gamma^2 \subset \Gamma^1$. The group $\Gamma^2$ acts thus by scalars of determinant $1$ in $P$, so we have $\Gamma^2 = {\rm Z}(\Gamma)$. Moreover, the finite group $\Gamma/\Gamma^2$ is isomorphic to a subgroup of ${\rm SO}({\rm ad} \,P)$;  it contains $\Gamma^0/\Gamma^2$, which is a subgroup of index $2$ isomorphic to $\got{A}_4$ and acts irreducibly on ${\rm ad}\, P$, hence with a trivial centralizer in ${\rm SO}({\rm ad} \,P)$. As $\got{A}_4$ has an automorphism group isomorphic to $\got{S}_4$, this forces $\Gamma/\Gamma^2$ itself to be isomorphic to $\got{S}_4$. We have thus proved that $\Gamma$ is a central extension of $\got{S}_4$ by $\Z/2\Z$ which does not split over $\got{A}_4$. \ps\ps
We know from Schur that, up to isomorphism, there are exactly two central extensions $G$ of $\got{S}_4$ by $\Z/2\Z$ which do not split over $\got{A}_4$: we have either $G \simeq \GL_2(\Z/3\Z)$, or $G$ is isomorphic to the group usually denoted $\widetilde{\got{S}_4}$ (see e.g. \cite[\S 9.1.3]{serregalois}). In any case, character theory shows that such a $G$ has exactly $3$ isomorphism classes of irreducible $k$-linear representations which are nontrivial on the center $\Z/2\Z$: they have the form $J, J\otimes c$ and $J \otimes H$ where $\dim J = 2$, $c$ is the non-trivial (order $2$) character of $\got{S}_4$, and $H$ is the $2$-dimensional irreducible representation of $\got{S}_4$ (inflated from a surjective homomorphism $\got{S}_4 \rightarrow \got{S}_3$). The groups $\GL_2(\Z/3\Z)$ and $\widetilde{\got{S}_4}$ are distinguished as follows: we have $\det J = c$ in the first case (and thus $J^\ast \simeq J \otimes c$), and $\det J =1$ in the second (and then $J^\ast \simeq J$). As a consequence, we have $\Gamma \simeq \GL_2(\Z/3\Z)$ and assertion (iv) holds. \end{pf}
\ps\ps

\begin{pf} Let us prove Theorem \ref{mainthmb} when $\Gamma$ has Witt index $2$. An immediate consequence of Proposition \ref{rhogl23} is that we are in case (iv) of Proposition \ref{lemmaindex2} if, and only if, we are in case (b) of the statement of Theorem \ref{mainthmb} with $\Gamma=\rho({\rm GL}_2(\Z/3\Z))$. It only remains to show that we are in case (iii) of the proposition if, and only if, we are in case (d) of Theorem \ref{mainthmb} with $\epsilon(\Gamma')=\{\pm 1\}$. But this follows from a similar argument as the one given in the last paragraph of \S \ref{parcasewindex3}, by Example \ref{exempled8} and the following lemma. Note that the characters $c$ and $\epsilon$ here play respectively the role of the characters $\epsilon$ and $\epsilon\mu$ of subgroups of ${\rm O}_2^\pm(k)$.\end{pf}
\ps\ps

\begin{lemme} Let $\Gamma \subset {\rm O}_2^\pm(k)$ be a subgroup on which both $\mu$ and $\epsilon$ are nontrivial. Then $\Gamma$ is isomorphic to ${\rm D}_8$ if, and only if, the $k[\Gamma]$-module ${\rm H}(k)$ is both irreducible and selfdual. 
\end{lemme}
\ps\ps

\begin{pf} We have already seen in Example \ref{exempled8} that if $\Gamma$ is isomorphic to ${\rm D}_8$ then ${\rm H}(k)$ is both irreducible and selfdual. Assume conversely that the $k[\Gamma]$-module $P={\rm H}(k)$ has these two properties. By irreducibility of $P$, $\epsilon\mu$ is nontrivial on $\Gamma$, and thus $\epsilon,\mu$ and $\epsilon\mu$ are of order $2$, hence distinct. As $P$ is selfdual  we have $P^\ast \simeq P \otimes \chi$ for all $\chi$ in $\{1,\epsilon,\mu,\epsilon\mu\}$. In particular ${\rm Sym}^2\,P$ is isomorphic to $1 \oplus \mu \oplus \epsilon\mu$, which implies that $|\Gamma|$ divides $8$. As $P$ is irreducible $\Gamma$ is nonabelian, so we have $|\Gamma|=8$ and $P$ is the unique $2$-dimensional irreducible representation of $\Gamma$ up to isomorphism. But if $\Gamma$ was the ``quaternion group'' of order $8$ we would have $\det P = 1$. So $\Gamma$ is isomorphic to ${\rm D}_8$, and we are done. 
\end{pf}

\subsection{The case of Witt index $1$}

We shall now prove that there is no subgroup $\Gamma$ of ${\rm SO}(E)$ of Witt index $1$ which satisfies assumptions (i) and (ii) of Theorem \ref{mainthmb}. We start with a lemma on representations of Witt index $0$. \ps\ps

 If $V$ is a finite dimensional $k$-linear representation of $\Gamma$, we shall denote by $V^\Gamma$ the subspace $\{v \in V \, |\, \, \gamma v = v \, \, \,\,\forall \gamma \in \Gamma\}$ of $\Gamma$-invariants in $V$, and we set ${\rm t}(V) = \dim_k V^\Gamma$. If $U,V$ are finite dimensional $k$-linear representations of $\Gamma$, the natural isomorphism $U \otimes V \isomo {\rm Hom}(U^\ast,V)$ shows the equality $${\rm t}(U \otimes V) = \dim_k {\rm Hom}_{k[\Gamma]}(U^\ast,V),$$
 that we shall freely use below. It implies that if both $U$ and $V$ are irreducible, then we have ${\rm t}(U \otimes V) \leq 1$. \ps\ps

\begin{lemme} \label{lemmaind0} Let $V$ be a nondegenerate quadratic space over $k$, $\Gamma$ a group, $\rho : \Gamma \longrightarrow {\rm O}(V)$ a semisimple representation of Witt index $0$, $r$ the number of irreducible summands of the $k[\Gamma]$-module $V$ and $c : \Gamma \rightarrow k^\times$ a group homomorphism. \ps\ps\begin{itemize}
\item[(i)] We have ${\rm t}(\Lambda^2 \, V)=0$ and ${\rm t}({\rm Sym}^2 V)=r$.\ps\ps
\item[(ii)] We have ${\rm t}(\Lambda^2\, V \otimes c^{-1}) \leq r $ and ${\rm t}(\Lambda^2\, V \otimes c^{-1}) \leq \frac{\dim V}{2}$.\ps\ps
\item[(iii)] If $\dim V$ is odd and ${\rm t}(\Lambda^2 \,V \otimes c^{-1}) \neq 0$ then we have $r\geq 2$. \ps \ps

\item[(iv)] If $\dim V \leq 5$ then we have ${\rm t}(\Lambda^2 \,V \otimes c^{-1}) \leq 2$. \ps\ps
 \end{itemize}
 
Assume furthermore $4 \leq \dim V \leq 5$ and ${\rm t}(\Lambda^2 \,V \otimes c^{-1})= 2$. \begin{itemize}\ps \ps
 
\item[(v)] We have $c^2 =1$ and a $k[\Gamma]$-module decomposition 
\begin{equation}\label{casdeg2} V =  U \oplus U' \oplus V' \end{equation}
with $\dim U = \dim U'=2$ and $\det U = \det U' = c$. \ps\ps
\item[(vi)] Assume $r<\dim V$ and ${\rm t}(\Lambda^2 V \otimes \chi^{-1})\neq 0$ for some character $\chi : \Gamma \rightarrow k^\times$ with $\chi \neq c$. Then we have ${\rm t}(\Lambda^2 V \otimes \chi^{-1})=1$, and if $r=\dim V -2$ then we have 
$U' \simeq U \otimes \chi$ in \eqref{casdeg2}. \ps\ps
\end{itemize}

\end{lemme}

\begin{pf} We may write $V = \oplus_{i \in I} V_i$ where the $V_i$ are irreducible, selfdual and pairwise distinct, and with $r=|I|$. For each $i$ in $I$, the unique $\Gamma$-invariant line of ${\rm End}\, V_i \simeq V_i \otimes V_i \simeq {\rm Sym}^2 \,V_i \oplus \Lambda^2 \,V_i$ lies in ${\rm Sym}^2 \,V_i$, so we have ${\rm t}(\Lambda^2\, V_i) = 0$ and ${\rm t}({\rm Sym}^2 \,V_i)=1$. Moreover, if $i\in I$ we have $V_i \simeq V_i^\ast$ so ${\rm t}(V_i \otimes V_j) = 0$ if $i \neq j$. Choose a total ordering $<$ on $I$.  We have a $k[\Gamma]$-linear isomorphism $$\Lambda^2\, V \simeq \left( \,\bigoplus_{i \in I} \Lambda^2 \,V_i\, \right) \oplus \left( \bigoplus_{i<j} V_i \otimes V_j \right),$$
as well as a similar isomorphism with each $\Lambda^2$ replaced by ${\rm Sym}^2$. This proves assertion (i). \ps\ps

Let $c : \Gamma \rightarrow k^\times$ be a character and consider the sets  
$${\rm I}(c) = \{ i \in I \, \,\, |\, \, \,  {\rm t}(\Lambda^2 \, V_i \otimes c^{-1}) \neq 0\}$$ $${\rm J}(c) = \{ (i,j) \in I \times I\, \, \, |\, \, \,  i<j, \,\,\,{\rm t}(V_i \otimes V_j \otimes c^{-1}) \neq 0\}.$$\ps\ps

Recall that $\Lambda^2 V_i$ injects into $V_i \otimes V_i$, and that ${\rm t}(V_i \otimes V_j) \leq 1$ for all $i,j$ by irreducibility of the $V_j$, so that we also have ${\rm t}(\Lambda^2 \, V_i \otimes c^{-1}) \leq 1$. In particular, we have the equality 
\begin{equation}\label{formtij} {\rm t}(\Lambda^2 \,V \otimes c^{-1}) = |{\rm I}(c)|+|{\rm J}(c)|.\end{equation}
We now make a couple of simple observations. \ps\ps

-- {\rm (O1)} if $i \in {\rm I}(c)$ then $\dim V_i$ is even, and we have $\det V_i = c^{\frac{\dim V_i}{2}}$ and $c^{\dim V_i} = 1$. In particular, $\det V_i = c$ if $\dim V_i=2$.\ps\ps

Indeed, if $i \in {\rm I}(c)$ then there is a nonzero alternating pairing on $V_i$ such that $\Gamma$ acts by symplectic similitudes of factor $c$. Such a pairing is necessarily nondegenerate as $V_i$ is irreducible, so ${\rm (O1)}$ follows from general properties of similitude symplectic groups. As a consequence, if $\dim V$ is odd and $V$ is irreducible, then $\Lambda^2 V$ does not contain any character, which proves assertion (iii). \ps\ps

--  {\rm (O2)} If $(i,j) \in {\rm J}(c)$ then $\dim V_i = \dim V_j$ and $c^{2 \dim V_i} = 1$. \ps\ps

Indeed, if $(i,j) \in {\rm J}(c)$ then there is an isomorphism $V_i \simeq V_j \otimes c$.\ps\ps

--  {\rm (O3)} if $i \in {\rm I}(c)$ then there is no $j \in I$ such that $(i,j) \in {\rm J}(c)$ or $(j,i) \in {\rm J}(c)$ (since $V_i \simeq V_i \otimes c \simeq V_i \otimes c^{-1}$). \ps\ps

--  {\rm (O4)} if $(i,j) \in {\rm J}(c)$ then $(r,s) \in {\rm J}(c)$ implies $\{i,j\} \cap \{r,s\} = \emptyset$ or $\{i,j\} = \{r,s\}$ (observe that if $i,j \in I$ and $V_i \simeq V_j \otimes c$ then $V_j \simeq V_i \otimes c$). \ps\ps

It follows from those observations that we have the inequalities 
$$\begin{array}{ccc} \dim V & \geq & \sum_{i \in {\rm I}(c)} \dim V_i + 2\,  \sum_{(i,j) \in {\rm J}(c)} \dim V_i \\ \\
&  \geq & 2 \,|{\rm I}(c)| + 2\,|{\rm J}(c)| = 2\, {\rm t}(\Lambda^2 \, V \otimes c^{-1}), \,\,\,\,{\rm and}\end{array}$$\ps
$$r \geq |{\rm I}(c)| + 2\,|{\rm J}(c)| = |{\rm J}(c)|+{\rm t}(\Lambda^2 \, V \otimes c^{-1}) \geq {\rm t}(\Lambda^2 \, V \otimes c^{-1}).$$ 
\ps\ps

This proves assertions (ii). Assertion (iv) is a special case of assertion (ii). \ps\ps

Assume now $4 \leq \dim V \leq 5$ and ${\rm t}(\Lambda^2 \, V \otimes c^{-1}) = 2$. We have $2 = |{\rm I}(c)|+|{\rm J}(c)|$ by \eqref{formtij}. There are three cases: \begin{itemize}\ps\ps

\item[(1)] $|{\rm I}(c)|=2$ and $|{\rm J}(c)|=0$. Then {\rm (O1)} shows $\dim V_i=2$ for $i \in {\rm I}(c)$, $c^2=1$, hence assertion (v) holds with $U$ and $U'$ irreducible of dimension $2$ and $r = \dim V -2$. \ps\ps

\item[(2)] $|{\rm J}(c)|=2$ and $|{\rm I}(c)|=0$. Then {\rm (O2)} and {\rm (O4)} show $\dim V_i=\dim V_j=1$ and $\det (V_i \oplus V_j) = V_i \otimes V_j = c$ if $(i,j) \in {\rm J}(c)$, $c^2=1$, and $V = V' \oplus \bigoplus_{(i,j) \in {\rm J}(c)} (V_i \oplus V_j)$. So assertion (v) holds with $U$ and $U'$ reducible of dimension $2$ and $r = \dim V$.  \ps\ps

\item[(3)] $|{\rm I}(c)|=|{\rm J}(c)|=1$. Then  {\rm (O1)}, {\rm (O2)} and {\rm (O3)} show $\dim V_i=2$ for $i \in {\rm I}(c)$ as well as $\dim V_i = \dim V_j = 1$ and $\det (V_i \oplus V_j)=c$ if $(i,j) \in {\rm J}(c)$. So assertion (v) holds with either $U$ or $U'$ reducible (but not both) and $r = \dim V-1$. \end{itemize}

Let now $\chi : \Gamma \rightarrow k^\times$ be a character different from $c$. We have ${\rm I}(\chi)=\emptyset$ by \eqref{casdeg2} and {\rm (O1)}, so ${\rm t}(\Lambda^2 \, V \otimes \chi^{-1})=|{\rm J}(\chi)|$. By {\rm (O2)} and {\rm (O4)}, we have $|{\rm J}(\chi)|\leq 2$ and if the equality holds then $V$ contains at least $4$ distinct characters, {\it i.e} $r=\dim V$. If we have $r=\dim V-2$, then $U$ and $U'$ are irreducible of determinant $\neq \chi$ in \eqref{casdeg2}.  Thus ${\rm J}(\chi) \neq \emptyset$ if, and only if, we have $U \simeq U' \otimes \chi$. This proves assertion (vi). \end{pf}
\ps\ps
\begin{prop}\label{lemmaindex1} Let $E$ be a $7$-dimensional nondegenerate quadratic space over $k$. Then there is no subgroup $\Gamma \subset {\rm SO}(E)$ of Witt index $1$ satisfying assumptions (i) and (ii) of the statement of Theorem \ref{mainthmb}.
\end{prop}

\begin{pf} Let $\Gamma \subset {\rm SO}(E)$ be of Witt index $1$ and satisfying assumptions (i) and (ii) of the statement of Theorem \ref{mainthmb}.
There is a character $\chi : \Gamma \rightarrow k^\ast$, as well as a representation $\Gamma \rightarrow {\rm SO}(V)$ of dimension $5$ and Witt index $0$, such that 
$E \simeq \,\chi \oplus \chi^{-1} \,\oplus \, V $. As $V$ is selfdual and $\det V=1$, we have $\Lambda^3 V \simeq \Lambda^2 V$, and one easily checks that equation \eqref{eqfundg2} 
is equivalent to the following identity
{\small \begin{equation} \label{eqind1} \Lambda^2 V  \, \otimes (1 \oplus \chi \oplus \chi^{-1}) \simeq 1 \oplus \chi \oplus \chi^{-1} \oplus \chi^2 \oplus \chi^{-2} \oplus {\rm Sym}^2 V \oplus V \otimes (\chi \oplus \chi^{-1}).\end{equation}}
Denote by $r$ the number of irreducible constituents of $V$. By Lemma~\ref{lemmaind0} (i) and (iv), we have ${\rm t}({\rm Sym}^2\, V)=r$, ${\rm t}(\Lambda^2\,V)=0$ and ${\rm t}(\Lambda^2 V \otimes \chi^{\pm 1}) \leq 2$. As $V\otimes \chi$ and $\Lambda^2 V \otimes \chi$ are semisimple, since so are the representations occurring in equation \eqref{eqind1} by (S), and as $V$ is selfdual, we also have ${\rm t}(\Lambda^2 V \otimes \chi)={\rm t}(\Lambda^2 V \otimes \chi^{-1})$ and ${\rm t}(V \otimes \chi)={\rm t}(V \otimes \chi^{-1})$. We obtain $$2\, {\rm t}(\Lambda^2 V \otimes \chi) = 1 + r + 2\, ({\rm t}(\chi)+{\rm t}(\chi^2)+{\rm t}(V \otimes \chi)).$$
We deduce from this identity: ${\rm t}(\Lambda^2 V \otimes \chi) >0$, $r$ is odd and $1+r+2 \, {\rm t}(\chi^2) \leq2\, {\rm t}(\Lambda^2 V \otimes \chi) \leq  4$. But Lemma \ref{lemmaind0} (iii) shows $r\geq 2$. The only possibility is thus ${\rm t}(\Lambda^2 V \otimes \chi)=2$, $r=3$ and $\chi^2 \neq 1$, which contradicts Lemma \ref{lemmaind0} (v). \end{pf}

\subsection{The case of Witt index $0$} \label{parcasewindex0}

In order to prove Theorem \ref{mainthmb}, it only remains to prove the following Theorem.\ps\ps

\begin{thm} \label{thmindex0} Let $E$ be a $7$-dimensional nondegenerate quadratic space over $k$ and $\Gamma \subset {\rm SO}(E)$ a subgroup of Witt index $0$ satisfying assumptions (i) and (ii) of the statement of Theorem \ref{mainthmb}. Then $\Gamma$ is contained in a ${\rm G}_2$-subgroup of ${\rm SO}(E)$.
\end{thm}
\ps\ps
We shall divide the proof into several parts. \ps\ps

\begin{prop}\label{propso4} For $i=3,4$, let $V_i$ be a nondegenerate quadratic space over $k$ of dimension $i$, set $V = V_3 \bot V_4$, and let $\Gamma \subset {\rm SO}(V)$ be a subgroup of Witt index $0$ such that $\Gamma \subset {\rm SO}(V_3) \times {\rm SO}(V_4)$. Then $\Gamma$ is contained in a ${\rm G}_2$-subgroup of ${\rm SO}(V)$ if, and only if, there is an injective $k[\Gamma]$-linear morphism $V_3 \rightarrow \Lambda^2 V_4$. \ps
If the $k[\Gamma]$-module $V_3$ has an irreducible summand of dimension $\geq 2$ which injects into $\Lambda^2 V_4$, then these properties are satisfied.
\end{prop}

\begin{pf} Let $W$, $W_3$ and $W_4$ be spinor modules for $V$, $V_3$ and $V_4$ respectively, and write $W_4 = W_4^+ \oplus W_4^{-}$ as the direct sum of its two half-spinor modules $W_4^{\pm}$. Recall from the sequence \eqref{rhoUV} that we have a natural homomorphism $\rho_{V_3;V_4} : {\rm GSpin}(V_3) \times {\rm GSpin}(V_4) \rightarrow {\rm GSpin}(V)$, whose image is the inverse image of ${\rm SO}(V_3) \times {\rm SO}(V_4)$ under the natural homomorphism $\pi_V : {\rm GSpin}(V) \rightarrow {\rm SO}(V)$. We may thus view $W$, $W_3$ and $W_4^{\pm}$ as $k[{\rm GSpin}(V_3) \times {\rm GSpin}(V_4)]$-modules. As such, there is an isomorphism \begin{equation}\label{f0so4}W \simeq W_3 \otimes (W_4^+\oplus W_4^-)\end{equation} by Proposition \ref{decospin}. As is well-known, the image of the natural injective homomorphism ${\rm GSpin}(V_4) \rightarrow \GL(W_4^+) \times \GL(W_4^-)$ is the subgroup of elements $(g^+,g^-) \in\GL(W_4^+) \times \GL(W_4^-)$ such that $\det g^+=\det g^-$. Moreover, as a representation of ${\rm GSpin}(V_4)$ we have 
\begin{equation}\label{f1so4} \det W_4^{\pm} = \nu_{V_4} \hspace{1 cm} {\rm and} \hspace{1cm} W_4^+ \otimes W_4^{-}  \simeq V_4 \otimes \nu_{V_4}.\end{equation} Besides, as already explained in the first paragraph of the proof of Proposition \ref{lemmaindex2}, the natural morphism ${\rm GSpin}(V_3) \rightarrow {\rm GL}(W_3)$ is an isomorphism, and as a representation of ${\rm GSpin}(V_3)$ we have 
\begin{equation}\label{f2so4}\det W_3 = \nu_{V_3} \hspace{1 cm} {\rm and} \hspace{1cm} {\rm ad}\,  W_3 \simeq V_3.\end{equation}\ps

Set $\widetilde{\Gamma} = \pi_V^{-1}(\Gamma) \cap {\rm Spin}(E)$. We may view $W, W_3, W_4^{\pm 1}, V, \nu_V, V_3, \nu_{V_3}$, $V_4, \nu_{V_4}$ as a representations of $\widetilde{\Gamma}$, and as such we have $\nu_{V_3}\nu_{V_4} = \nu_V =1$. We set $\nu = {\nu_{V_4}}_{|\widetilde{\Gamma}}$, we thus have $\det W_4^{\pm } = \nu$ and $\det W_3=\nu^{-1}$. \ps\ps

We claim that $W$, $W_3$ and $W_4^{\pm}$ are semisimple $k[\widetilde{\Gamma}]$-modules. We apply for this the main result of \cite{serre} as well as \cite[Thm. 2.4]{serre2}, using ${\rm char}\, k \neq 2$. Indeed, the semisimplicity of $V_3$ and $V_4$ implies that of $W_3 \otimes W_3$, $W_3$ and $W_4^{\pm}$ by \eqref{f2so4}, \eqref{f1so4} and \cite{serre2}, which in turn implies the semisimplicity of $W$ by \eqref{f0so4} and \cite{serre}.\ps\ps

Let us now check that $W_3$ and $W_4^{\pm}$ are simple $k[\widetilde{\Gamma}]$-modules. If $W_3$ is reducible, then we have $W_3 \simeq \chi \oplus \nu^{-1} \chi^{-1}$ for some character $\chi : \widetilde{\Gamma} \rightarrow k^\times$ as $W_3$ is semisimple. This implies $V_3 \simeq {\rm ad} \, W_3 \simeq 1 \oplus \chi^2 \nu \oplus \chi^{-2}\nu^{-1}$, a contradiction as $V_3$ has Witt index $0$ as a representation of $\Gamma$. Assume now that $W_4^+$ is reducible (the argument for $W_4^-$ will be similar), and write $W_4^+ \simeq \chi \oplus \chi^{-1} \nu$ as above. This implies $V_4 \simeq \chi \nu^{-1} \otimes W_4^- \oplus \chi^{-1} \otimes W_4^-$. Observe that the dual of $\chi \nu^{-1} \otimes W_4^-$ is $\chi^{-1} W_4^{-}$ as $\det W_4^- = \nu$. It follows that $V_4$ has Witt index $2$ as a representation of $\Gamma$, a contradiction. \ps\ps

By the same argument as the one given in the proof of Proposition \ref{lemmaindex2}, page \pageref{pp1}, we obtain the equivalences between: \begin{itemize} \ps\ps
\item[(a)] $\Gamma$ is contained in a ${\rm G}_2$-subgroup of ${\rm SO}(V)$, \ps\ps
\item[(b)] there exists a character $\beta : \widetilde{\Gamma} \rightarrow k^\times$ with $\beta^2=1$ such that $\beta$ is a summand of $W_{|\widetilde{\Gamma}}$, or which is the same, such that $W_3^\ast$ is isomorphic to $W_4^+ \otimes \beta$ or to $W_4^- \otimes \beta$,\ps\ps
\item[(c)] we have ${\rm ad}\, W_3 \simeq {\rm ad}\, W_4^+$ or ${\rm ad}\, W_3 \simeq {\rm ad}\, W_4^-$ as $k[\widetilde{\Gamma}]$-modules. \ps\ps
\end{itemize}

As we have $V_4 \simeq W_4^+ \otimes W_4^{-} \otimes \nu^{-1}$, we obtain a natural $k[\widetilde{\Gamma}]$-linear isomorphism
$\Lambda^2 \, V_4 \simeq \Lambda^2 W_4^+ \otimes {\rm Sym}^2 W_4^- \otimes \nu^{-2} \oplus \Lambda^2 W_4^- \otimes {\rm Sym}^2 W_4^+ \otimes \nu^{-2}$, which also writes as 
\begin{equation} \label{finpropso4} \Lambda^2 \, V_4 \simeq {\rm ad}\, W_4^- \oplus {\rm ad}\, W_4^+.\end{equation}
We are now ready to prove the proposition. Assume that $\Gamma$ is contained in a ${\rm G}_2$-subgroup of ${\rm SO}(V)$. Assertion (c) above, relation \eqref{f2so4} and the isomorphism \eqref{finpropso4} show that $V_3$ embeds in $\Lambda^2 \, V_4$. \ps\ps

Conversely, assume either that $V_3$ embeds in $\Lambda^2 \, V_4 \simeq  {\rm ad}\, W_4^+ \oplus {\rm ad}\, W_4^-$, or that $V_3$ has an irreducible summand of dimension $\geq 2$ which injects into $\Lambda^2 V_4$. There exists in both cases a $k[\Gamma]$-submodule $H$ of $V_3$ of dimension $\geq 2$ such that $H$ embeds either in ${\rm ad}\, W_4^+$ or in ${\rm ad}\, W_4^-$. 
As both $V_3$ and ${\rm ad}\, W_4^+$ are semisimple of dimension $3$ and determinant $1$, this implies that $V_3$ is isomorphic either to ${\rm ad}\, W_4^+$ or to ${\rm ad}\, W_4^-$, thus (c) holds.  \end{pf}

\begin{cor}\label{propso4cor}  Let $V$ be a $7$-dimensional nondegenerate quadratic space over $k$ and $\Gamma \subset {\rm SO}(V)$ a subgroup of Witt index $0$ satisfying assumptions (i) and (ii) of the statement of Theorem \ref{mainthmb}. Assume furthermore that the $k[\Gamma]$-module $V$ has a summand of dimension $3$ (or $4$) and determinant $1$. Then $\Gamma$ is contained in a ${\rm G}_2$-subgroup of ${\rm SO}(V)$.
\end{cor}

\begin{pf} We have a $k[\Gamma]$-module decomposition $V \simeq V_3 \oplus V_4$ with $\dim V_i = i$ and $\det V_i=1$. As the $V_i$ are selfdual of determinant $1$ we easily checks that equation \eqref{eqfundg2} simplifies as 
$$1  \oplus \Lambda^2 \,V_4 \otimes V_3 \simeq V_3\oplus {\rm Sym}^2 \,V_3 \oplus {\rm Sym}^2 \,V_4.$$
Let $r_i$ be the number of irreducible summands of $V_i$. By Lemma \ref{lemmaind0} (i), we have ${\rm t}(V_3)  = {\rm t}(\Lambda^2 V_3) = 0$  and 
$${\rm t}(\Lambda^2 \,V_4 \otimes V_3) = r_3+r_4 -1.$$
In particular ${\rm t}(\Lambda^2 \,V_4 \otimes V_3)>0$, so there is a nonzero $k[\Gamma]$-linear morphism $V_3 \rightarrow \Lambda^2 \,V_4$. 
If $V_3$ is irreducible we are done by Proposition \ref{propso4}. \ps\ps

Consider a character $\epsilon : \Gamma \rightarrow k^\ast$ and assume $V_3 \simeq \epsilon \oplus H$ with $H$ irreducible of dimension $2$ and determinant $\epsilon$. If $H$ occurs in $\Lambda^2 V_4$ then we are done again by the last assertion of Proposition \ref{propso4}.
We may thus assume ${\rm t}(\Lambda^2 \,V_4 \otimes H)=0$. We have then ${\rm t}(\Lambda^2 \,V_4 \otimes V_3) = {\rm t}(\Lambda^2 \,V_4 \otimes \epsilon) = 1+ r_4  \geq 2$. But thus is a contradiction as by Lemma \ref{lemmaind0} (ii), we have ${\rm t}(\Lambda^2 V_4 \otimes \epsilon) \leq 2$, with $r_4\geq 2$ if the equality holds. \ps\ps 

We may thus assume $V_3 = \oplus_{i=1}^3 \epsilon_i$ for some {\it distinct} characters $\epsilon_i$. Then:
$$\sum_{i=1}^3 {\rm t}(\Lambda^2 \,V_4 \otimes \epsilon_i) = 2 + r_4 \geq 3.$$
If each term of the sum on the left is nonzero then $V_3$ embeds in $\Lambda^2\, V_4$ so we are done by Proposition~\ref{propso4}, thus we may assume that one of them is $0$. If one term is $\geq 2$ then it is equal to $2$ and we have $r_4 \geq 2$ by  Lemma \ref{lemmaind0} (ii), so the two nonzero  $ {\rm t}(\Lambda^2 \,V_4 \otimes \epsilon_i)$ are equal to $2$ and we have $r_4=2$, which contradicts Lemma \ref{lemmaind0} (vi). \end{pf}

\begin{lemme}\label{lemmaindex01} Let $E$ be a $7$-dimensional nondegenerate quadratic space over $k$ and $\Gamma \subset {\rm SO}(E)$ of Witt index $0$ satisfying assumptions (i) and (ii) of the statement of Theorem \ref{mainthmb}. Then the $k[\Gamma]$-module $E$ has no irreducible summand of dimension $5$ and at most $1$ irreducible summand of dimension $3$ and nontrivial determinant. 
\end{lemme}

\begin{pf} Assume first that we have a $\Gamma$-stable decomposition $E = P \bot V$ where $P$ has dimension $2$ and $V$ is irreducible of dimension $5$. The character $c=\det P$ satisfies $c^2=1$, and if we set $\Gamma^0 = {\rm ker}\, c$ then we have a natural morphism $\Gamma^0 \rightarrow {\rm SO}(V)$. In particular, there is a character $\chi : \Gamma^0 \rightarrow k^\times$ such that $P_{|\Gamma^0} \simeq \chi \oplus \chi^{-1}$. As $\Gamma^0$ has index $\leq 2$ in $\Gamma$, and as $\dim V = 5$ is odd, observe that $V_{|\Gamma^0}$ is irreducible. It follows that $\Gamma^0$ has Witt index $1$ and satisfies assumptions (i) and (ii) of Theorem \ref{mainthmb}. This impossible by Proposition \ref{lemmaindex1}. \ps \ps
 
 Assume we have $E \simeq \epsilon \oplus V \oplus V'$ with $V,V'$ distinct, irreducible, selfdual, of dimension $3$. Set $c= \det V$ and $c'=\det V'$, we assume that $c$ and $c'$ are nontrivial, so the relation $cc'\epsilon=1$ implies $c \neq \epsilon$ and $c'\neq \epsilon$. We have $\Lambda^2 V \simeq V \otimes c$, $\Lambda^2 V' \simeq V' \otimes c'$, thus $$\Lambda^3 E \simeq c \oplus c' \oplus V \otimes V' \otimes (\epsilon \oplus c \oplus c') \oplus V \otimes c' \oplus V' \otimes c.$$
Let $\chi \in \{c,c',\epsilon\}$ be such that ${\rm t}(V \otimes V' \otimes \chi) \neq 0$. Then we have $\chi^2=1$, $V \simeq V' \otimes \chi$ and
${\rm t}(V \otimes V' \otimes \chi) = 1$. Taking determinant shows $c=c'\chi^3=c'\chi$, so $\chi = \epsilon$. This implies ${\rm t}(\Lambda^3\,E) = {\rm t}(V \otimes V' \otimes \epsilon) \leq 1$.  On the other hand, equation \eqref{eqfundg2} and Lemma \ref{lemmaind0} (i) imply ${\rm t}(\Lambda^3\,E) = {\rm t}(\epsilon) + 3$, a contradiction.\end{pf}

\begin{lemme} Theorem \ref{thmindex0} holds when the $k[\Gamma]$-module $E$ contains at least two characters.
\end{lemme}

\begin{pf} Assume we have $E \simeq \epsilon_1 \oplus \epsilon_2 \oplus V$ for some characters $\epsilon_i : \Gamma \rightarrow k^\times$. As $\dim V = 5$, $V^\ast \simeq V$ and $\det V = \epsilon_1 \epsilon_2$, we have $\Lambda^3 V \simeq \Lambda^2\, V \otimes \epsilon_1\epsilon_2$, hence an isomorphism $\Lambda^3 \,E \simeq  \epsilon_1 \epsilon_2 \otimes V \oplus \Lambda^2\, V \otimes (\epsilon_1 \epsilon_2 \oplus \epsilon_1 \oplus \epsilon_2)$. If ${\rm t}(\epsilon_1 \epsilon_2 \otimes V ) \neq 0$, then $\epsilon_1\epsilon_2$ occurs in $V$, and $\epsilon_1 \oplus \epsilon_2 \oplus \epsilon_1 \epsilon_2$ is a $3$-dimensional summand of $E$ of determinant $1$, so we are done by Corollary \ref{propso4cor}. We may thus assume ${\rm t}(V \otimes \epsilon_1 \epsilon_2)=0$, so that we have
$${\rm t}(\Lambda^3 E) = {\rm t}(\Lambda^2\, V \otimes \epsilon_1 ) + {\rm t}(\Lambda^2\, V \otimes \epsilon_2 ) + {\rm t}(\Lambda^2\, V \otimes \epsilon_1 \epsilon_2 ),$$
and each term in the sum on the right is $\leq 2$ by Lemma \ref{lemmaind0} (iv). On the other hand we have ${\rm t}(\Lambda^3 E) ={\rm t}(E) + 2+ r$ where $r$ is the number of irreducible summands of $V$ by Equation \eqref{eqfundg2} and Lemma \ref{lemmaind0} (i). In particular, there exists a character $\chi$ of $\Gamma$ such that ${\rm t}(\Lambda^2\, V \otimes \chi) \neq 0$. This implies $r\geq 2$ by Lemma \ref{lemmaind0} (iii), so there even exists a character $\chi$ of $\Gamma$ such that ${\rm t}(\Lambda^2\, V \otimes \chi) = 2$. But this implies that $\det V$, that is $\epsilon_1 \epsilon_2$, is a summand of $V$ by the same lemma assertion (v) (note $V' \simeq \det V$ in the notations {\it loc. cit.}), which contradicts  ${\rm t}(V \otimes \epsilon_1 \epsilon_2)=0$.
\end{pf}

\begin{lemme} Theorem \ref{thmindex0} holds when the $k[\Gamma]$-module $E$ contains at least two irreducible summands of dimension $2$.
\end{lemme}

\begin{pf} Assume we have $E \simeq P_1 \oplus P_2 \oplus V$ where each $P_i$ is irreducible of dimension $2$. Set $c_i = \det P_i$. As $V$ is selfdual of dimension $3$ and determinant $c_1 c_2$ we easily check the isomorphism
$$\Lambda^3 E \simeq c_1c_2 \oplus P_1 \otimes c_2 \oplus P_2 \otimes c_2 \oplus V \otimes (c_1 \oplus c_2 \oplus P_1 \otimes c_1c_2 \oplus P_2 \otimes c_1c_2 \oplus P_1\otimes P_2).$$
By Corollary \ref{propso4cor} we may assume that $E$ has no summand of dimension $3$ and determinant $1$. In particular, we have
$c_1 c_2  \neq 1$ (consider $V$), ${\rm t}(V \otimes c_i) =0$ (consider $c_i \oplus P_i$) and ${\rm t}(V \otimes P_i \otimes c_1 c_2)=0$ otherwise we have $V \simeq P_i \otimes c_1c_2 \oplus c_j$ with $j\neq i$, and thus ${\rm t}(V \otimes c_j) \neq 0$. We obtain 
$${\rm t}(\Lambda^3 E) =  {\rm t}(V \otimes P_1\otimes P_2).$$ 
On the other hand we have ${\rm t}(\Lambda^3 E) ={\rm t}(V) + 2+ r \geq 3$, where $r$ is the number of irreducible summands of $V$. As $P_1$ is irreducible of dimension $2$ and $\dim V \otimes P_2 = 6$, we have ${\rm t}(V \otimes P_1\otimes P_2) = \dim {\rm Hom}_{k[\Gamma]}(P_1,V\otimes P_2) \leq 3$. It follows that $r=1$, {\it i.e. } $V$ is irreducible. As $\dim P_1 \otimes P_2 =4$, we thus have $3 \leq {\rm t}(V \otimes P_1\otimes P_2)= \dim {\rm Hom}_{k[\Gamma]}(V,P_1\otimes P_2) \leq 1$, a contradiction. \end{pf}

\begin{lemme} Theorem \ref{thmindex0} holds when the $k[\Gamma]$-module $E$ contains exactly one irreducible summand of each dimension $1, 2$ and $4$.
\end{lemme}

\begin{pf} Assume we have $E \simeq \epsilon \oplus P \oplus V$ where $P$ and $V$ are irreducible of respective dimension $2$ and $4$. Set $c = \det P$. It is an order $2$ character of $\Gamma$ as ${\rm SO}(P)$ acts reducibly on $P$. By Corollary \ref{propso4cor} we may assume $c\epsilon\neq 1$. We easily check $\Lambda^3 E \simeq c\epsilon \oplus V \otimes (\epsilon c  \oplus c \oplus P\otimes \epsilon ) \oplus \Lambda^2\, V \otimes (\epsilon \oplus P)$. By a standard argument (from now on) we thus get the equality
\begin{equation}\label{flemmpass} {\rm t}(\epsilon) + 3 = {\rm t}(\Lambda^2 V \otimes \epsilon) + {\rm t}(\Lambda^2 V \otimes P).\end{equation}
The number ${\rm t}(\Lambda^2 V \otimes P)$ is the multiplicity of $P$ as a submodule $\Lambda^2 V$. So we have ${\rm t}(\Lambda^2 V \otimes P) \leq \frac{\dim \Lambda^2 V}{\dim P} = 3$. As ${\rm t}(\Lambda^2 \,V )=0$ by Lemma \ref{lemmaind0} (i), the identity \eqref{flemmpass} implies $\epsilon \neq 1$. The irreducibility of $V$ implies ${\rm t}(\Lambda^2 V \otimes \epsilon)\leq 1$, so $P$ occurs at least twice in $\Lambda^2 V$. As $\det \Lambda^2 \,V = (\det V)^3 = \epsilon c$, and $\epsilon \neq 1$, we cannot have $\Lambda^2\, V \simeq P \oplus P \oplus P$. We have thus $\Lambda^2\, V \simeq P \oplus P \oplus \epsilon \oplus c$. \ps\ps

On the other hand the character $\epsilon$ is a summand of $E$, hence of $\Lambda^3\, E$ as well by equation \eqref{eqfundg2}. As $1$ does not occur in $\Lambda^2 V$, the only possibility is that $\epsilon$ is a summand  of $\Lambda^2\, V \otimes P$. But we have 
$\Lambda^2\, V \otimes P \simeq P \otimes (P \oplus P) \oplus P \otimes (\epsilon \oplus c)$. As a consequence, $\epsilon$ occurs in $P \otimes P \simeq c \oplus {\rm Sym}^2\, P$, hence in ${\rm Sym}^2\, P$. As ${\rm Sym}^2\, P$ is semisimple of determinant $c^3=c$, we have proved ${\rm Sym}^2\, P \simeq 1\oplus \epsilon \oplus \epsilon c$. \ps\ps

Let $\Gamma^0$ be the kernel of the natural morphism $\Gamma \rightarrow {\rm O}(P)$. What we have just proved about ${\rm Sym}^2\, P$ implies that $\Gamma/\Gamma^0$ is a (nonabelian) finite group of order $8$ and $\epsilon(\Gamma^0)=1$. Moreover, $\Gamma^0$ acts trivially on $\Lambda^2 V$, hence by scalars of square $1$ on $V$, so $|\Gamma^0|\leq 2$ and $|\Gamma|\leq 16$. This contradicts the inequality $|\Gamma| \geq 1+1+ (\dim P)^2 + (\dim V)^2 = 22$.\end{pf}

\begin{lemme} Theorem \ref{thmindex0} holds when the $k[\Gamma]$-module $E$ contains exactly one irreducible summand of each dimension $3$ and $4$.
\end{lemme}

\begin{pf} Assume we have $E \simeq V_3 \oplus V_4$ with each $V_i$ irreducible of dimension $i$. Set $c=\det V_3 = \det V_4$. We have an  isomorphism $\Lambda^3 \, E \simeq c \oplus V_4 \otimes c \oplus V_3 \otimes V_4 \otimes c \oplus V_3 \otimes \Lambda^2\, V_4$ and ${\rm t}(\Lambda^3 \, E)=2$. We thus get
	$$2 = {\rm t}(c) + {\rm t}(V_3 \otimes \Lambda^2\, V_4 ).$$
We have ${\rm t}(V_3 \otimes \Lambda^2\, V_4 )\leq \frac{\dim \Lambda^2 \, V_4}{\dim V_3}=2$. If $c \neq 1$ we obtain ${\rm t}(V_3 \otimes \Lambda^2\, V_4)=2$, hence $\Lambda^2\, V_4 \simeq V_3 \oplus V_3$. But this implies $\det \Lambda^2 \,V_4 = 1$, which contradicts $\det \Lambda^2 \,V_4 = (\det V_4)^3 = c^3 = c$. So we have $c  = \det \, V_3 = 1$ and we conclude by Corollary \ref{propso4cor}. 
\end{pf}

\begin{lemme}Theorem \ref{thmindex0} holds when the $k[\Gamma]$-module $E$ contains exactly one irreducible summand of each dimension $1$ and $6$.
\end{lemme}

\begin{pf} Assume we have $E \simeq \epsilon \oplus V$ with $V$ irreducible of dimension $6$. Let us show first that $\epsilon$ has order $2$. If $\epsilon = 1$ then we have ${\rm t}(E \oplus {\rm Sym}^2 E)={\rm t}(\Lambda^3\, E) = 3$ and $\Lambda^3 E \simeq \Lambda^3 V \oplus \Lambda^2 V$. As ${\rm t}(\Lambda^2 V) =0 $ we obtain ${\rm t}(\Lambda^3 V) = 3$. As $V$ is selfdual, the $k[\Gamma]$-module $\Lambda^3 V$ embeds in $V \otimes \Lambda^2 V$. In particular we have ${\rm t}(V \otimes \Lambda^2\,V) \geq 3$. This contradicts the inequality $\dim \Lambda^2 V = 15 < 3 \dim V$. \ps\ps

Set $\Gamma^0 = {\rm ker }\, \epsilon$. As $\epsilon$ has order $2$ this is a subgroup of index $2$ in $\Gamma$ and $E_{|\Gamma^0}$ still satisfies assumption (i) and (ii) of Theorem~\ref{mainthmb}. The previous paragraph thus shows that $V_{|\Gamma^0}$ is reducible. By Clifford theory it is thus a direct sum of two nonisomorphic irreducible $3$ dimensional representations of $\Gamma^0$. Those representations are not selfdual by the last assertion of Lemma \ref{lemmaindex01}. We have thus $V_{|\Gamma^0} \simeq W \oplus W^\ast$ for some irreducible, nonselfdual, representation $W$ of $\Gamma^0$.  As a consequence, $\Gamma^0$ has Witt index $3$, and we are in {\it Case 1.} of the proof of Proposition \ref{lemmaindex3}, whose conclusion is $\det W = 1$. We conclude the proof by the following lemma. \end{pf}

 Let $I$ be a finite dimensional $k$-vector space. Define ${\rm N}(I)$ as the subgroup of ${\rm O}({\rm H}(I))$ of elements $g$ such that the couple of subspaces $(g(I),g(I^\ast))$ is either $(I,I^\ast)$ or $(I^\ast,I)$. The action of ${\rm N}(I)$ on the set with two elements $\{I,I^\ast\}$ defines a group homomorphism ${\rm N}(I) \rightarrow \Z/2\Z$. The map $g \mapsto g_{|I}$ identifies the kernel of this morphism with the group $\GL(I)$, and the inverse image of $1 \in \Z/2\Z$ with the set of $k$-linear isomorphisms $I \rightarrow I^\ast$. If $\varphi$ is a symmetric isomorphism $I \rightarrow I^\ast$ and $g \in \GL(I)$ then $\varphi g^{-1} \varphi^{-1}$ is the adjoint of $g$ with respect to $\varphi$ and $\varphi$ has order $2$ in ${\rm N}(I)$. In particular we have ${\rm N}(I) \simeq {\rm GL}(I) \rtimes \Z/2\Z$. Set $E = {\rm H}(I) \bot k$. The morphism $g \mapsto g \bot \det g$ allows to see ${\rm O}({\rm H}(I))$ as a subgroup of ${\rm SO}(E)$. \ps
  
\begin{lemme} Let $I$ be a vector space of dimension $3$ over the algebraically closed field $k$ and $E = {\rm H}(I) \bot k$. Let ${\rm N}(I) \subset {\rm SO}(E)$ be as above and $\Gamma$ a subgroup of ${\rm N}(I)$ such that $\Gamma \cap {\rm \GL}(I) \subset {\rm SL}(I)$. Then $\Gamma$ is contained in a ${\rm G}_2$-subgroup of ${\rm SO}(E)$. 
\end{lemme}
\ps
\begin{pf} We have defined above natural inclusions $\GL(I) \subset {\rm N}(I) \subset {\rm SO}(E)$. Let $H$ be the inverse image of $\GL(I)$ under the surjective morphism $\pi_E : {\rm GSpin}(E) \rightarrow {\rm SO}(E)$ and let $W$ be a spinor module for $E$.  Recall from \S \ref{parclag} and Propositions \ref{cliffordhyp} and \ref{decospin} that there is a natural isomorphism $\epsilon : k^\times \times \GL(I) \isomo H$ such that the map ${\pi_E} \circ \epsilon : k^\times \times \GL(I) \rightarrow \GL(I)$ coincides with the projection onto the second factor and such that the restriction of $W$ to $H$ is isomorphic to the representation $\alpha \otimes \Lambda\,I$ where $\alpha : H \rightarrow k^\times$ is the character obtained as the composition of $\epsilon^{-1}$ and of the first projection $k^\times \times \GL(I) \rightarrow k^\times$. \ps\ps

Let $G$ be the inverse image of ${\rm N}(I)$ under $\pi_E$ and fix a symmetric isomorphism $\varphi : I \rightarrow I^\ast$. Then $H$ has index $2$ in $G$ and we may choose $\tau \in G-H$ such that the restriction to $I \subset E$ of the element $\pi_E(\tau) \in {\rm N}(I)$ coincides with $\varphi$. 
By the observations preceding the lemma there is a unique map $\psi : \GL(I) \rightarrow k^\times$ with $$ \tau \,\epsilon(\lambda, g) \,\tau^{-1}\, = \,\epsilon(\lambda \psi(g),{}^{\rm t}\!g^{-1})$$
for all $\epsilon(\lambda,g) \in H$, where ${}^{\rm t}\!g$ is the adjoint of $g$ with respect to $\varphi$. This implies that $\psi$ is a group homomorphism. Applying $\nu_E$ we obtain the relation $(\det \, \cdot \,\psi^{-1})^2=1$, hence $\psi = \det $ as $k$ is algebraically closed. \ps\ps

The $k[H]$-module isomorphism $W \simeq \alpha \otimes \Lambda\, I$ shows that there are exactly $2$ lines in $W$ which are stable by $H$, let us call $P$ their direct sum. Then $P$ is stable by $G$ and we have a $k[H]$-linear isomorphism $P \simeq \alpha \oplus \alpha \det I$. The previous paragraph shows that the outer conjugation of $G$ on $H$ exchanges the (distinct) characters $\alpha$ and $\alpha \det I$. As $\dim I=3$, recall from \S \ref{octoalg} that $W$ may be equipped with a structure of nondegenerate quadratic space over $k$ such that ${\rm GSpin}(E)$ acts on $W$ as orthogonal similitudes of factor $\nu_E$. By definition, $P$ is necessarily nondegenerate, and as neither $\alpha$ nor $\alpha \det I$ has square $1$ on $H$ it follows that the two stable lines of $P$ are its isotropic lines, and that $\tau$ exchanges them. \ps \ps

So far the group $\Gamma$ did not play any role. Set now $\widetilde{\Gamma} = \pi_E^{-1}(\Gamma) \cap {\rm Spin}(E)$ and $\widetilde{\Gamma}^0 = H \cap \widetilde{\Gamma}$. We may and will view $\alpha$ and $I$ as representations of $\widetilde{\Gamma}$.  We have a natural morphism $\rho : \widetilde{\Gamma} \rightarrow {\rm O}(P)$, with $P_{|\widetilde{\Gamma}^0} \simeq \alpha_{|\widetilde{\Gamma}} \oplus \alpha_{|\widetilde{\Gamma}} \det I $. But we also have the relation $\alpha^2 \det I = \nu_E = 1$ on $\widetilde{\Gamma}$, as well as $\det I = 1$ on $\widetilde{\Gamma}^0$ by assumption. As a consequence, $\widetilde{\Gamma}^0$ acts by $\pm {\rm id}$ on $P$. If $\widetilde{\Gamma}^0 \neq \widetilde{\Gamma}$, the previous paragraph shows $\rho(\widetilde{\Gamma}) \subsetneq {\rm SO}(P)$, thus the kernel of any orthogonal symmetry in $\rho(\widetilde{\Gamma})$ is a nondegenerate one-dimensional subspace of $P$ on which $\widetilde{\Gamma}$ acts via an order $2$ character. We conclude the proof by Proposition~\ref{propgspin7} (c). \end{pf}

\ps
This ends the proof of Theorem \ref{thmindex0}, hence of Theorem \ref{mainthmb}.
\ps
\section{Proof of Theorem \ref{thmintroc}}\label{pfthmc}

The existence, and uniqueness up to conjugacy, of the morphisms denoted by $\alpha, \beta, \gamma$ before the statement of Theorem \ref{thmintroc} are consequences of \S \ref{parexamples} and assertion (b) of the following lemma applied to $K={\rm SO}(7)$ and $G = {\rm SO}(E \otimes \C)$.

\begin{lemme}\label{lemmecompactcomplex}  Let $K$ be a compact connected Lie group and $\iota : K \hookrightarrow G$ a complexification of $K$. \begin{itemize}
\item[(a)] If we have $x,y \in K$ and $g \in G$ with $g\,\iota(x)\,g^{-1}\,=\,\iota(y)$, and if $g\,=\,\iota(k)\,s$ is the polar {\rm (}{\it i.e.} Cartan{\rm)} decomposition of $g$ with respect to $\iota(K)$, then we have $k \,x\, k^{-1}\, =\, y$. \ps
\item[(b)] Let  $\Gamma$ be a compact group and $\rho : \Gamma \rightarrow G$ a continuous group homomorphism. Then there exists a continuous group morphism $\rho' : \Gamma \rightarrow K$ such $\iota \circ \rho'$ is $G$-conjugate to $\rho$, and such a $\rho'$ is unique up to $K$-conjugacy. 
\end{itemize}
\end{lemme}

{\scriptsize \begin{pf} (see e.g. \cite[p.1]{serrelie}) In order to prove (a), we rewrite the hypothesis as an equality $\iota(kx)\,(\iota(x)^{-1}\,s\,\iota(x)) \,=\, \iota(yk)\, s$ and conclude by the uniqueness of the polar decomposition. Let us prove (b). The group $\rho(\Gamma)$ is a compact subgroup of $G$ and a well-known result of Cartan asserts that the maximal compact subgroups of $G$ are the conjugate of $\iota(K)$. This shows the existence of $\rho'$, and its uniqueness up to $K$-conjugacy follows from (a).\end{pf}}

Applying this lemma again both to the compact group ${\rm SO}(7)$ and its ${\rm G}_2$-subgroups we conclude that Theorem \ref{thmintroc} of the introduction follows from the case $k=\C$ of Theorem \ref{mainthmb}. $\square$

\section{Automorphic and Galois representations}\label{autgalrep}
\ps 

\subsection{}\label{notationsaut} Let $F$ be a number field, $\mathcal{O}_F$ its ring of integers, and $\AAA_F$ its ad\`ele ring. For each place $v$ of $F$ we denote by $F_v$ the associated completion of $F$ and $\mathcal{O}_v \subset F_v$ its ring of integers. \ps\ps

Let $G$ be a connected linear algebraic group defined and reductive over $\mathcal{O}_F$. We refer to \cite{boreljacquet} (see also \cite[\S 3]{cogdell} when $G={\rm GL}_n$) for the notion of {\it cuspidal automorphic representation} of $G(\AAA_F)$. Such a representation $\pi$ has, for each place $v$ of $F$, a {\it local component} $\pi_v$ which is irreducible and well-defined up to isomorphism. When $v$ is a finite place, $\pi_v$ is a complex {\it smooth representation} of $G(F_v)$. When $v$ is an Archimedean place, $\pi_v$ is a complex $(\mathfrak{g}_v,K_v)$-{\it module} in the sense of Harish-Chandra, where $K_v$ is a fixed maximal compact subgroup of the Lie group $G(F_v)$ and $\mathfrak{g}_v$ is the Lie algebra of $G(F_v)$. \ps \ps

Let $v$ be a finite place of $F$. Recall that a complex, irreducible, smooth representation $\varpi$ of the locally profinite group $G(F_v)$ is called {\it unramified} if its underlying space has a nonzero vector invariant by the compact open subgroup $G(\mathcal{O}_v)$. According to Satake and Langlands, the set of isomorphism classes of unramified representations of $G(F_v)$ is in canonical bijection with the set of conjugacy classes of semisimple elements in $\widehat{G}$ \cite{euler,borel}. Here, $\widehat{G}$ is a dual group for $G$ in the sense of Langlands, namely a reductive (connected) complex linear algebraic group whose based root datum is equipped with an isomorphism with the dual based root datum of $G$. In particular we may and shall take $\widehat{\GL_n}=\GL_n(\C)$. We shall denote by ${\rm c}(\varpi) \subset \widehat{G}$ the semisimple conjugacy class associated to $\varpi$. The unramified representation $\varpi$ of $G(F_v)$ is said {\it tempered} if ${\rm c}(\varpi)$ is the conjugacy class of an element lying in a compact subgroup of $\widehat{G}$ (viewed as a Lie group). \ps\ps

If $\pi$ is a cuspidal automorphic representation of $G(\AAA_F)$, then $\pi_v$ is unramified for all but finitely many finite places $v$ of $F$. Let ${\rm Ram}(\pi)$ denote the (finite) set of places $v$ of $F$ which are either Archimedean, or such that $v$ is finite and $\pi_v$ is not unramified. If $S$ is a finite set of places of $F$ containing ${\rm Ram}(\pi)$, we say that $\pi$ is $S$-tempered if $\pi_v$ is tempered for all $v \notin S$. We say that $\pi$ is tempered if it is $S$-tempered for some such $S$. Conjecturally, ``$S$-tempered'' implies ``${\rm Ram}(\pi)$-tempered''. \ps\ps

Let $\pi$ be a cuspidal automorphic representation of $G(\AAA_F)$, $S$ a finite set of places of $F$ containing ${\rm Ram}(\pi)$, and $\rho : \widehat{G} \rightarrow \GL_n(\C)$ a polynomial representation. Following Langlands \cite{euler} consider the Euler product
$${\rm L}^S(s,\pi,\rho) = \prod_{v \notin S} \frac{1}{\det(1- q_v^{-s} \rho({\rm c}(\pi_v)))},$$
where $q_v$ denotes the cardinality of the residue field of $\mathcal{O}_v$. If $\pi$ is $S$-tempered, this product is well-defined and absolutely convergent in the half plane ${\rm Re}\, s >1$. In general, Langlands has shown that ${\rm L}^S(s,\pi,\rho)$ is absolutely convergent for ${\rm Re}\, s$ big enough. \ps\ps

Let $v$ be a real place of $F$, {\it i.e.} $F_v=\R$. For any irreducible $(\mathfrak{g}_v,K_v)$-module $\varpi$, the center $\mathfrak{z}_v$ of the enveloping algebra of $\mathfrak{g}_v$ acts by scalars on $\varpi$. The resulting $\C$-algebra homomorphism $\mathfrak{z}_v \rightarrow \C$ is called the {\it infinitesimal character} of $\varpi$; it may be viewed according to Harish-Chandra and Langlands as a semisimple conjugacy class ${\rm c}(\varpi)$ in the complex Lie algebra $\widehat{\mathfrak{g}}$ of $\widehat{G}$ \cite{langlandspb,euler,borel}.  When $G = \GL_n$, in which case we have $\widehat{\mathfrak{g}} = {\rm M}_n(\C)$, we shall say that $\varpi$ is {\it algebraic} if the eigenvalues of ${\rm c}(\varpi)$ are integers, which are then called the {\it weights} of $\varpi$. An algebraic $\varpi$ is {\it regular} if its weights are distinct. \ps\ps

\subsection{} \label{algreg}Assume $F$ is totally real and $n\geq 1$ is an odd integer. Let $\pi$ be a cuspidal automorphic representation of ${\rm GL}_n(\AAA_F)$ such that $\pi_v$ is algebraic regular for each real place $v$ of $F$. Then there exists a number field $E \subset \C$, called {\it a coefficient field for $\pi$}, such that 
for each finite place $v$ of $F$ the representation $\pi_v$ is defined over $E$ \cite[\S 3.5]{clozel}. In particular, if $\pi_v$ is unramified then the characteristic polynomial $\det(t -{\rm c}(\pi_v))$ belongs to $E[t]$. Let $\ell$ be a prime number, $\lambda$ a place of $E$ above $\ell$, and $\overline{E_\lambda}$ an algebraic closure of $E_\lambda$. Assume furthermore that $\pi$ is {\it selfdual}, {\it i.e.} isomorphic to its contragredient $\pi^\vee$. Then by \cite{harristaylor,shin,chharris} we know that there exists a continuous semisimple representation 
\begin{equation} \label{galoisrepeq} r_{\pi,\lambda} : {\rm Gal}(\overline{F}/F) \longrightarrow {\rm GL}_n(\overline{E_\lambda}),\end{equation}
unique up to isomorphism, satisfying the following property: for each finite place $v$ of $F$ which is prime to $\ell$, and such that $\pi_v$ is unramified, the Galois representation of $r_{\pi,\lambda}$ is unramified at $v$ (see e.g. \cite[\S 2.1]{serreabelian}) and the characteristic polynomial of a Frobenius element ${\rm Frob}_v$ at $v$ satisfies $\det (t -r_{\pi,\lambda}({\rm Frob}_v)) =  \det ( t - {\rm c}(\pi_v))$. 
\ps\ps
\subsection{} \label{thmsgalois} We shall denote by ${\rm G}_2$ the automorphism group scheme of the standard split octonion algebra over $\Z$. This is a linear algebraic group, which is split and semisimple over $\Z$, and {\it of type ${\bf G}_2$}. Recall that for any algebraically closed field $k$ of characteristic $0$, there is up to isomorphism a unique $7$-dimensional irreducible $k$-linear algebraic representation of ${\rm G}_2(k)$ (see Lemma \ref{lemmerepg2}). We shall denote by $\rho : {\rm G}_2(k) \rightarrow {\rm GL}_7(k)$ such a representation. 
\ps\ps

\begin{thm}\label{thmgalois} Let $F$ be a totally real number field, $\pi$ a cuspidal automorphic representation of ${\rm GL}_7(\AAA_F)$ such that $\pi_v$ is algebraic regular for each real place $v$ of $F$, and $E$ a coefficient field for $\pi$. Assume that $\pi$ satisfies the following condition: \ps 

{\rm (G)} {\it For all but finitely many places $v$ of $F$ such that $\pi_v$ is unramified, then ${\rm c}(\pi_v)$ is the conjugacy class of an element in $\rho({\rm G}_2(\C))$.} \ps\ps

Then for any prime $\ell$ and any place $\lambda$ of $E$ above $\ell$, there is a morphism $$\widetilde{r}_{\pi,\lambda}: {\rm Gal}(\overline{F}/F) \longrightarrow {\rm G}_2(\overline{E_\lambda}),$$ 
unique up to ${\rm G}_2(\overline{E_\lambda})$-conjugacy, such that the representation $\rho \circ \widetilde{r}_{\pi,\lambda}$ is isomorphic to the representation $r_{\pi,\lambda}$ introduced in \eqref{galoisrepeq} above.
\end{thm}
\ps\ps

By the Cebotarev theorem, an equivalent way to state the theorem, but without any reference to $r_{\pi,\lambda}$, is as follows.\ps\ps

\begin{cor} \label{cor1thmgalois} Let $F$, $\pi$, $E$, $\ell$ and $\lambda$ be as in the statement of Theorem \ref{thmgalois}. Then there exists a continuous semisimple morphism $$\widetilde{r}_{\pi,\lambda}: {\rm Gal}(\overline{F}/F) \longrightarrow {\rm G}_2(\overline{E_\lambda}),$$ unique up to ${\rm G}_2(\overline{E_\lambda})$-conjugacy, satisfying the following property: for each finite place $v$ of $F$ which is prime to $\ell$, and such that $\pi_v$ is unramified, the morphism $\widetilde{r}_{\pi,\lambda}$ is unramified at $v$ and we have $\det (t - \rho(\widetilde{r}_{\pi,\lambda}({\rm Frob}_v))) = \det (t - {\rm c}(\pi_v))$. 
\end{cor}
\ps
In this statement, the term {\it semisimple} means that the neutral component of the Zariski-closure of the image of $\widetilde{r}_{\pi,\lambda}$ is reductive. It is equivalent to ask that for some (resp. any) faithful polynomial representation $f$ of ${\rm G}_2$ over $\overline{E_\lambda}$ the representation $f \circ  \widetilde{r}_{\pi,\lambda}$ is semisimple. 
\ps\ps

\begin{pf} (of Theorem \ref{thmgalois})  Assumption (G) on $\pi$ implies that for all but finitely many places $v$ of $F$ such that $\pi_v$ is unramified, the conjugacy class ${\rm c}(\pi_v) \subset {\rm GL}_7(\C)$ is equal to its inverse ${\rm c}(\pi_v^\vee)$. As $\pi^\vee$ is a cuspidal automorphic representation of $\GL_7(\AAA_F)$, the strong multiplicity one theorem of Piatetski-Shapiro, Jacquet and Shalika shows that $\pi$ is selfdual, so that the discussion of \S\ref{algreg} applies to $\pi$. \ps\ps
Consider a Galois representation $r_{\pi,\lambda}$ as in \eqref{galoisrepeq}. By the main theorem of \cite{bc}, there is a structure of non-degenerate quadratic space on the underlying space $V = \overline{E_\lambda}^7$ of $r_{\pi,\lambda}$ such that the image $\Gamma$ of $r_{\pi,\lambda}$ is a subgroup of ${\rm SO}(V)$. To conclude the existence of a morphism $\widetilde{r}_{\pi,\lambda}$ such that $\rho \circ \widetilde{r}_{\pi,\lambda}$ is isomorphic to $r_{\pi,\lambda}$ it is enough to show that $\Gamma$ is contained in a ${\rm G}_2$-subgroup of ${\rm SO}(V)$ (\S \ref{defg2sb}).\ps\ps

The semisimplicity of $r_{\pi,\lambda}$ implies that $\Gamma$ satisfies assumption (i) of Theorem \ref{mainthmb} (i.e. assumption {\rm (S)}). The assumption (G) on $\pi$, Cebotarev's theorem and Remark \ref{rempotg2} ensure that $\Gamma$ satisfies as well assumption (ii) of this theorem. So Theorem \ref{mainthmb} applies, and it is enough to show that $\Gamma$ is not as in case (a), (b), (c) or (d) of its statement. By the concrete description of these specific cases, observe that it is enough to show that there is no open subgroup $\Gamma'$ of the profinite group $\Gamma$ such that the subspace of invariants $V^{\Gamma'} \subset V$ has dimension $\geq 3$ over $\overline{E_\lambda}$.\ps\ps

Let $w$ be a finite place of $F$ dividing $\ell$. The restriction to a decomposition group at $w$ of the representation $r_{\pi,\lambda}$, that we may view as a representation $r_w: {\rm Gal}(\overline{F_w}/F_w) \rightarrow \GL_7(\overline{E_\lambda})$, is known to be of Hodge-Tate type \cite{shin,chharris}. More precisely, the Hodge-Tate weights of $r_w$, relative to any field embedding $F_w \longrightarrow \overline{E_\lambda}$, are the weights of the algebraic regular representation $\pi_v$, for a suitable real place $v$ of $F$ (depending on the choice of that field embedding). In particular, those weights are distinct integers and only one of them is zero (note that if $k$ is weight then so is $-k$ by selfduality of $r_{\pi,\lambda}$). As a general fact, the same properties trivially hold for the restriction of $r_w$ to ${\rm Gal}(\overline{F_w}/L)$ for any finite extension $F_w \subset L \subset \overline{F_w}$. In particular, if $L$ is any such extension, the multiplicity of the Hodge-Tate weight $0$ for ${r_w}_{|{\rm Gal}(\overline{F_w}/L)}$, relative to any field embedding $L \rightarrow
\overline{E_\lambda}$, is at most one. But this implies that the dimension of the space of invariants of $V$ under $r_w({\rm Gal}(\overline{F_w}/L))$ is at most one as well, and we are done. \ps\ps

It only remains to justify the uniqueness assertion. But this is a consequence of a result of Griess \cite[Thm. 1]{griess} (see also \cite[Prop. 2.8]{larsen}).  \end{pf}

\begin{remark} {\rm If $r_{\pi,\lambda}$ is irreducible, which is conjectured to always be the case but is still an open problem at the moment, we can use in the proof above Theorem \ref{thmbirrcase} instead of Theorem \ref{mainthmb}, whose proof is much easier. The known compatibilty of $r_{\pi,\lambda}$ with the local Langlands correspondence shows that this irreducibility assumption is automatically satisfied if, for some finite place $v$, the representation $\pi_v$ is square integrable (modulo the center). In the special case where $\pi_v$ is the Steinberg representation for some $v$, Theorem \ref{thmgalois} even has a shorter proof, as already observed by Magaard and Savin in \cite[\S 7]{ms}. See also \cite{kls} for some quite specific special cases of Theorem \ref{thmgalois}. }
\end{remark}
In the irreducible case, we can draw slightly further conclusions about the field of definition of $\widetilde{r}_{\pi,\lambda}$.
\ps\ps

\begin{cor} \label{cor2thmgalois} Let $F$, $\pi$, $E$, $\ell$ and $\lambda$ be as in the statement of Theorem \ref{thmgalois}. Assume furthermore that $r_{\pi,\lambda}$ is irreducible. Then the conclusion of Corollary \ref{cor1thmgalois} holds with ${\rm G}_2(\overline{E_\lambda})$ replaced {\rm (}twice{\rm )} by ${\rm G}_2(E_\lambda)$.
\end{cor}

 {\scriptsize\begin{pf} Choose $\widetilde{r}_{\pi,\lambda}$ as in Theorem \ref{thmgalois} and denote by $C$ the centralizer of the image of $\widetilde{r}_{\pi,\lambda}$ in ${\rm G}_2(\overline{E_\lambda})$. The irreducibility of $r_{\pi,\lambda}$ and Schur's lemma imply that $C$ is central in ${\rm G}_2(\overline{E_\lambda})$, hence the equality $C=\{1\}$. Set $k=E_\lambda$, $\overline{k}=\overline{E_\lambda}$, $\Gamma= {\rm Gal}(\overline{E_\lambda}/E_\lambda)$ and $G={\rm G}_2$. Then $\Gamma$ naturally acts on $G(\overline{k})$, with fixed points $G(k)$. Let $\sigma \in \Gamma$. Observe that we have $$\det(t - \rho \circ \sigma \circ \widetilde{r}_{\pi,\lambda}(g))=\sigma(\det(t-r_{\pi,\lambda}(g)))=\det(t-r_{\pi,\lambda}(g))=\det(t - \rho \circ \widetilde{r}_{\pi,\lambda}(g))$$ for all $g \in {\rm Gal}(\overline{F}/F)$. The uniqueness part of the statement of Theorem \ref{thmgalois} implies the existence of an element $P_\sigma \in G(\overline{k})$ such that we have $\sigma \circ \widetilde{r}_{\pi,\lambda} = {\rm int}_{P^{-1}_\sigma} \circ \widetilde{r}_{\pi,\lambda}$ (where ${\rm int}_g$ denotes the inner automorphism $x \mapsto gxg^{-1}$ of $G(\overline{k})$). As $C=\{1\}$ the element $P_\sigma$ is uniquely determined by this relation. This implies that the map $\sigma \mapsto P_\sigma$ is a $1$-cocycle of $\Gamma$ acting on $G(\overline{k})$; it is continuous for the discrete topology on $G(\overline{k})$ as we have ${\rm Im}\, \widetilde{r}_{\pi,\lambda} \subset G(K)$ for some finite extension $k \subset K \subset \overline{k}$ by Baire's theorem. The corresponding cohomology group classifies the isomorphism classes of octonion $k$-algebras. But since $k$ is a nonarchimedean local field, there is a unique such $k$-algebra up to isomorphism, as any $8$-dimensional non-degenerate quadratic space does represent $0$ \cite[\S 2,\S 3]{blijspringer1}. This shows the existence of $P \in G(\overline{k})$ with $P_\sigma = P^{-1}\sigma(P)$ for all $\sigma \in \Gamma$, and then that ${\rm int}_{P} \circ \widetilde{r}_{\pi,\lambda}$ is fixed under $\Gamma$. This shows the existence assertion. The uniqueness follows from the uniqueness part of the statement of Theorem \ref{thmgalois} and from $C=\{1\}$. \end{pf}}
\ps\ps

Our next application concerns the local components of a cuspidal automorphic representation $\pi$ of $\GL_7(\AAA_F)$ satisfying assumption (G) of Theorem \ref{thmgalois}. We refer to \cite[\S 4]{tate} for the various notions of Weil groups, Weil-Deligne groups and representations. Let $v$ be a finite place of $F$ and ${\rm W}_{F_v} \subset {\rm Gal}(\overline{F_v}/F_v)$ the Weil group of $F_v$. Recall that the local Langlands correspondence, proved by Harris and Taylor \cite{harristaylor}, may be viewed as a natural bijection $\varpi \mapsto {\rm L}(\varpi)$ from the set of isomorphism classes of irreducible smooth complex representations of ${\rm GL}_n(F_v)$ and onto the set of isomorphism classes of continuous semisimple representations ${\rm W}_{F_v} \times {\rm SU}(2) \rightarrow \GL_n(\C)$. \ps\ps

\begin{cor} \label{corcplx} Let $F$ be a totally real number field and $\pi$ a cuspidal automorphic representation of ${\rm GL}_7(\AAA_F)$ such that $\pi_v$ is algebraic regular for each real place $v$ of $F$. The following properties are equivalent: \begin{itemize}\ps
\item[(i)] $\pi$ satisfies property {\rm (G)} of the statement of Theorem \ref{thmgalois},\ps\ps
\item[(ii)] for any finite place $v$ of $F$, there exists a continuous semisimple morphism $\phi_v : {\rm W}_{F_v} \times {\rm SU}(2) \rightarrow {\rm G}_2(\C)$, unique up to ${\rm G}_2(\C)$-conjugacy, such that $\rho \circ \phi_v$ is isomorphic to ${\rm L}(\pi_v)$.\ps
\end{itemize}
\end{cor}

\begin{pf} The nontrivial assertion is (i) $\Rightarrow$ (ii). Assume (i) holds. Let $v$ be a finite place of $F$, $E$ a coefficient field of $\pi$, $\ell$ a prime such that $v$ is not above $\ell$, $\lambda$ a place of $E$ above $\ell$, set $k=\overline{E_\lambda}$. Let $\widetilde{r}_{\pi,\lambda}$ be given by Theorem \ref{thmgalois}. \ps\ps

Recall that if $G$ is a connected reductive group over $k$, there is a natural map from the set of $G(k)$-conjugacy classes of continuous morphisms $r : {\rm W}_{F_v} \rightarrow G(k)$, where $G(k)$ is equipped with the $\ell$-adic topology, onto the set of $G(k)$-conjugacy classes of semisimple morphisms $\varphi : {\rm W}_{F_v} \times {\rm SL}_2(k) \rightarrow G(k)$ which are $k$-algebraic on the second factor and trivial on $H \times 1$ for some open subgroup $H$ of (the inertia subgroup) of ${\rm W}_{F_v}$. The morphism $\varphi$ is a version of the so-called {\it Frobenius-semisimple Weil-Deligne representation} associated to $r$; its construction, which is standard at least in the $G=\GL_n$ case, is an application of Grothendieck's $\ell$-adic monodromy theorem and of the Jacobson-Morosov theorem (see \cite[\S 4]{tate} and \cite[Prop. 2.2]{grossgamma}). \ps\ps

The restriction of $\widetilde{r}_{\pi,\lambda}$ to a decomposition group at $v$ gives rise to a continuous morphism $\widetilde{r}_v: {\rm W}_{F_v} \rightarrow {\rm G}_2(k)$. If $\widetilde{\varphi}$ is an associated Frobenius-semisimple representation, then $\varphi = \rho \circ \widetilde{\varphi}$ is a Frobenius-semisimple representation associated to $\rho \circ r_v$. It is convenient to choose a field embedding $\sigma : k \longrightarrow \C$ whose restriction to $E$, which is naturally embedded in $E_\lambda$, is the given inclusion $E \subset \C$. The known compatibility of $r_{\pi,\lambda}$ with the local Langlands correspondence \cite{shin,chharris} asserts that $\varphi \otimes_k \C$ (where the scalar extension is made through $\sigma$) is the complexification of ${\rm L}(\pi_v)$. We can thus take for $\phi_v$ the morphism $\widetilde{\varphi} \otimes_k \C$, restricted to ${\rm W}_{F_v} \times {\rm SU}(2)$. The uniqueness assertion is again a consequence of \cite[Thm. 1]{griess} . \end{pf}

\begin{remark} \label{remconjal1}  {\rm A similar application of Theorem \ref{thmbirrcase}, but arguing with the hypothetical Langlands group instead of a Galois group as on page \pageref{lanintro} of the introduction, shows that properties (i) and (ii) of corollary \ref{corcplx} should be equivalent without assuming $F$ is totally real or anything about the Archimedean places of $\pi$. Moreover, those assumptions should be equivalent to the existence of a cuspidal tempered automorphic representation $\pi'$ of ${\rm G}_2(\AAA_F)$ such that for all but finitely many finite places $v$ of $F$, the conjugacy class of $\rho({\rm c}(\pi'_v))$ is ${\rm c}(\pi_v)$ (recall that we have ${\rm c}(\pi'_v) \subset \widehat{{\rm G}_2}={\rm G}_2(\C)$).  }
\end{remark}

\subsection{Examples}\label{exgalois} Examples of automorphic representations $\pi$ satisfying the assumptions of Theorem \ref{thmgalois} can be constructed using Arthur's recent results on automorphic representations of classical groups \cite{arthurbook}. \ps\ps

One method, developed by Gross and Savin \cite{gsduke}, is to use an exceptional theta correspondence between anisotropic forms of ${\rm G}_2$ and ${\rm PGSp}_6$. This is also the point of view of Magaard and Savin in \cite{ms}, who give some interesting examples. Another method, that we hope to pursue in the future, is to use triality for the algebraic group ${\rm PGSO}_8$ over $F$ (or for suitable forms of it). \ps\ps

Let us mention that in \cite[\S 8]{chrenard} the authors give an explicit conjectural formula for the number of cuspidal automorphic representations $\pi$ of $\GL_7(\AAA_\Q)$ satisfying assumption {\rm (G)} and such that:
\begin{itemize}\ps
\item[(i)] $\pi_p$ is unramified for each prime $p$, \ps\ps
\item[(ii)] $\pi_\infty$ is algebraic of weights $0, \pm u, \pm v, \pm (u+v)$ with $0<u<v$;\ps\ps
\end{itemize}
this is the number denoted ${\rm G}_2(2v,2u)$ {\it loc. cit}. For instance, by \cite[Table 11]{chrenard}, there should be $3$ such representations with $u+v \leq 13$ (each one having $\Q$ as coefficient field), namely with $$(u,v)=(4,8), \,\,\,\,\,(3,10)\,\,\,\,\,{\rm and}\,\,\,\,\, (5,8).$$
Some informations about the Satake parameters of these representations $\pi$, such as ${\rm c}(\pi_2)$ and ${\rm trace}\, {\rm c}(\pi_p)$ for $p \leq 13$, have also been recently obtained by Megarbane  \cite[Tables 8 \& 9]{megarbane}. In particular, it should be possible to study the images of their associated Galois representations $\widetilde{r}_{\pi,\lambda} : {\rm Gal}(\overline{\Q}/\Q) \rightarrow {\rm G}_2(\overline{\Q_\ell})$ given by Theorem \ref{thmgalois}.\ps\ps 

\begin{remark} {\rm Let $\pi$ be a cuspidal automorphic representation of $\GL_7(\AAA_\Q)$ satisfying assumption {\rm (G)}, as well as (i) and (ii) above. The crystalline property of the restriction of $r_{\pi,\lambda}$ at a decomposition group at a prime $\ell$ with $\lambda\, |\, \ell$, as well as Minkowski's theorem that any nontrivial number field has a ramified prime, imply that the Zariski-closure of the image of $r_{\pi,\lambda}$ is connected. In particular, in order to show the existence of $\widetilde{r}_{\pi,\lambda}$ in this case one may invoke (the much more elementary) Corollary \ref{corconnected} instead of Theorem \ref{mainthmb}.}
\end{remark}

\subsection{A local characterization of Langlands' functorial lifting between  ${\rm G}_2$ and ${\rm PGSp}_6$}\label{localcarg2sp6}

In this last paragraph, we fix an injective algebraic group morphism $$\rho : \widehat{{\rm G}_2}={\rm G}_2(\C) \longrightarrow \widehat{{\rm PGSp}_6} = {\rm Spin}_7(\C).$$ Let $\pi$ be a cuspidal tempered automorphic representation of ${\rm PGSp}_6(\AAA_F)$. As already explained on page \pageref{lanintro} of the introduction, if we assume the existence of the Langlands group of the number field $F$ then Theorem \ref{thmintroa} shows the equivalence between the following properties: 
\begin{itemize}\ps \ps
\item[(P1)] for all but finitely many places $v$ of $F$, the Satake parameter ${\rm c}(\pi_v)$ is the conjugacy class of an element of a ${\rm G}_2$-subgroup of ${\rm Spin}_7(\C)$,\ps\ps 
\item[(P2)] there exists a cuspidal automorphic representation $\pi'$ of ${\rm G}_2(\AAA_F)$ such that for all but finitely many places $v$ of $F$, the conjugacy class of $\rho({\rm c}(\pi'_v))$ is ${\rm c}(\pi_v)$.\ps\ps\end{itemize}

As $\pi$ was assumed to be tempered, note that if $\pi'$ is as in property (P2) then $\pi'$ is tempered as well (see \S \ref{notationsaut} for our convention on the meaning of the term {\it tempered} in this paper).\ps\ps 

It turns out that the equivalence between these properties can be proved unconditionally. The following theorem grew out from a discussion with Wee Teck Gan and Gordan Savin, after the first version of this paper appeared on arXiv. \ps\ps

\begin{thm}\label{thmCGS} Let $\pi$ be a cuspidal automorphic representation of ${\rm PGSp}_6(\AAA_F)$. Assume that $\pi$ is tempered, or more generally that $\pi$ is {\rm nearly generic}. Then $\pi$ satisfies {\rm (P1)} if, and only if, it satisfies {\rm (P2)}. \end{thm}

Let us first define the term {\it nearly generic} which occurs in this statement. If $\pi$ is an admissible irreducible representation of ${\rm GSp}_{2g}(\AAA_F)$ we shall denote by ${\rm res}\, \pi$ its restriction to ${\rm Sp}_{2g}(\AAA_F)$. It is an admissible semisimple representations of ${\rm Sp}_{2g}(\AAA_F)$ : this follows from \cite[Lemma 2.3.1]{Xu} and the equality $\dim \pi_v^{{\rm Sp}_{2g}(\mathcal{O}_v)} =1$ when we have $\dim  \pi_v^{{\rm GSp}_{2g}(\mathcal{O}_v)} =1$. Denote by $$\eta : \widehat{{\rm GSp}_{2g}} = {\rm GSpin}_{2g+1}(\C) \longrightarrow {\rm SO}_{2g+1}(\C) = \widehat{{\rm Sp}_{2g}}$$ the natural morphism. Let $\sigma$ be an irreducible constituent of ${\rm res}\, \pi$. If $\pi$ and $\sigma$ are unramified at the finite place $v$, which holds for almost all $v$, then the compatibility of the Satake isomorphism to isogenies (proved by Satake) shows the equality
\begin{equation}\label{satiso} \eta({\rm c}(\pi_v))={\rm c}(\sigma_v). \end{equation}
(see also \cite[Lemma 5.2]{Xupack} for another proof). If $\pi$ is cuspidal, then the restriction of functions along the inclusion ${\rm Sp}_{2g}(\AAA_F) \rightarrow {\rm GSp}_{2g}(\AAA_F)$ shows that a nonzero quotient of ${\rm res}\, \pi$ is a direct sum of cuspidal automorphic representations, by Proposition \ref{extcuspform} (ii). By the equality \eqref{satiso}, all the cuspidal constituents of ${\rm res}\, \pi$ belong to a same global packet of the form $\widetilde{\Pi}_{\psi}$ in the sense of Arthur \cite[p. 45]{arthurbook}. We shall say that $\pi$ is {\it nearly generic} if this packet (or $\psi$) is generic in the sense of {\it loc. cit.} p. 49. This condition is automatically satisfied if $\pi$ is tempered by the Jacquet-Shalika estimates. Note also that if $\pi$ is tempered then so is any constituent of ${\rm res}\, \pi$, by the equality \eqref{satiso} and our definition of the term tempered in \S \ref{notationsaut}. The terminology {\it nearly generic} will be justified by Proposition \ref{propgengen}.

\begin{pf} (of Theorem \ref{thmCGS}) The implication (P2) $\Rightarrow$ (P1) is trivial. Assume now that (P1) holds. Denote by $({\rm spin},W)$ a spinor representation of ${\rm Spin}_7(\C)$, and $({\rm st},E)$ its standard representation. Let $S$ be a finite set of places of $F$ containing ${\rm Ram}(\pi)$, and such that for all $v \notin S$ then $\pi_v$ is unramified and ${\rm c}(\pi_v)$ is the conjugacy class of an element of a ${\rm G}_2$-subgroup of ${\rm Spin}_7(\C)$. Recall that we have a $\C[H]$-linear isomorphism $W \simeq 1 \oplus E$ over any ${\rm G}_2$-subgroup $H$ of ${\rm Spin}_7(\C)$. It follows that for ${\rm Re} \,s$ big enough we have a decomposition $${\rm L}^S(s,\pi,{\rm spin}) = \zeta_F^S(s) {\rm L}^S(s,\pi,{\rm st})$$
 where $\zeta_F^S(s)={\rm L}^S(s,\pi,1)$ is the partial zeta function of $F$. \ps\ps
 
 Let $\sigma$ be a cuspidal automorphic representation of ${\rm Sp}_6(\AAA_F)$ which is a constituent of ${\rm res}\, \pi$. Up to enlarging $S$ if necessary we may assume that $\sigma_v$ is unramified for $v \notin S$. Formula \eqref{satiso} shows ${\rm L}^S(s,\pi,{\rm st})={\rm L}^S(s,\sigma,{\rm st})$. As $\pi$ is nearly generic, the main global result of Arthur \cite[Thm. 1.5.2]{arthurbook} asserts that this ${\rm L}$-function is a finite product of Godement-Jacquet ${\rm L}$-functions of certain cuspidal unitary automorphic representations of the general linear groups over $F$. Such an ${\rm L}$-function has a meromorphic continuation to $\C$, and does not vanish at $s=1$ by a classical result of Jacquet and Shalika. It follows that ${\rm L}^S(s,\pi,{\rm spin})$ has a meromorphic continuation to $\C$ with a pole at $s=1$, coming from $\zeta_F^S(s)$. \ps\ps
If $\pi$ is {\it globally generic} (see below), this last assertion and the main result of Ginzburg and Jiang \cite[Thm. 1.1]{gjiang} show that (P2) holds; the sought after automorphic representation of ${\rm G}_2$ is obtained there by an exceptional theta correspondence. In general, we reduce to this case by Proposition \ref{propgengen} below.
\end{pf}
\ps\ps
Recall, following e.g. \cite[\S 1]{CKPS}, that a cuspidal automorphic representation $\pi$ of a split reductive group $H$ is said {\it globally generic} if there exists an irreducible $H(\AAA_F)$-submodule of the space of cuspforms on $H(\AAA_F)$ which is isomorphic to $\pi$, and which is generic in the traditional sense (recalled in the appendix, before the statement of Proposition \ref{extcuspformgen}).\ps\ps

\begin{prop}\label{propgengen} Let $g\geq 1$ be an integer and $\pi$ a cuspidal automorphic representation of ${\rm GSp}_{2g}(\AAA_F)$. Then $\pi$ is nearly generic if, and only if, there exists a globally generic cuspidal automorphic representation $\varpi$ of ${\rm GSp}_{2g}(\AAA_F)$ with the same central character as $\pi$, and which is {\it nearly equivalent} to $\pi$, {\rm i.e.} such that $\varpi_v \simeq \pi_v$ for all but finitely many places $v$ of $F$.
\end{prop}

\begin{pf} Set $G = {\rm Sp}_{2g}$ and $\widetilde{G}={\rm GSp}_{2g}$. Assume first that $\pi$ is nearly generic. We denote by $\widetilde{c}$ the central character of $\pi$ (a Hecke character of $F$). \ps\ps

Let $\sigma$ be a cuspidal automorphic constituent of ${\rm res} \,\,\pi$. Arthur's theorem \cite[Thm. 1.5.2]{arthurbook} and the work of Cogdell, Kim, Piatetski-Shapiro and Shahidi \cite[Thm. 7.2]{CKPS} show that there is a globally generic cuspidal automorphic representation $\sigma'$ of $G(\AAA_F)$ which is nearly equivalent to $\sigma$ (see also \cite[Prop. 8.3.2]{arthurbook}). We shall denote by $V$ a generic irreducible subspace of the space of cuspforms on $G(\AAA_F)$ giving rise to $\sigma'$. \ps\ps

According to Arthur, $\sigma$ and $\sigma'$ are in a same (generic) global packet; in particular, they have the same central character (see the paragraph before \cite[Prop. 6.27]{Xupack}), namely the restriction $c$ of $\widetilde{c}$ to the center of $G(\AAA_F)$. Denote by $\rho$ the natural restriction map $\mathcal{A}_{\rm cusp}(\widetilde{G})_{\widetilde{c}} \longrightarrow \mathcal{A}_{\rm cusp}(G)_c$ studied in Proposition \ref{extcuspform}. By the assertion (iii) of this proposition, this map is surjective, so there is an irreducible admissible $\widetilde{G}(\AAA_F)$-subrepresentation $\widetilde{V} \subset \mathcal{A}_{\rm cusp}(\widetilde{G})_{\widetilde{c}}$ whose image in $\mathcal{A}_{\rm cusp}(G)_c$ contains the space $V$ defined above. We denote by $\pi'$ the isomorphism class of the cuspidal automorphic $\widetilde{G}(\AAA_F)$-representation $\widetilde{V}$. By Proposition \ref{extcuspformgen}, $\widetilde{V}$ is generic as so is $V$, so $\pi'$ is globally generic. \ps\ps

By construction, ${\rm res}\, \pi$ has a cuspidal automorphic constituent (namely $\sigma$) which is nearly equivalent to some cuspidal automorphic constituent of ${\rm res}\, \pi'$ (namely $\sigma'$). Let $\nu : \widetilde{G}(\AAA_F) \rightarrow  \AAA_F^\times$ denote the similitude factor. By Xu's theorem \cite[Thm. 1.8]{Xupack}, there exists a Hecke character $\omega$ of $F$ such that $\pi$ and $\pi' \otimes \omega \circ \nu$ belong to the same ``global $L$-packet'' that he defines. In particular, they are nearly equivalent and have the same central character (see \cite[Thm. 1.2]{Xupack}). The representation $\varpi := \pi' \otimes \omega \circ \nu$ of $\widetilde{G}(\AAA_F)$ is cuspidal, globally generic, nearly equivalent to $\pi$, and with central character $c$, and we are done.  \ps\ps

Conversely, up to replacing $\pi$ by $\varpi$ if necessary we may assume that $\pi$ is globally generic. In order to show that $\pi$ is nearly generic, it is enough to prove that ${\rm res}\,\pi$ has a cuspidal globally generic irreducible constituent, by \cite[Thm. 7.2]{CKPS}. But this follows from Proposition \ref{extcuspformgen}.
\end{pf}

\appendix

\section{Some basic facts about restrictions and extensions of cusp forms}\label{appresext}
\ps\ps\ps
{\scriptsize 
The results of this appendix are probably well known to some specialists, however we have not been able to find their proofs in the litterature. Some of the arguments given below already appear in some form in Labesse-Langlands' paper \cite[p. 49-52]{LL} and in \cite[Lemma 5.3]{Xulift}. Our main reference for what follows is Borel-Jacquet's paper \cite{boreljacquet}. \ps\ps

Let $F$ be a number field and $H$ a reductive linear algebraic group defined over $F$. We set $H_\infty=H(F \otimes_\Q \R)$, the direct product of the Lie groups $H(F_v)$ for $v$ Archimedean, and denote by $H_f$ the locally profinite group $H(\AAA_F^f)$ where $\AAA_F^f$ is the $F$-algebra of finite ad\`eles of $F$: we have $H(\AAA_F) = H_\infty \times H_f$. If $c$ is a character of the center of $H(\AAA_F)$, we denote by $\mathcal{A}_{\rm cusp}(H)_c$ the space of cuspidal automorphic functions $H(\AAA_F) \rightarrow \C$ with central character $c$ and with respect to a choice of maximal compact subgroup $K_\infty$ of $H_\infty$ (in the sense of \cite[\S 4.4]{boreljacquet}). This is an admissible, semisimple, $({\rm Lie}\, H_\infty,K_\infty) \times H_f$-module. For short, a $({\rm Lie}\, H_\infty,K_\infty) \times H_f$-module will be called an {\it $H(\AAA_F)$-{\rm gk}-module}. \ps\ps

\begin{prop}\label{extcuspform} Let $\widetilde{G}$ be a reductive linear algebraic group, $T$ a torus, $\nu : \widetilde{G} \rightarrow T$ a surjective morphism, $\widetilde{Z}$ the center of $\widetilde{G}$, and $G$ the kernel of $\nu$. We assume that $\widetilde{G}$, $T$ and $\nu$ {\rm (}hence $\widetilde{Z}$ and $G${\rm )} are defined over the number field $F$, and that $G$ is connected. Let $\widetilde{c}$ be an automorphic character of $\widetilde{Z}(\AAA_F)$ and $c$ its restriction to the center of $G(\AAA_F)$. Then the following assertions hold: \begin{itemize}\ps\ps
\item[(i)] the restriction map $\rho : \varphi \mapsto \varphi_{|G(\AAA_F)}$ induces an equivariant morphism $\rho_{\rm cusp} : {\rm res}\, \mathcal{A}_{\rm cusp}(\widetilde{G})_{\widetilde{c}} \longrightarrow \mathcal{A}_{\rm cusp}(G)_c$ of $G(\AAA_F)$-{\rm gk}-modules, \ps\ps
\item[(ii)] if $V$ is a nonzero $\widetilde{G}(\AAA_F)$-{\rm gk}-submodule of $\mathcal{A}_{\rm cusp}(\widetilde{G})_{\widetilde{c}}$ then we have $\rho_{\rm cusp}(V) \neq 0$,\ps\ps
\item[(iii)] if furthermore $\widetilde{Z}$ is a split $F$-torus, then $\rho_{\rm cusp}$ is surjective. \ps\ps
\end{itemize}
\end{prop}

\begin{pf} As $\nu$ is trivial on the derived subgroup of $\widetilde{G}$, we have $\nu(\widetilde{Z})=T$. We may thus find a torus $C \subset \widetilde{Z}$ defined over $F$ such that $\nu : C \rightarrow T$ is an isogeny. We shall denote by $$\iota : C \times G \rightarrow \widetilde{G}$$ 
the isogeny $(z,g) \mapsto zg$. For each place $v$ of $F$, note that $G(F_v)$ is a normal subgroup of $\widetilde{G}(F_v)$ with abelian quotient. Moreover, $\nu(C(F_v))$ is an open subgroup of finite index in $T(F_v)$, by the finiteness of the Galois cohomology group $H^1(F_v, C \cap G)$ (Tate). In particular, the image of the morphism $\iota_v : C(F_v) \times G(F_v) \rightarrow \widetilde{G}(F_v)$ induced by $\iota$ is an open subgroup of finite index in $\widetilde{G}(F_v)$. For each Archimedean place $v$, we choose a maximal compact subgroup $K_v$ of $G(F_v)$, as well as a maximal compact subgroup $\widetilde{K}_v$ of $\widetilde{G}(F_v)$ containing $K_v$. We have thus $K_v = G(F_v) \cap \widetilde{K}_v$. We define $K_\infty$ (resp. $\widetilde{K}_\infty$) as the product of the $K_v$ (resp. $\widetilde{K}_v$) for $v$ Archimedean. We may and do choose these maximal compact subgroups in order to define  $\mathcal{A}_{\rm cusp}(\widetilde{G})_{\widetilde{c}}$ and $\mathcal{A}_{\rm cusp}(G)_c$. If we denote by $\mathfrak{z}(H)$ the $\C$-algebra defined as the center of the complex enveloping algebra ${\rm U}({\rm Lie}\, H)$ of the Lie group $H$, then for each Archimedean $v$ the local isomorphism $\iota_v$ defines a $\C$-algebra isomorphism  $\mathfrak{z}(C(F_v)) \otimes \mathfrak{z}(G(F_v)) \isomo \mathfrak{z}(\widetilde{G}(F_v)) $.  In particular, the natural inclusion ${\rm U}({\rm Lie}\, G(F_v)) \subset {\rm U}({\rm Lie}\, \widetilde{G}(F_v))$ induces an injection $\mathfrak{z}(G(F_v)) \subset \mathfrak{z}(\widetilde{G}(F_v))$. \ps\ps


Let us check assertion (i). Let $\varphi \in \mathcal{A}_{\rm cusp}(\widetilde{G})_{\widetilde{c}}$. There is a compact open subgroup $\widetilde{K}_f \subset \widetilde{G}_f$ such that $\varphi$ is $\widetilde{K}_f$-invariant on the right. Moreover, for any $g_f \in  \widetilde{G}_f$ the function $g_\infty \mapsto \varphi(g_\infty \times g_f), \widetilde{G}_\infty \rightarrow \C$, is smooth, $\widetilde{K}_\infty$-finite on the right, $\mathfrak{z}(\widetilde{G}_\infty)$-finite, and slowly increasing (see \cite[\S 1.2]{boreljacquet}). Using a model over an open subset of ${\rm Spec}\,\mathcal{O}_F$ of the closed morphism $G \subset \widetilde{G}$,  we see that $K_f= \widetilde{K}_f \cap G_f$ is a compact open subgroup of $G_f$. By the inclusions $K_f \subset \widetilde{K}_f$, $K_\infty \subset \widetilde{K}_\infty$, and $\mathfrak{z}(G_\infty) \subset \mathfrak{z}(\widetilde{G}_\infty)$ (justified above), it follows that $\rho(\varphi)$ is $K_f$-invariant and $K_\infty$-finite on the right, as well as $\mathfrak{z}(G_\infty)$-finite. It is also trivially slowly increasing (note {\it e.g.} that the restriction to $G_\infty$ of a {\it norm} on $\widetilde{G}_\infty$ is still a norm, in the sense of  \cite[\S 1.2]{boreljacquet}). Of course, $\rho(\varphi)$ is $G(F)$-invariant on the left and has central character $c$. We have proved that $\rho(\varphi)$ is an automorphic form.
It only remains to check that it is cuspidal. But for any parabolic subgroup $P$ of $G$ defined over $F$, there is a unique parabolic subgroup $\widetilde{P}$ of $\widetilde{G}$ such that $\widetilde{P} \cap G = P$, namely $\widetilde{P}=\widetilde{Z} \, P$. We conclude as the unipotent radical of $\widetilde{P}$ and $P$ coincide (in particular, they are included in $G$), and as $f$ is cuspidal. This shows that $\rho_{\rm cusp}$ is well defined. By construction, it is $({\rm Lie} \,G_\infty, K_\infty)\times G_f$-equivariant.\ps\ps

We now prove assertion (ii).  Let $V$ be as in the statement and $\varphi \in V - \{0\}$. Choose $g = g_\infty g_f \in \widetilde{G}(\AAA_F)$ such that $\varphi(g) \neq 0$. Up to replacing $\varphi$ by $g_f^{-1} \cdot \varphi \in V$ we may assume $g_f=1$. As the image $\iota_v$ is an open subgroup of $\widetilde{G}(F_v)$ for $v$ Archimedean, and as $\widetilde{K}_\infty$ meets every connected component of $\widetilde{G}_\infty$ by Cartan's decomposition, we have:
$$G_\infty \cdot \widetilde{Z}_\infty \cdot \widetilde{K}_\infty = \widetilde{G}_\infty.$$
It follows that up to replacing $\varphi$ by $h \cdot \varphi$ with $h \in \widetilde{Z}_\infty \cdot \widetilde{K}_\infty$, which preserves $V$, we may assume that we have $g \in G_\infty \subset G(\AAA_F)$. In particular, we have $\rho(\varphi) \neq 0$. \ps\ps

Let us now check (the main) assertion (iii). Denote by $Z$ the center of $G$; we have $Z = \widetilde{Z} \cap G$ as $\widetilde{G}=\widetilde{Z}\,G$. Consider the subgroups $G(\AAA_F) \subset H_1 \subset H_2 \subset \widetilde{G}(\AAA_F)$ with $H_1=\widetilde{Z}(\AAA_F) G(\AAA_F)$ and $H_2=\widetilde{G}(F) H_2$. They are normal in $\widetilde{G}(\AAA_F)$ (with abelian quotient), as so is $G(\AAA_F)$.  The assumption $\widetilde{c}_{|Z(\AAA_F)}=c$ and the obvious inclusion $\widetilde{Z}(\AAA_F)\cap G(\AAA_F) \subset Z(\AAA_F)$ show that there is a unique function $\varphi_1 : H_1 \rightarrow \C$ such that $\varphi_1(zg) = \widetilde{c}(z)\varphi(g)$ for all $(z,g) \in \widetilde{Z}(\AAA_F) \times G(\AAA_F)$. \ps\ps

We claim that there is a unique function $\varphi_2 : H_2 \rightarrow \C$ such that $\varphi_2(\gamma h)=\varphi_1(h)$ for all $(\gamma,h) \in \widetilde{G}(F) \times G(\AAA_F)$. Indeed, such a function exists if, and only if, $\varphi_1$ is left-invariant under $\widetilde{G}(F) \cap H_1$. It is thus enough to show the inclusion $\widetilde{G}(F) \cap H_1 \subset \widetilde{Z}(F) G(F)$, or which is the same, the inclusion $\nu(\widetilde{G}(F) \cap H_1) \subset \nu( \widetilde{Z}(F))$. It suffices to prove \begin{equation}\label{incltoprove} T(F) \cap \nu(\widetilde{Z}(\AAA_F)) \subset \nu ( \widetilde{Z}(F)).\end{equation} 

We have not used that $\widetilde{Z}$ is an $F$-split torus so far, but we shall do so now. As the morphism $\nu$ induces a surjective homomorphism $\widetilde{Z} \rightarrow T$ defined over $F$, it implies that $T$ is $F$-split. Moreover, all the subtori of $\widetilde{Z}$ are defined and split over $F$ as well. Let $Z^0$ be the neutral component of $Z=G \cap \widetilde{Z}$. Up to changing $C$ if necessary, we may thus assume that we have $\widetilde{Z} = C \times Z^0$. In particular, the following equalities hold
\begin{equation} \label{simplzc} \nu(\widetilde{Z}(\AAA_F)) = \nu(C(\AAA_F)), \, \, \, \widetilde{Z}(\AAA_F)G(\AAA_F) = C(\AAA_F) G(\AAA_F), \, \,\,  \widetilde{Z}(F)G(F) = C(F)G(F).\end{equation}
By the theory of elementary divisors, we may also assume that we have $T=\mathbb{G}_m^r$  for some integer $r\geq 0$, as well as an isomorphism $\mu :  \mathbb{G}_m^r \isomo C$ over $F$, and integers $m_1,\dots,m_r$, such that we have $\nu \circ \mu (z_1,\dots,z_r) = (z_1^{m_1},\dots,z_r^{m_r})$. But it is well-known that if $m \in \Z$ and $x \in F^\times$ is an $m$-th power in $F_v^\times$ for each $v$, then $x$ is an $m$-th power in $F^\times$.  This shows \eqref{incltoprove}, hence the existence of $\varphi_2$. \ps\ps

Let $S$ be a finite set of places of $F$ containing the Archimedean places. We denote by $\mathcal{O}_{F,S} \subset F$ the subring of $S$-integers. Up to enlarging $S$ if necessary, we may assume that $\widetilde{G}$ is an affine group schemes of finite type defined over $\mathcal{O}_{F,S}$, that $C \subset \widetilde{Z} \subset \widetilde{G}$ are closed subgroup schemes defined over $\mathcal{O}_{F,S}$, that $\nu : G \rightarrow \mathbb{G}_m$ is a group scheme homomorphism defined over $\mathcal{O}_{F,S}$, and that $\mu$ is a group scheme isomorphism $\mathbb{G}_m^r \isomo C$ over $\mathcal{O}_{F,S}$. We define again $G$ as the kernel of $\mu$; this is an affine group scheme of finite type over $\mathcal{O}_{F,S}$. By definition of the adelic topology, if we have a collection of open subsets $U_v \subset \widetilde{G}(F_v)$ for each $v \in S$, then $(\prod_{v \in S} U_v) \times (\prod_{v \notin S} \widetilde{G}(\mathcal{O}_v))$ is an open subset of $\widetilde{G}(\AAA_F)$ (and a similar assertion holds with $\widetilde{G}$ replaced by $G$). For any $v \notin S$ we shall set $\widetilde{K}_v = \widetilde{G}(\mathcal{O}_v)$, and we shall also set $\widetilde{K}^S = \prod_{v \notin S} \widetilde{K}_v$. Up to enlarging $S$ if necessary, we may assume that $\varphi$ is $G(\mathcal{O}_v)$-invariant on the right, as well as $\widetilde{c}_v(C(\mathcal{O}_v))=1$, for all $v \notin S$. We may also assume $m_i \in \mathcal{O}_{F,S}^\times$ for $i=1,\dots,r$. \ps\ps

Define $M$ as the maximum of the $m_i$ for $ 1\leq i \leq r$. By Hermite's theorem, there are only finitely many field extensions $F'/F$ of degree $\leq M$ and which are unramified everywhere. Up to enlarging $S$ if necessary, we may thus assume that the following extra property holds : if $F'/F$ is such an extension which is furthermore split at all places above $S$, then we have $F'= F$. In particular, for each $i=1,\dots, r$, if we have $x \in F^\times$ and $u \in \AAA_F^\times$ such that $u_v=1$ for $v \in S$, $u_v \in \mathcal{O}_v^\times$ for $v \notin S$, and such that $xu_v$ is an $m_i$-th power in $F_v^\times$ for each $v$, then $x$ is an $m_i$-th power in $F^\times$.  This property shows
$\widetilde{G}(F) \cap (C(\AAA_F) G(\AAA_F)  \widetilde{K}^S) \subset C(F) G(F)$, hence $\widetilde{G}(F) \cap (\widetilde{Z}(\AAA_F) G(\AAA_F)  \widetilde{K}^S) \subset \widetilde{Z}(F) G(F)$ by \eqref{simplzc}, which implies in turn the equality \begin{equation}\label{eq1phi3} ( \widetilde{G}(F)\widetilde{Z}(\AAA_F) G(\AAA_F) ) \cap \widetilde{K}^S = (\widetilde{Z}(\AAA_F) G(\AAA_F)) \cap  \widetilde{K}^S. \end{equation}
Observe moreover that for $v \notin S$ we have the inclusion \begin{equation}\label{eq2phi3}(\widetilde{Z}(F_v) G(F_v)) \cap \widetilde{G}(\mathcal{O}_v) \subset C(\mathcal{O}_v) G(\mathcal{O}_v).\end{equation} Indeed, we have $\widetilde{Z}(F_v) G(F_v) = C(F_v) G(F_v)$, and if $x \in \mathcal{O}_v^\times$ is an $m$-th power in $F_v^\times$, then $x$ is an $m$-th power in $\mathcal{O}_v^\times$. \ps\ps

Consider the subgroup $H_3 = H_2 \widetilde{K}^S \subset \widetilde{G}(\AAA_F)$. The inclusions \eqref{eq1phi3} and \eqref{eq2phi3} show that there is a unique function $\varphi_3 : H_3 \rightarrow \C$ such that $\varphi_3(hk) = \varphi_2(h)$ for all $h \in H_2$ and $k \in \widetilde{K}^S$. Note that $H_3$ is open in $\widetilde{G}(\AAA_F)$. Indeed, it contains the image of $\iota_v$ for each $v \in S$, which is open in $\widetilde{G}(F_v)$ by the first paragraph of the proof, as well as $\widetilde{K}^S$.  Define now $\psi : \widetilde{G}(\AAA_F) \rightarrow \C$ by $\psi_{|H_3} = \varphi_3$ and $\psi(h) = 0$ for $h \notin H_3$. This function 
has central character $c$, is $\widetilde{K}^S$-invariant on the right, $\widetilde{G}(F)$-invariant on the left, and satisfies $\rho(\psi)=f$. We claim that we have $\psi \in  \mathcal{A}_{\rm cusp}(\widetilde{G})_{\widetilde{c}}$. One first observes that the $K_\infty$-finiteness of $\varphi$ implies the one of $\psi_{|H_3}$, hence the one of $\psi$ (use that $\psi$ vanishes outside $H_3$ and that $K_\infty$ is a subgroup of $H_3$). As the normal subgroup $K_\infty$ of $\widetilde{K}_\infty$ is easily seen to be of finite index, this implies the $\widetilde{K}_\infty$-finiteness of $\psi$.  Let now $g \in \widetilde{G}(\AAA_F)$. Observe that the map $\psi_g : C_\infty \times G_\infty \rightarrow \C$, $(z,h) \mapsto \psi(gzh)$, is either identically $0$ or a constant multiple of $\psi_{g'}$ with $g' \in G(\AAA_F)$. For such a $g'$ we have $\psi_{g'}(z,h) =\widetilde{c}_\infty(z)\varphi(g'h)$. As $\iota_v$ is a smooth finite covering for each Archimedean $v$, this implies that $\psi$ is smooth, as so are $\varphi$ and $\widetilde{c}$. As $\varphi$ and $\widetilde{c}$ are slowly increasing, and as $C_\infty \times G_\infty$ has a finite index in $\widetilde{G}_\infty$, this also shows that $\psi$ is slowly increasing. We have proved that $\psi$ is an automorphic form. \ps\ps
It only remains to check that $\psi$ is a cuspform. Let $\widetilde{P}$ be a parabolic subgroup of $\widetilde{G}$ defined over $F$, and $N$ its unipotent radical. We have $N \subset G$, hence $N(\AAA_F) \subset H_3$, so if $g \notin H_3$ we have $\psi(N(\AAA_F)g) =0$. If $g =\gamma h z k$ with $\gamma \in \widetilde{G}(F)$, $h \in G(\AAA_F)$, $z \in \widetilde{Z}(\AAA_F)$ and $k \in \widetilde{K}^S$, then we have $$\int_{N(F) \backslash N(\AAA_F)} \psi(ng) dn = \widetilde{c}(z) \int_{N(F) \backslash N(\AAA_F)} \varphi (\gamma^{-1}n\gamma h) dn,$$
which vanishes as $\varphi$ is a cuspform and $\gamma^{-1} N \gamma$ is the unipotent radical of the $F$-parabolic subgroup of $\gamma^{-1} \widetilde{P} \gamma \cap G$ of $G$. 
\end{pf}

Let $H$ be a linear reductive group defined and quasi-split over $F$. Recall that a {\it Whittaker datum} for $H$ is a quadruple $D=(B,T,U,\chi)$ where $B$ is a Borel subgroup of $H$ defined over $F$, $T$ is a maximal torus of $B$ defined over $F$, $U$ is the unipotent radical of $B$, and $\chi : U(\AAA_F) \rightarrow \C^\times$ is a continuous character which is trivial on $U(F)$ and {\it nondegenerate} (see {\it e.g} \cite[p. 54]{KS}). When $H$ is split, it means {\it e.g.} that $\chi$ is nontrivial on the root subgroup $U_\alpha \subset U$ for each simple root $\alpha$ relative to $(B,T)$. Recall that an irreducible $H(\AAA_F)$-{\rm gk}-submodule $V$ of the space of cuspforms on $H(\AAA_F)$ is said {\it generic} if there exists a Whittaker datum $D=(B,T,U,\chi)$, such that for all $\varphi \in V-\{0\}$ the function on $H(\AAA_F)$ defined by $${\rm W}_\varphi^D(g)=\int_{N(F)\backslash N(\AAA_F)} \varphi(ng) \chi(n) dn$$ is not identically zero.  As the subspace of $\varphi \in V$ with ${\rm W}^{D}_\varphi =0$ is a ${\rm gk}$-submodule of $V$, this latter condition is equivalent to ask that we have ${\rm W}^{D}_\varphi \neq 0$ for {\it some} $\varphi \in V$, by irreducibility of $V$.\ps\ps

\begin{prop}\label{extcuspformgen} Keep the assumptions of Proposition \ref{extcuspform} and assume furthermore that $G$ is quasi-split over $F$. Let $\widetilde{V} \subset \mathcal{A}_{\rm cusp}(\widetilde{G})_{\widetilde{c}}$ be a nonzero irreducible $\widetilde{G}(\AAA_F)$-{\rm gk}-submodule. The following properties are equivalent :
\begin{itemize} \ps \ps
\item[(i)] $\widetilde{V}$ is generic, \ps\ps
\item[(ii)] $\rho_{\rm cusp}(\widetilde{V})$ contains a nonzero generic irreducible $G(\AAA_F)$-{\rm gk}-sub\-module. 
\end{itemize}
\end{prop}

\begin{pf} Assume first that (i) holds. Let  $\widetilde{D}=(\widetilde{B},\widetilde{T},U,\chi)$ be a Whittaker datum such that for all $\varphi \in \widetilde{V}-\{0\}$ we have ${\rm W}_{\varphi}^{\widetilde{D}} \neq 0$. Then $D=(\widetilde{B}\cap G, \widetilde{T}\cap G, U, \chi)$ is a Whittaker datum for $G$. Choose a nonzero $\varphi \in \widetilde{V}$. Up to replacing $\varphi$ by some $\widetilde{K}_\infty \times \widetilde{G}(\AAA_F^f)$-translate as in the proof of Proposition \ref{extcuspform} (ii), we may assume that we have ${\rm W}^{\widetilde{D}}_\varphi \neq 0$ on $G(\AAA_F)$. This implies that the cuspform $f := \rho(\varphi)=\varphi_{|G(\AAA_F)}$ is nonzero and ${\rm W}_{f}^D \neq 0$. The $G(\AAA_F)$-{\rm gk}-submodule of $\rho_{\rm cusp}(\widetilde{V}) \subset \mathcal{A}_{\rm cusp}(G)_{c}$ generated by $f$ is semisimple, hence a finite direct sum of irreducible {\rm gk}-modules $V_i$. Write $f = \sum_i f_i$ with $f_i \in V_i$. We have ${\rm W}_{f}^D = \sum_i {\rm W}_{f_i}^{D} \neq 0$, so there exists $i$ such that ${\rm W}_{f_i}^D \neq 0$, and $V_i$ is a generic irreducible constituent of $\rho_{\rm cusp}(\widetilde{V})$. \ps\ps
  Assume now that $\rho_{\rm cusp}(\widetilde{V})$ contains some generic irreducible $G(\AAA_F)$-submodule $V \neq 0$. Let $D=(B,T,U,\chi)$ be a Whittaker datum for $G$ such that ${\rm W}_\varphi^D$ is nonzero for all $\varphi \in V -\{0\}$. There exists a unique Whittaker datum $\widetilde{D}$ for $\widetilde{G}$ of the form $(\widetilde{B},\widetilde{T},U,\chi)$ with $\widetilde{B} \cap G = B$ and $\widetilde{T} \cap G = T$. As we have $0 \neq V \subset \rho(\widetilde{V})$, we may find some $\varphi \in \widetilde{V}$ such that the cuspform $f=\rho(\varphi)=\varphi_{|{\rm G}(\AAA_F)}$ is nonzero and belongs to $V$. It is enough to show ${\rm W}^{\widetilde{D}}_\varphi \neq 0$. But for $g \in G(\AAA_F)$ we have the identity ${\rm W}^{\widetilde{D}}_\varphi(g) = {\rm W}^{D}_{f}(g)$, and we are done. 
\end{pf}
  }

\end{document}